\documentclass[11pt]{article}
\usepackage{latexsym}
\usepackage{amsmath, amsfonts,amssymb, amsthm, euscript,makeidx,color,mathrsfs}
\usepackage{comment}

\newtheorem{assumption}{Assumption}
\def\qed{ \ \vrule width.2cm height.2cm depth0cm\smallskip}

\newcommand{\la}{\langle}
\newcommand{\ra}{\rangle}

\newcommand{\ba}{\begin{array}}
\newcommand{\ea}{\end{array}}
\newcommand{\be}{\begin{equation}}
\newcommand{\ee}{\end{equation}}
\newcommand{\bea}{\begin{eqnarray}}
\newcommand{\eea}{\end{eqnarray}}
\newcommand{\beaa}{\begin{eqnarray*}}
\newcommand{\eeaa}{\end{eqnarray*}}

\def\dbE{\mathbb{E}}
\def\dbF{\mathbb{F}}

\def\dbL{\mathbb{L}}

\def\dbP{\mathbb{P}}

\def\a{\alpha}
\def\b{\beta}

\def\d{\delta}
\def\e{\varepsilon}

\def\k{\kappa}
\def\l{\lambda}

\def\si{\sigma}
\def\t{\tau}
\def\f{\varphi}
\def\th{\theta}
\def\o{\omega}

\def\G{\Gamma}
\def\D{\Delta}
\def\Th{\Theta}
\def\L{\Lambda}

\def\O{\Omega}

\def\cF{{\cal F}}

\def\no{\noindent}

\def\ms{\medskip}

\def\q{\quad}
\def\qq{\qquad}

\def\pa{\partial}
\def\cd{\cdot}
\def\cds{\cdots}

\def\tr{\hbox{\rm tr}}

\def\qed{ \hfill \vrule width.25cm height.25cm depth0cm\smallskip}

\newcommand{\basa}{\begin{assumption}}
\newcommand{\easa}{\end{assumption}}

\newcommand{\bas}{\begin{assum}}
\newcommand{\eas}{\end{assum}}

\def\pa{\partial}

\def\cd{\cdot}
\def\cds{\cdots}

\def\tr{\hbox{\rm tr$\,$}}

\def\ol{\overline}

\def\dis{\displaystyle}

\def\cad{c\`{a}dl\`{a}g}

\def\1{{\bf 1}}

\def\:{\!:\!}
\def\reff#1{{\rm(\ref{#1})}}
\def \proof{{\noindent \bf Proof\quad}}
  
at 9pt
\newtheorem{thm}{Theorem}[section]
\newtheorem{lem}[thm]{Lemma}

\newtheorem{prop}[thm]{Proposition}
\newtheorem{rem}[thm]{Remark}
\newtheorem{eg}[thm]{Example}
\newtheorem{defn}[thm]{Definition}
\newtheorem{assum}[thm]{Assumption}

\begin{document}

\title{\textbf{A Martingale Approach for Fractional Brownian Motions and
Related Path Dependent PDEs}}
\author{ Frederi Viens\thanks{
\noindent Michigan State University, Department of Statistics and Probability, East Lansing, MI 48824. E-mail: viens@msu.edu. This author is supported in part by NSF DMS awards \#1407762 and \#1734183.} ~ and ~{Jianfeng Zhang}\thanks{%
\noindent Department of Mathematics, University of Southern California, Los
Angeles, CA 90089. E-mail: jianfenz@usc.edu. This author is supported in
part by NSF grant \#1413717. }}
\date{\today}
\maketitle

\begin{abstract}
In this paper we study dynamic backward problems, with the computation of conditional expectations as a sepcial objective, in a framework where the (forward) state process satisfies a Volterra type SDE, with fractional Brownian motion as a typical example. Such processes are neither Markov processes nor semimartingales, and most notably, they feature a certain time inconsistency which makes any direct application of Markovian ideas, such as flow properties, impossible without passing to a path-dependent framework. Our main result is a functional It\^{o} formula, extending the seminal work of Dupire \cite{Dupire} to our more general framework. In particular, unlike in \cite{Dupire} where one needs only to consider the stopped paths, here we need to concatenate the observed path up to the current time with a certain smooth observable curve derived from the distribution of the future paths. This new feature is due to the time inconsistency involved in this paper. We then derive the path dependent PDEs for the backward problems. Finally, an application to option pricing and hedging in a financial market with rough volatility is presented.    
\end{abstract}



\vfill\eject

\section{Introduction}

\label{sect-introduction} \setcounter{equation}{0}

\subsection{Background and heuristic description of the main ideas}

This paper introduces a new technique for analyzing functionals of
non-Markov processes using ideas which are borrowed from the Markovian case,
but necessarily require taking into account path dependence. 
In the Markovian case, consider the conditional expection
\begin{equation*}
Y _{t}=\mathbb{E}\left[ g\left( X_{T}\right) |\mathcal{F}_{t}^{X}\right],\q 0\le t\le T,
\end{equation*}%
where $g$ is a continuous function, $X$ is a Markov diffusion process, and $\{
\mathcal{F}_{t}^{X}\}_{0\le t\le T} $ is the filtration generated by $X$,  it is well known that $Y_{t}$ is a
deterministic function of $X_{t}$ only
\begin{equation*}
Y _{t}=u\left( t,X_{t}\right) ,
\end{equation*}
and that function $u$ solves a parabolic PDE, at least in weak form.

When $g(X_T)$ is replaced with a more general $\cF^X_T$-measurable random variable $\xi$ and/or $X$ is a non-Markov diffusion process, $Y_t$ will depend on the entire path of $X$ up to time $t$: 
\begin{equation}
Y_t = u(t,  X_{[0, t]}),\q \mbox{where}\q  X_{[0, t]} := \{X_r\}_{0\le r\le t},\label{Ypath}
\end{equation}
 and $u$ solves a so called Path Dependent PDE (PPDE, for short), which was first proposed by Peng \cite{Peng}. A powerful tool to study PPDEs is Dupire \cite{Dupire}'s functional It\^{o} calculus, see also Cont \& Fournie \cite{CF1, CF2, Cont}. A successful viscosity theory for (fully) nonlinear PPDEs was established by Ekren, Keller, Ren, Touzi, \& Zhang \cite{EKTZ, ETZ1, ETZ2, RTZ}. We also refer to Lukoyanov \cite{Lukoyanov},  Peng \& Wang \cite{PW}, Peng \& Song \cite{PS}, and Cosso \& Russo \cite{CR2} for some different approaches, and the book Zhang \cite{Zhang} for more references in this direction. We shall emphasize that, in all above works, the state process $X$ is a diffusion process, in particular, it is a semimartingale under all involved probabilities.

In this paper, we are interested in extending the above path-dependent analysis to more "heavily" non-Markov processes $X$, beyond the semimartingale framework. Typical examples of such non-Markov $X$ are Gaussian processes with memory properties, such as the fractional Brownian motion (fBm). 
When $X$ is fBm or similar processes, if there is a hope to replicate PDE-type ideas for representing $Y _{t}$ even in the state dependent case  $\xi=g(X_T)$, then any representation using a deterministic function $u$ will necessarily depend on  the entire path of $X$ up to $t$, namely in the form of \reff{Ypath}. 

How to figure out this dependence, and how to find the deterministic function $u$ in as explicit a way as possible (theoretically or numerically), is what this paper is about. To illustrate our main idea, in Section \ref{sect-fBM} we consider a special case that:
\begin{equation} 
X_t =\int_{0}^{t}K\left( t,r\right) dW_r,\label{XK}
\end{equation}
where $W$ is a Brownian motion and $K$ is a deterministic kernel, and  thus $X$ is a Volterra type Gaussian process. This is in particular the case for fBM. Our main idea is to introduce the following simple but crucial auxiliary process $\Th$ with two time variables:
\begin{equation}
\Theta _{s}^{t}:=\int_{0}^{t}K\left( s,r\right) dW_r.
\label{ThetaGauss}
\end{equation}
This process $\Th$ enjoys many nice properties, as we explain in details below.

We first note that, for any fixed $s$, the process $t\in [0,s]\mapsto \Th^t_s$ is a martingale. The existing theory of PPDEs relies heavily on the semimartingale theory, while $X$ is not a semimartingale. So by adding $\Th$ as our "state process", we will be be able to exploit its (semi)martingale property and thus recover the PPDE language.  We remark that,   for $s>t$, we have the orthogonal decomposition: $X_s = \Th^t_s + [X_s - \Th^t_s]$. 
This elementary property is a common computational tool in stochastic analysis, used in many studies regarding fBm and related
processes (see the textbook Nualart \cite[Chapter 5]{Nbook}). However, we believe our paper
is the first instance where this property is applied in the context of
PPDEs; the reason for this may be that the property is usually invoked to
exploit the independent part $X_s - \Th^t_s$ of the decomposition, and the martingale
property of $\Th$ is rarely exploited. 

Next, the process $X$ typically violates the standard flow property, which is another major obstacle for using PDEs and PPDEs. The introduction of the martingale component $\Th$ is the key for recovering the flow property, or say $(X, \Th)$ together will enjoy certain "Markov" property. More precisely, we shall rewrite \reff{Ypath} as
\begin{equation}
Y_t = u\big(t, X_{[0, t)} \otimes \Th^t_{[t, T]}\big),\label{Ypath2}
\end{equation}
show that this function $u$ satisfies a PPDE. While in the standard literature on PPDEs as mentioned earlier, $u$ depends only on the stopped path $X_{[0, t]}$, in our situation  $u$ will depend on the path $\Th^t_{[t, T]}$ as well. This is the major difference between Dupire's functional It\^o calculus and our extension.  We also note that, when $K(t,r) = K(r)$,  then $\Th^t_s = X_t$ for $t\le s\le T$, and thus \reff{Ypath2} reduces back to \reff{Ypath}.

Moreover, the introduction of $\Th$ is also crucial for numerical computation of the $u$ in \reff{Ypath2}. On one hand,  writing $u$ as the solution to a PPDE  enables us to  extend the existing numerical methods for standard PDEs to PPDEs naturally, which will be carried out in a separate project. On the other hand, we note that the function $u$ in \reff{Ypath2} is continuous under mild conditions, which is important for numerical purpose. However, $\Th$ is typically discontinuous in $X$, so if we write $Y$ as a function of $X$ only, the function $u$ in \reff{Ypath} could be discontinuous and thus its numerical methods would be less efficient. 

Finally, we discuss the tractability of the process $\Th$, which is important both for numerical purpose and for applications, and we shall discuss it more when we consider a financial application in the next subsection and Section \ref{sect-Appli}. First we note that $\Th^t_s$ is $\cF^{W}_t$-measurable, so mathematically all our analysis will have no measurability issues. However, in many applications people may not observe  $W$. Fortunately, for the models we will consider,  $\Th^t_s$ will be measurable to the observed information. In particular, when $X$ is an fBM, actually we have $\dbF^X=\dbF^W$ since one can represent $W$ as a function of $X$ through certain transform operator, see the textbook Nualart  \cite[Chapter 5]{Nbook}; also see Mocioalca \& Viens \cite{MoVi} for the case that $K$ is of convolution form, which covers the so-called Riemann-Liouville fBm.  Then  it is legitimate to use $\Th^t_s$ at time $t$, provided we observe $X_{[0, t]}$.

Having explained the ideas in details, in Section 3 we turn to the general framework where $X$  solves a Volterra SDE, see \reff{X} below. In this case the corresponding $\Th$ will be a semimartingale, and still shares all the nice properties discussed above.  Our main technical result is a functional It\^{o} formula for functions $u$   of the form \reff{Ypath2}, extending Dupire  \cite{Dupire}'s functional It\^{o} calculus which involves only the stopped process $X_{[0, t]}$. We remark that in Dupire's calculus the spatial derivatives involve only perturbations of $X_t$, but not of  the path before $t$: $X_{[0, t)}$. The main feature of our extension is again that the state variable contains the auxiliary process $\Th$, and our path derivatives will involve only the perturbation of $\Th^t_{[t, T]}$, which is exactly in the spirit of Dupire. This is important because on one hand $\Th^\cd_s$ is a semimartingale so an It\^{o} formula involving its derivatives is possible, and on the other hand $X_\cd$ is not a semimartingale so its derivatives are not helpful for It\^{o} calculus and should be avoided . 
%
 We note that Dupire's calculus serves as an alternative to the Malliavin
calculus, and appears as a simpler calculus of variations,
when questions of measurability with respect to current information are
crucial.  Said differently, the Malliavin calculus can be viewed as an
overkill from the standpoint of keeping track of this adaptability, since it
applies equally well to anticipating processes.

Section \ref{sect-PPDE} applies our functional It\^{o} formula to solve the backward problems in such framework and obtain naturally the PPDEs. We shall formulate it as a nonlinear Backward SDE (BSDE), whose linear version is essentially the conditional expectation $\dbE[g(X_\cd) |\cF_t]$ as discussed in Section \ref{sect-fBM}. Such nonlinear problems have many applications, especially in finance, stochastic control, and probabilistic numerical methods. 
We identify the corresponding semilinear PPDE, and assuming a classical solution exists, it yields immediately a solution to the BSDE. This strategy, known as a nonlinear Feynman-Kac formula, goes back to the original work of Peng \cite{Peng1}. 
 Section \ref{sect-PPDE} also provides a brief discussion of fully nonlinear PPDE, corresponding to stochastic optimization problems with diffusion controls in our framework. 

In Section \ref{sect-Appli}, we apply our methodology to the option pricing and hedging problem in a rough volatility model, motivated by the recent work El Euch \& Rosenbaum \cite{ER}. We discuss this and more general financial applications in the next subsection. 

\subsection{Application in rough volatility models }
\label{Fintro}

Consider a standard stock price model under risk neutral probability: 
\begin{equation}
dS_{t}= 
\sigma _{t}S_{t}dB_{t},  \label{BS}
\end{equation}
where $B$ is a Brownian motion, $\sigma $ is the volatility process, and we are assuming zero interest rate for simplicity.  
A number of recent studies have questioned the possibility
of assuming that $\sigma $ is a Markov process. In continuous time, the
first paper to work with this assumption is Comte \& Renault \cite{CR}, in which $\sigma $ is
assumed to be driven by an fBm. So-called continuous-time long-memory models
of that sort have been the main source of highly non-Markov volatility
model. The paper Chronopoulou \& Viens \cite{CV} can be consulted for references to such
works, and pertains to validating and calibrating this type of model from
option data.  Fractional stochastic volatility models continue to draw lots of 
interest.   A notable work is  Gatheral, Jaisson, \& Rosenbaum \cite{GJR}, which finds
market evidence that 
volatility's high-frequency behavior could be modeled as a rough path, e.g. based on fBm
with $H\in (0,1/2)$, and thus introduced the rough volatility models; see also  Bennedsen, Lunde, \& Pakkanen \cite{BLP}.  Among many others,   Abi Jaber, Larsson, \& Pulido \cite{ALP} and  Gatheral \& Keller-Ressel \cite{GK} studied affine variance models, where $\si^2$ is modeled as a convolution type linear Volterra SDE, which in particular includes the rough Heston model studied in El Euch \& Rosenbaum \cite{ER0, ER}; Bayer, Friz, \& Gatheral \cite{BFG} studied a rough Bergomi model; Cuchiero \& Teichmann \cite{CT} studied  affine Volterra processes with jumps; Gulisashvili, Viens, \& Zhang \cite{GVZ} provided an asymptotic analysis
applying to short-time fBm-modeling of volatility for fixed-income
securities near maturity; and  Fouque \& Hu  \cite{FH} studied an portfolio optimization problem in a model with fractional Ornstein-Uhlenbeck process. 

We remark that our model of general Volterra SDEs covers all the models mentioned above, except the jump model in \cite{CT} which we believe can be dealt with by extending our work to PPDEs on {\cad} paths, see Keller \cite{Keller} where the state process is a standard jump diffusion. Several works in the literature, e.g. \cite{BFG, ER, GK}, have already used the forward variance processes $\hat \Th^t_s:= \dbE[ \si^2_s |\cF_t]$, which is closely related to our process $\Th^t_s$. Indeed, in the affine variance models, they can be transformed from one to the other, as we will see in Section \ref{sect-Appli}.  However, for general models, especially when the drift of the Volterra SDE is nonlinear, we believe our process $\Th$ is intrinsic and is more convenient. Moreover, we note that most works in the literature either focus on modeling the financial market, or on pricing contingent claims in rough volatility models. We shall provide a systematic study on  the hedging  of contingent claims in general rough volatility models, motivated by the work \cite{ER}. Finally, we allow the backward process to be nonlinear, which appears naturally in applications, for example when the borrowing and lending  interest rates are different or when a control is involved; and we allow the payoff to be path dependent, thus including Asian options and lookback options.

We now focus on the hedging issue. Consider the model \reff{BS} and assume $\si^2$ is modeled through certain Volterra SDE which induce the crucial auxiliary process $\Th$. Then the price of the contingent claim will take the form $Y_t = u\big(t, S_t, \si^2_{[0, t)} \otimes_t \Th^t_{[t, T]}\big)$ and $u$ solves a PPDE.
We note that the market is typically incomplete when $S$ is the only tradable instruments. Applying our functional It\^{o} formula, we will see in Section \ref{sect-Appli} that the contingent claims can be hedged by using $S$ and $\Th$ (in appropriate sense),  and the hedging portfolios are exactly in the spirit of the $\D$-hedging, as derivatives of the function $u$  with respect to $S$ and $\Th$, respectively.  Then the issue boils down to whether or not we can hedge $\Th$ by using tradable assets in the market. Note that $\Th$ is defined through $\si^2$, so the key is to understand the variance $\si^2$ or the volatility $\si$. 

First note that, by observing $S$ continuously in time, mathematically we may compute $\si$ from it, or say $\si$ is $\dbF^S$-measurable. In fact, based on this fact, in the financial market there are 
proxies for $\sigma $, such as the VIX index from the Chicago Board of Options Exchange, which is a proxy for the
volatility on the S\&P500 index. This VIX index has become so mature that even skeptics
when it comes to volatility quotes will argue that if $S$ in (\ref{BS}) is a model for the S\&P500, then both $S$ and $\sigma $ are observable stochastic
processes. Consequently,  depending on the market and on the assumptions one is willing to make regarding volatility quotes, 
we may assume that  $\left(S,\sigma \right)$ is observed. This clearly has advantage compared to computing $\si$ from $S$ when numerical methods are concerned. 

Next, note that in the simple setting \reff{XK}-\reff{ThetaGauss} and assume for simplicity that $\si = X$,  we have $\Th^t_s = \dbE[\si_s |\cF_t]$. Mathematically, to compute $\Th$ from the observed information $\si$, one needs to first compute $W$ from $\si$ by using certain transfer operator as mentioned before, and then do the pathwise stochastic integration in \reff{ThetaGauss}. Numerically this will be very expensive. 
%
 However, in the financial context, the expression $\dbE[\si_s |\cF_t]$
above is also known as the forward volatility at time $t$ with horizon 
$s>t$, and that is a market observable. For the S\&P500 and many other
equities, it is directly quoted by use of implied volatility. 

Finally, for the rough volatility models we will consider in Section \ref{sect-Appli}, we will be able to replicate $\Th^t_s$ by using the forward variance $\dbE[\si^2_s |\cF_t]$ (similar to the forward volatility as we just discussed), which can be further replicated (approximately)  by using the variance swaps.  See more details in Section \ref{sect-Appli}. 
We can therefore conclude that we are able to hedge the contingent claims by using the tradable assets $S$ and the forward variance $\dbE[\si^2_s |\cF_t]$.

\section{The flow property of fBm}
\label{sect-fBM} 
\setcounter{equation}{0}

This section provides simple heuristics in the case of $X=$ fBm for easily
tractable examples.

\subsection{Martingale decomposition of fBm}

For simplicity in this section we restrict to one-dimensional processes
only. Let $B^{H}$ be a fBm with Hurst parameter $H\in \left( 0,1\right) $.
As we mentioned in the introduction, Chapter 5 in the textbook Nualart \cite{Nbook}
explains that there exists a Brownian motion (standard Wiener process) $W$
and an explicitly known deterministic kernel $K(t,s)>0$ such that 
\begin{equation}
B_{t}^{H}=\int_{0}^{t}K(t,r)dW_{r}\quad \mbox{and}\quad \mathbb{F}^{B^{H}}=%
\mathbb{F}^{W}=:\mathbb{F},  \label{BH}
\end{equation}%
where the notation $\mathbb{F}^{X}$ is the filtration of $X:~\mathbb{F}%
^{X}=\left\{ \mathcal{F}_{t}^{X}:t\geq 0\right\} $. The inclusion $\mathbb{F}%
^{B^{H}}\subset \mathbb{F}^{W}$ is immediate. The reverse inclusion comes
from the existence of a bijective transfer operator to express $W$ as a
Wiener integral with respect to $B^H$. We remark that, among others, one main feature of fBM is the violation of the standard flow property which can be viewed as certain time inconsistency: for $0\le t < s \le T$,
\bea
\label{timeinconsistency}
B^H_s \neq \tilde B^{t,H}_s,\q\mbox{where}\q \tilde B^{t,H}_s := B^H_t +  \int_{t}^{s}K(s,r)dW_{r}.
\eea

We are interested in backward problems. Let $\xi =g(B_{\cdot }^{H})\in \mathbb{L}^{2}(\mathcal{F}_{T})$, and denote 
\begin{equation}
Y_{t}:=\mathbb{E}[\xi |\mathcal{F}_{t}].  \label{Yt}
\end{equation}%
Clearly $Y$ is a martingale. Our goal is to characterize $Y$ from the PPDE
point of view. Due to \textrm{(\ref{BH})}, it is clear that 
\begin{equation}
\label{u12}
Y_{t}=u_{1}\big(t,B^{H}_{[0, t]}\big)=u_{2}\big(t,W_{[0, t]}\big),
\end{equation}%
for some measurable functions $u_{1},u_{2}$. Since $B^{H}$ is not a
semimartingale (when $H\neq {\frac{1}{2}}$), it is difficult to derive a PPDE
for $u_{1}$. On the other hand, since $W$ is a standard Brownian motion, formally $u_{2}$
should satisfy a path dependent heat equation. Indeed, this is true if $%
u_{2}$ is continuous in $W$ in the topology of uniform convergence. However,
provided $g$ is continuous in $B^{H}$, since $B^{H}$ is discontinuous in $W$
in pathwise sense, it is unlikely that $u_{2}$ will have desired pathwise regularity.

To get around of this, we will utilize the following simple but crucial
decomposition: 
\begin{equation}
B_{s}^{H}=\Theta
_{s}^{t}+I_{s}^{t}:=\int_{0}^{t}K(s,r)dW_{r}+\int_{t}^{s}K(s,r)dW_{r},\q 0\leq t\leq s\leq T,
\label{mart}
\end{equation}%
where 
\begin{equation}
\Theta _{s}^{t}~:=~ \int_{0}^{t}K(s,r)dW_{r} ~=~\mathbb{E}_{t}[B_{s}^{H}]  \label{Th} 
\end{equation}%
is $\mathcal{F}_{t}$ measurable, and $I_{s}^{t}$ is independent of $\mathcal{%
F}_{t}$. We note that when $s$ is fixed,  the process $t\in [0, s]\mapsto  \Theta _s^t$ is an $\mathbb{F}$-martingale.

\begin{rem}
\label{rem-rough} 
 (i) If we use the rough path norms instead of the
uniform convergence topology, then $B^{H}$ can be continuous in $W$.
However,  this requires a weaker regularity on the function $u_2$ in  \reff{u12}, which would induce serious difficulties for the functional It\^{o} formula \reff{FIto} below, and we do not see how to use the existing PPDE theory
to exploit the rough-path dependence of  $B^{H}$ on $W$.
Further exploration in this direction may be worthwhile, but is beyond the
scope of this paper.

(ii) Alternatively one may view the $u_2$ in \reff{u12} as a weak solution to the path dependent heat equation, in Sobolev sense without requiring pathwise regularity as in  Cont \cite{Cont}. However, the regularity itself is interesting and important, e.g. when one considers numerical methods. Moreover, in applications typically one observes $B^H$. Although in theory one may obtain $W$ from $B^H$ through the bijective transfer operator, such operator is not explicit and thus in practice it is not convenient to use the information $W$. As we will see soon, we shall express $Y$ through $\Th$ which is more trackable  in many applications.
\qed
\end{rem}

\subsection{A state dependent case}

\label{ss-state}To begin by illustrating our idea of using the process $%
t\mapsto \Theta _{s}^{t}$ in the simplest possible context, we consider a
very special case: 
\begin{equation*}
\xi =g(B_{T}^{H}).
\end{equation*}%
By the martingale/orthogonal decomposition (\ref{mart}) we have 
\begin{equation*}
Y_{t}=u\big(t,\Theta _{T}^{t}\big),\quad \mbox{where}\quad u(t,x):=\mathbb{E}%
\big[g\big(x+I_{T}^{t}\big)\big].
\end{equation*}%
Assuming $g$ is smooth, then the regularity of $u$ is clear. Moreover, since 
$t\mapsto \Theta^t_T$ is a martingale, applying the standard
It\^{o} formula we obtain: 
\begin{equation*}
dY_{t}=du(t,\Theta _{T}^{t})=[\partial _{t}u(t,\Theta _{T}^{t})+{\frac{1}{2}}%
\partial _{xx}^{2}u(t,\Theta _{T}^{t})K^{2}(T,t)]dt+\partial _{x}u(t,\Theta
_{T}^{t})K(T,t)dW_{t}.
\end{equation*}%
Noticing that $Y$ is a martingale by definition (\ref{Yt}), this implies
that the drift term above must vanish,  that is
\begin{equation}
\partial _{t}u(t,x)+{\frac{1}{2}}K^{2}(T,t)\partial _{xx}^{2}u(t,x)=0,\quad
u(T,x)=g(x).  \label{simpleheat}
\end{equation}%
This very simple (backward) heat equation with a time-dependent diffusion
coefficient shows how to compute the martingale $Y$ by tracking the observed
martingale process $t\mapsto \Theta _{T}^{t}$ and solving a deterministic
problem \reff{simpleheat} for $u$.

\begin{rem}
\label{rem-PDE} \textrm{If one were merely interested in finding a PDE
representation of }$Y$\textrm{\ at a given time, say }$Y_{0},$\textrm{\ a
natural solution, in analogy to the Brownian case, would emerge. Define 
\begin{equation*}
\tilde{u}(t,x):=\mathbb{E}[g(x+B_{T}^{H}-B_{t}^{H})]=\mathbb{E}%
[g(x+B_{T-t}^{H})],
\end{equation*}%
where the second equality holds because of the stationarity of increments of
fBm. Note that $x+B_{T-t}^{H}\sim Normal(x,(T-t)^{2H})$, then 
\begin{equation*}
\tilde{u}(t,x)=\int_{\mathbb{R}}g(y)p^{H}(T-t,y-x)dt,\quad \mbox{where}\quad
p^{H}(t,x):={\frac{1}{\sqrt{2\pi }t^{H}}}e^{-{\frac{x^{2}}{2t^{2H}}}}.
\end{equation*}%
One may check straightforwardly that 
$\partial _{t}p^{H}(t,x)=Ht^{2H-1}\partial _{xx}p^{H}(t,x).$
Then $\tilde{u}$ solves the following PDE 
\begin{equation*}
\partial _{t}\tilde{u}+H(T-t)^{2H-1}\partial _{xx}\tilde{u}(t,x)=0,\quad \tilde{u%
}(T,x)=g(x).
\end{equation*}%
This PDE  was already obtained by Decreusefond $\&$ Ustunel \cite{DU}, see also Baudoin $\&$ Coutin \cite{BC} for a more general result in this direction,   and it does not look very different from the previous one we identified.
Note that we still have $\tilde{u}(T,B_{T}^{H})=\xi $ and $\tilde{u}%
(0,B_{0}^{H})=Y_{0}$, however, $\tilde{u}(t,B_{t}^{H})$ is not a martingale
and in particular $\tilde{u}(t,B_{t}^{H})\neq Y_{t}$ for $0<t<T$. This is due to the
fact that the natural decomposition $%
B_{T}^{H}=B_{t}^{H}+(B_{T}^{H}-B_{t}^{H})$ is not an orthogonal
decomposition, namely $B_{t}^{H}$ and $B_{T}^{H}-B_{t}^{H}$ are not
independent.}

\textrm{The standard technique in the Brownian case }$W$\textrm{, for which
this decomposition over increments works so well, happens to coincide with
the use of }$\Theta $\textrm{\ since }$\Theta _{T}^{t}=W_{t}$\textrm{\ and }$%
I _{T}^{t}=W_{T}-W_{t}$\textrm{. But when trying the same increments
trick with $B_{t}^{H}$\ instead of }$W$\textrm{, since $\tilde{u}%
(t,B_{t}^{H})\neq Y_{t}$, the PDE above does not help track the value of
conditional expectations dynamically. Thus the orthogonal decomposition (\ref%
{mart}) is preferable for purposes, such as in stochastic finance, where }$%
t\mapsto Y_{t}$ \textrm{needs to be evaluated dynamically: we want to use a
single PDE (or PPDE) to represent all values of }$Y$\textrm{\ and for this,
we need to track }$\Theta $\textrm{.\hfill \vrule width.25cm height.25cm
depth0cm\smallskip }
\end{rem}

\subsection{A simple path dependent case}
\label{ss-path}

We now consider the case 
\begin{equation}
\xi =g(B^H_{T})+\int_{0}^{T}f(t,B^H_{t})dt.  \label{linear} 
\end{equation}%
As we alluded to in the introduction, this is a typical example useful in
finance, for instance as a model of a portfolio utility with legacy ($g$)
and consumption ($f$) terms, or as a contingent claim with straightforward
path dependence such as in Asian options. Because of the explicit path
dependence, we contend that our framework based on tracking $\Theta $ can
handle this dependence without any additional effort beyond what needs to be
deployed to handle the stochastic path dependence in $B^H$.

\begin{rem}
\label{rem-Markov}\textrm{In the special case $H={1/2}$, i.e. }$B^H$\textrm{\
is a standard Brownian motion, which we denote by }$W$\textrm{, we have 
\begin{equation*}
\tilde{Y}_{t}:=Y_{t}-\int_{0}^{t}f(s,W_{s})ds=\mathbb{E}\Big[%
g(W_{T})+\int_{t}^{T}f(s,W_{s})ds\big|\mathcal{F}_{t}\Big]=\tilde{u}(t,W_{t})
\end{equation*}%
where $\tilde{u}$ satisfies a standard backward heat equation with additive
forcing: 
\begin{equation*}
\qquad \qquad \qquad \qquad \qquad \partial _{t}\tilde{u}+{\frac{1}{2}}%
\partial _{xx}^{2}\tilde{u}+f(t,x)=0,\quad \tilde{u}(T,x)=g(x).\qquad \qquad
\qquad \qquad \quad \hfill \vrule width.25cmheight.25cmdepth0cm\smallskip
\end{equation*}%
}
\end{rem}

In the general fBm case, however, the above $\tilde{Y}_{t}$ is not Markovian
anymore. Instead, we use the same idea that leads to $\Theta $ to the
expression for $\xi $: we decompose the conditional expectation of the
integral part in $Y_{t}$ from (\ref{Yt}) into a term which is observable at
time $t$, and other terms. Those other terms are not independent of the
past, but we can apply the results of the previous section directly to them, whether at the terminal time, or for each $s$ in the integral from $t$
to $T$. Thanks to the explicit PDE (\ref{simpleheat}) from that section, we
obtain%
\begin{eqnarray*}
Y_{t} &=&\int_{0}^{t}f(s,B_{s}^{H})ds+\mathbb{E}[g(B_{T}^{H})|\mathcal{F}%
_{t}]+\int_{t}^{T}\mathbb{E}[f(s,B_{s}^{H})|\mathcal{F}_{t}]ds \\
&=&\int_{0}^{t}f(s,B_{s}^{H})ds+u_{g}(T;t,\Theta
_{T}^{t})+\int_{t}^{T}u_{f}(s;t,\Theta _{s}^{t})ds,
\end{eqnarray*}%
where 
\begin{equation}
\left. 
\begin{array}{c}
\partial _{t}u_{g}(T;t,x)+{\frac{1}{2}}K^{2}(T,t)\partial
_{xx}^{2}u_{g}(T;t,x)=0,~0\leq t\leq T,\quad u_{g}(T;T,x)=g(x); \\ 
\partial _{t}u_{f}(s;t,x)+{\frac{1}{2}}K^{2}(s,t)\partial
_{xx}^{2}u_{f}(s;t,x)=0,~0\leq t\leq s,\quad u_{f}(s;s,x)=f(s,x).%
\end{array}%
\right.  \label{uv}
\end{equation}%
Therefore, we discover that $Y_{t}$ can be expressed as a single
deterministic function $u\left( t,\cdot \right) $ of the concatenated path
which equals $B^{H}$ up to time $t$ and equals $\Theta ^{t}$ afterwards:
\begin{equation}
\left. 
\begin{array}{c}
\displaystyle Y_{t}=u\big(t,\{\int_{0}^{t\wedge s}K(s,r)dW_{r}\}_{0\leq s\leq
T})=u(t,B^{H}\otimes _{t}\Theta ^{t}\big),\quad \mbox{where} \\ 
\displaystyle(\omega \otimes _{t}\theta )_{s}:=\omega _{s}\mathbf{1}%
_{[0,t)}(s)+\theta _{s}\mathbf{1}_{[t,T]}(s), \\ 
\displaystyle u(t,\omega \otimes _{t}\theta ):=\int_{0}^{t}f(s,\omega
_{s})ds+u_{g}(T;t,\theta _{T})+\int_{t}^{T}u_{f}(s;t,\theta _{s})ds.%
\end{array}%
\right.  \label{Yu}
\end{equation}%
The last line above clearly shows how $u\left( t,\cdot \right) $ is an
explicit function of the entire concatenated path $\omega \otimes _{t}\theta 
$. From that standpoint, at least in this example, we have succeeded in
representing the conditional expectation $Y$ dynamically thanks to a single
deterministic path-dependent functional, by tracking $\Theta $. 
 We remark that the above function $u$ is continuous
in $(t,\omega \otimes _{t}\theta )$ under mild and natural conditions. This
is reassuring for our goal, which is to introduce appropriate time and path
derivatives and then derive a path dependent PDE for the above $u$ and in
more general contexts. Moreover, such regularity is  important when one considers numerical methods, even though it is not the focus of this paper.

Let us first take a look heuristically at the above example. Differentiating 
$u$ formally:
\begin{equation*}
\partial _{t}u(t,\omega \otimes _{t}\theta )=f(t,\omega _{t})+\partial
_{t}u_{g}(T;t,\theta _{T})-u_{f}(t;t,\theta _{t})+\int_{t}^{T}\partial
_{t}u_{f}(s;t,\theta _{s})ds.
\end{equation*}%
Note that 
\begin{equation*}
u_{f}(t;t,\theta _{t})=f(t,\theta _{t})=f(t,\omega _{t}),\quad %
\mbox{provided}~\theta _{t}=\omega _{t},
\end{equation*}%
resulting a corresponding cancellation in $\partial _{t}u$. This condition $%
\theta _{t}=\omega _{t}$ simply requires continuity of the path $\omega
\otimes _{t}\theta $ at the point $s=t$. This is certainly the case for us
in this section since $\Theta _{t}^{t}=B_{t}^{H}$, and will remain true for
any continuous Gaussian process and for non-Gaussian processes of interest,
as can be seen  in (\ref{Th2}) below. Then, by \textrm{(\ref{uv}%
)} we have 
\begin{eqnarray}
\partial _{t}u(t,\omega \otimes _{t}\theta ) &=&\partial
_{t}u_{g}(T;t,\theta _{T})+\int_{t}^{T}\partial _{t}u_{f}(s;t,\theta _{s})ds
\notag \\
&=&-{\frac{1}{2}}K^{2}(T,t)\partial _{xx}^{2}u_{g}(T;t,\theta _{T})-{\frac{1%
}{2}}\int_{t}^{T}K^{2}(s,t)\partial _{xx}^{2}u_{f}(s;t,\theta _{s})ds  \label{PPDEsimple}
\\
&=&-{\frac{1}{2}}K^{2}(T,t)\partial _{\theta _{T}\theta _{T}}^{2}u(t,\omega
\otimes _{t}\theta )-{\frac{1}{2}}\int_{t}^{T}K^{2}(s,t)\partial _{\theta
_{s}\theta _{s}}^{2}u(t,\omega \otimes _{t}\theta )ds, \notag 
\end{eqnarray}%
with the terminal condition 
$u\left( T,\omega \otimes _{T}\theta \right) =g\left( \omega _{T}\right) +\int_{0}^{T}f\left( t,\omega _{t}\right) dt.
$
The last equality in (\ref{PPDEsimple}) comes from the expression for $u$ in
(\ref{Yu}). In that line, the notation $\partial _{\theta _{s}\theta _{s}}^{2}u$ means that this is a derivative with respect to the value of the path $\omega \otimes _{t}\theta $ at time $s\in [t,T]$.

In any case, at least heuristically, one sees that $u$ itself appears to
solve a PPDE in some sense.\ We will make this sense precise in the next
section, in Theorem \ref{thm-rep}. One observation is that the spatial
derivatives of $u$ above are only with respect to $\theta $, not  to $\omega $,
which is crucial in our context because $\omega $ corresponds to $B^{H}$
which is not a semimartingale, and we wish to use the semimartingale
property to transport stochastic objects (such as conditional expectations)
to their deterministic representations.

\begin{rem}
\label{rem-flow}
(i) The kernel $K$ involves two time variables, and thus $B^H$ is not a semimartingale (when $H\neq {1\over 2}$). Moreover, $B^H$ violates the standard flow property, see \reff{timeinconsistency}, and is by nature time inconsistent.  The introduction of the term $\Th$ is the key for recovering the flow property for $B^H$, which is crucial for deriving the corresponding PPDE. 

(ii) In the standard functional It\^{o} calculus of Dupire \cite{Dupire}, the function $u$ depends only on the stopped path $X_{t\wedge \cd}$, or equivalently, in their situation $\Th^t$ is flat. So our main result in the next section is indeed a nontrivial extension of Dupire's result. 

 (iii) We also note that the occurrence of
derivatives of $u$ with respect only to the \textquotedblleft $\theta $%
\textquotedblright\ portion of the path does not contradict Dupire's
functional It\^{o} calculus, because, as just mentioned,  in his setting  $\theta _{s}=\omega
_{t}$ for all $s\in \lbrack t,T]$, and thus the derivatives with respect to $%
\theta $ reduce to the derivative with respect to $\omega _{t}$ alone.
\qed
\end{rem}

\section{Functional It\^{o} formula}
\label{sect-FIO} \setcounter{equation}{0}

In this section, we expand our framework's reach by considering more general
processes $X$, beyond the Gaussian class. Assume $X$ is a solution to the
$d$-dimensional  Volterra SDE:  
\begin{equation}
X_{t}=x+\int_{0}^{t}b(t;r,X_{\cdot })dr+\int_{0}^{t}\sigma (t;r,X_{\cdot
})dW_{r},\quad 0\leq t\leq T,  \label{X}
\end{equation}%
where $W$ is a standard (possibly multidimensional) Wiener process, and $b$ and $\sigma $ have appropriate dimensions and are adapted in
the sense that $\varphi (t; r,X_{\cdot })=\varphi (t; r,X_{r\wedge \cdot })$
for $\varphi =b,\sigma $.  As in \reff{timeinconsistency}, one main feature of such SDE is that it violates the flow property. That is,  $X_s \neq \tilde X^t_s$ for $0\le t<s\le T$, where $\tilde X^t$ is the solution to the following SDE:
\beaa
\tilde X^t_s= X_t+\int_t^s b(s;r, X\otimes_t \tilde X^t_{\cdot })dr+\int_t^s \sigma (s;r,X\otimes_t \tilde X^t)dW_{r},\quad t\leq s\leq T.
\eeaa

Throughout the paper, the following assumption will always be in force. 
\begin{assum}
\label{assum-X}
(i) The  SDE \textrm{(\ref{X})} admits a weak solution $(X, W)$.

(ii)  $\mathbb{E}\big[\sup_{0\le t\le T}|X_t|^p\big] <\infty$ for all $p\ge 1$.
\end{assum}
The condition (ii) is technical, and in order not to distract our main focus, we postpone its discussion to Appendix. For condition (i), there have been many works on wellposedness of Volterra SDEs, see e.g. Berger \& Mizel \cite{BM1, BM2}. In this paper we prefer not to restrict to specific conditions so as to allow for the most generality, and in applications any reasonable model should  admit at least one solution. However, we would like to mention that, since in most applications $X$ is the  observable state process and $W$ is just used to model the distribution of $X$, it suffices to consider a weak solution. Moreover, no uniqueness of weak solution is needed. So from now on, we will always fix a weak solution $(X, W)$, and slightly unlike in Section \ref{sect-fBM}, we shall always use the full filtration:
\bea
\label{F}
\dbF = \dbF^{X, W}.
\eea
In this framework, the analogue of the martingale term in the
decomposition of $X$ can be defined using exactly the same idea as in the
Gaussian case, by basing it on (\ref{X}) rather than (\ref{BH}). Thus we
denote 
\begin{equation}
\Theta _{s}^{t}:=x+\int_{0}^{t}b(s;r,X_{\cdot })ds+\int_{0}^{t}\sigma
(s;r,X_{\cdot })dW_{r},\quad t\leq s\leq T, \label{Th2}
\end{equation}
where $t\mapsto \Th^t_s$ is a semimartinagle. 
To simplify the notation, quite often we omit the $X$ in $b$ and $\sigma $,
and simply write $\varphi (t,s)=\varphi (t;s,X_{\cdot })$ for $\varphi
=b,\sigma $.

\begin{rem}
\label{rem-regularity}
As we will see in the paper, we shall write the interested value process $Y_t$ as a function of the paths $X\otimes_t \Th^t$, which we observe, and the function $u$ is typically continuous under mild conditions.  We note that $\Th^t_s$ is also a function of $X_{[0, t]}$. However, this dependence is typically discontinuous under uniform convergence. For example, set 
\beaa
b=0;\q \si(t; s, \o) = 1,~ t\in [0,  {T\over 2}];\q \si(t; s, \o) = 1+[t-{T\over 2}] \si_0(s, \o),~ t\in ({T\over 2}, T],
\eeaa
 for some appropriate function $\si_0$. Then $X_{[0, {T\over 2}]} = W_{[0, {T\over 2}]}$, and, for $0<t\le {T\over 2}$,  $\Th^t_T = W_t+ {T\over 2} \int_0^t \si_0(s, W_\cd) dW_s$. This involves a stochastic integral and is typically discontinuous in pathwise sense.  Consequently, if we rewrite $Y_t$ as a function of $X_{[0, t]}$ only, the function could be discontinuous. Besides theoretical interest, such regularity is crucial when one studies numerical methods for the related problems.
\qed
\end{rem}

\subsection{The path derivatives\label{Derivatives}}

As in Dupire \cite{Dupire}, though in the end all paths are continuous,
since we employ some piecewise-continuous approximations, we must extend the
sample space to the {c\`{a}dl\`{a}g} space $D^{0}$. Denote 
\begin{equation*}
\left. 
\begin{array}{c}
\displaystyle\Omega :=C^{0}([0,T],\mathbb{R}^{d}),\quad \overline{\Omega }%
:=D^{0}([0,T],\mathbb{R}^{d}),\quad \Omega _{t}:=C^{0}([t,T],\mathbb{R}%
^{d});\medskip \\ 
\displaystyle\Lambda :=[0,T]\times \Omega ,\quad \overline{\Lambda }:=\Big\{%
(t,\omega )\in \lbrack 0,T]\times \overline{\Omega }:\omega |_{[t,T]}\in
\Omega _{t}\Big\};\medskip \\ 
\displaystyle\Vert \omega \Vert _{T}:=\sup_{0\leq t\leq T}|\omega
_{t}|,\quad \mathbf{d}((t,\omega ),(t^{\prime },\omega ^{\prime
})):=|t-t^{\prime }|+\Vert \omega -\omega ^{\prime }\Vert _{T}.%
\end{array}%
\right.
\end{equation*}%
Here we change our notation slightly compared to what we had used in the
illustrative examples of the previous section. Indeed, we see here that $%
\omega $ is defined on $[0,T]$, whereas previously we used the letter $%
\omega $ for paths on $[0,t)$. The correspondence between these two
conventions is that what we now call $\omega \mathbf{1}_{[0,t)}$ and $\omega 
\mathbf{1}_{[t,T]}$ correspond to the $\omega $ and $\theta $ in \textrm{(%
\ref{Yu})}, respectively. Though the old covention was natural because it
highlighted the concatenation of the path of $X$ up to $t$ with the path of
its observable martingale component $\Theta ^{t}$ after $t$, the new
convention does not presume that this is the structure of the full path $%
\omega $, and allows more compact notation. Moreover, we emphasize that in this subsection all the terms are deterministic and there is no probability involved.

The space we really care about is $\Lambda $, however, for technical reasons
we need to allow $\omega $ to be discontinuous on $[0,t]$. We note that
Dupire's framework is covered by our setup since Dupire's space is the
subset of those $\omega $ which are constant on $[t,T]$. We also note that,
for any $(t,\omega )\in \overline{\Lambda }$ and $t^{\prime }>t$, we have $%
(t^{\prime },\omega )\in \overline{\Lambda }$. Let $C^{0}(\overline{\Lambda }%
)$ denote the set of functions $u:\overline{\Lambda }\rightarrow \mathbb{R}$
continuous under $\mathbf{d}$. For $u\in C^{0}(\overline{\Lambda })$, define 
\begin{equation}
\partial _{t}u(t,\omega ):=\lim_{\delta \downarrow 0}{\frac{u(t+\delta
,\omega )-u(t,\omega )}{\delta }}\quad \mbox{for all}~(t,\omega )\in 
\overline{\Lambda },  \label{patu}
\end{equation}%
provided the limit exists. We note that here $\partial _{t}u$ is actually
the right time derivative.

We next define the spatial derivative with respect to $\omega $. Given $%
(t,\omega )\in \overline{\Lambda }$, we define $\partial _{\omega }u(t,\omega )$ as  the Fr\'{e}chet derivative with respect to $\omega 
\mathbf{1}_{[t,T]}$ which is a linear operator on $\Omega _{t}$:  
\begin{equation}
u(t,\omega +\eta \mathbf{1}_{[t,T]})-u(t,\omega ) = \langle \partial _{\omega }u(t,\omega ),\eta \rangle + o(\|\eta\mathbf{1}_{[t,T]} \|_T),\q \mbox{for any}~\eta \in \Omega _{t}. \label{pathu}
\end{equation}%
It is clear that this is equal to the Gateux derivative:
\begin{equation}
\langle \partial _{\omega }u(t,\omega ),\eta \rangle =\lim_{\varepsilon
\rightarrow 0}{\frac{u(t,\omega +\varepsilon \eta \mathbf{1}%
_{[t,T]})-u(t,\omega )}{\varepsilon }},\quad \mbox{for any}~\eta \in \Omega
_{t}.  \label{pathu2}
\end{equation}
We emphasize that the above perturbation is only on $[t,T]$, not on $[0,t)$.
This is consistent with Dupire's derivative.  For any $s<t$ and $\eta \in
\Omega _{s}$, we will take the convention that 
\begin{equation}
\langle \partial _{\omega }u(t,\omega ),\eta \rangle :=\langle \partial
_{\omega }u(t,\omega ),\eta  \mathbf{1}_{[t,T]}\rangle. \label{pathu0}
\end{equation}%

\begin{defn}
\label{defn-C0} Let $u\in C^0(\overline \L)$ such that $\pa_\o u$ exists for all $(t,\o) \in \overline\L$.

(i) We say $\pa_\o u$ has polynomial growth if there exist constants $C>0$, $\kappa >0$ such that
\begin{equation}
\Big|\langle \pa_\o u(t,\o), \eta \rangle \Big| \le C [1+\|\o\|_T^\kappa]~ \|\eta \mathbf{1}_{[t,T]} \|_T,\q\mbox{for all}~ (t,\o)\in \overline\L, \eta \in \O. \label{paoubound}
\end{equation} 

(ii) We say $\pa_\o u$ is continuous if, for any $\eta \in \Omega $,  the mapping $(t,\o) \in \overline \L \mapsto
\langle \partial _{\omega }u(t,\o),\eta \rangle$ is continuous under $\mathbf{d}$.
\end{defn}
Throughout the paper, we use $\kappa$ to denote a generic order of polynomial growth, which may vary from line to line. We note that, when $\partial _{\omega }u$ is continuous, it is clear that the mapping $\l\in [0,1]\mapsto u(t,\omega +\l\eta \mathbf{1}_{[t,T]})$ is continuously differentiable, and thus 
\begin{equation}
u(t,\omega +\eta \mathbf{1}_{[t,T]})-u(t,\omega ) = \int_0^1 \langle \partial _{\omega }u(t,\omega +\l\eta \mathbf{1}_{[t,T]}), \eta \rangle d\l,\q \mbox{for any}~\eta \in \Omega _{t}.
 \label{pathu3}
\end{equation}
Define further the second derivative $\partial
_{\omega \omega }^{2}u(t,\omega )$ as a bilinear operator on $\Omega
_{t}\times \Omega _{t}$: 
\begin{equation}
\langle \partial
_{\omega }u(t,\omega + \eta _{1}\mathbf{1}_{[t,T]}),\eta
_{2}\rangle -\langle \partial _{\omega }u(t,\omega ),\eta _{2}\rangle = 
\langle \partial _{\omega \omega }^{2}u(t,\omega ),(\eta _{1},\eta
_{2})\rangle + o(\|\eta_1 \mathbf{1}_{[t,T]} \|_T),
\label{path2u}
\end{equation}%
for any $\eta _{1},\eta _{2}\in \Omega _{t}$. 
Similarly define $\langle \partial _{\omega \omega }^{2}u(t,\omega ),(\eta
_{1},\eta _{2})\rangle$  for $\eta _{1},\eta _{2}\in \Omega_s$ as in \reff{pathu0},
and define the polynomial growth and continuity  in the spirit of Definition \ref{defn-C0}, with  $\|\eta \mathbf{1}_{[t,T]} \|_T$ in \reff{paoubound} replaced with $\|\eta_1 \mathbf{1}_{[t,T]} \|_T\|\eta_2 \mathbf{1}_{[t,T]} \|_T$.

\begin{defn}
\label{defn-C12bar} We say $u\in C^{1,2}(\overline{\Lambda })\subset
C^{0}(\overline{\Lambda })$ if $\partial _{t}u$, $\partial _{\omega
}u,\partial _{\omega \omega }^{2}u$ exist and are continuous on $\overline{%
\Lambda }$. Let $C_{+}^{1,2}(\overline{\Lambda })$ be the subset of $C^{1,2}(%
\overline{\Lambda })$ such that all the derivatives have polynomial growth, and $\langle\pa_{\o\o}^2 u, (\eta, \eta)\rangle$ is  locally uniformly continuous in $\o$ with polynomial growth, that is, there exist $\kappa>0$ and a bounded modulus of continuity function $\rho$ such that, for any $(t, \o), (t,\o') \in \overline \L$ and $\eta \in \O_t$,
\begin{equation}
\Big|\big\langle \partial^2_{\omega\omega}u(t,\omega ) -\partial^2_{\omega\omega}u(t,\omega ^{\prime }), (\eta, \eta)\big\rangle \Big|\leq
[1+\|\o\|_T^\kappa +\|\o'\|_T^\kappa  ] ~\|\eta \mathbf{1}_{[t,T]} \|_T^2  ~\rho(\Vert \omega -\omega ^{\prime }\Vert _{T}). \label{paooucont}
\end{equation}
\end{defn}

\begin{rem}
\label{rem-C12bar}
Cont  {\rm \&} Fournier \cite{CF2} established the functional It\^{o} formula in their framework for  all $u \in C^{1,2}(\overline \L)$, by using the standard localization techniques with stopping times. In their framework only $(t, \o\mathbf{1}_{[0,t]})$ is involved and thus it is sufficient to consider the stopped paths $\o_{t\wedge \cd}$. However, in our framework, the whole path of $\o$ on $[0, T]$ is involved, we have difficulty to apply the localization techniques directly. Thus in this paper we require slightly stronger conditions by restricting  $u$ to $C_{+}^{1,2}(\overline{\Lambda })$. We shall leave the possible relaxation of these conditions in future research.
\qed
\end{rem}

\begin{eg}
\label{eg-path} The first example below is in the framework of the
path-dependent case of Section \ref{ss-path}. The second covers Dupire's
case.

(i) If $u(t,\omega )=g(\omega _{T})+\int_{t}^{T}f(s,\omega _{s})ds$ and $%
f,g$ are smooth, then 
\begin{eqnarray*}
&\partial _{t}u(t,\omega )=-f(t,\omega _{t}),\quad \langle \partial _{\omega
}u(t,\omega ),\eta \rangle =\pa_x g(\omega _{T})\cdot \eta
_{T}+\int_{t}^{T}\pa_x f(s,\omega _{s})\cdot \eta _{s}ds,& \\
&\big\langle \partial _{\omega \omega }^{2}u(t,\omega),(\eta ^{1},\eta ^{2})\big\rangle
=\pa^2_{xx}g(\omega _{T}):[\eta _{T}^{1}(\eta _{T}^{2})^{\top
}]+\int_{t}^{T}\pa^2_{xx}f(s,\omega _{s}):[\eta _{s}^{1}(\eta _{s}^{2})^{\top }]ds.&
\end{eqnarray*}
Here $A_1 : A_2 := \tr(A_1 A_2^\top)$ for two matrices $A_1, A_2$. 

(ii) If $u$ is adapted in the sense that, after time $t$, the path of $%
\omega $ is frozen, i.e. $u(t,\omega )=v(t,\omega \mathbf{1}_{[0,t)}+\omega
_{t}\mathbf{1}_{[t,T]})$ for some function $v$, then 
\begin{eqnarray*}
&\partial _{t}u(t,\omega )=\partial _{t}v(t,\omega ),~\langle \partial
_{\omega }u(t,\omega ),\eta \rangle =\partial _{\omega }v(t,\omega )\cdot
\eta _{t},& \\
&\langle \partial _{\omega \omega }^{2}u(t,\omega ),(\eta ^{1},\eta
^{2})\rangle =\partial _{\omega \omega }^{2}v(t,\omega ):[\eta _{t}^{1}(\eta
_{t}^{2})^{\top }],&
\end{eqnarray*}%
where $\partial _{t}v,\partial _{\omega }v,\partial _{\omega \omega }^{2}v$
are Dupire's path derivatives.
\end{eg}

The following result is similar to Cont \& Fourni\'{e} \cite{CF1}.

\begin{prop}
\label{prop-path} Let $u_1, u_2 \in C_{+}^{1,2}(\overline\Lambda)$. Assume $%
u_1=u_2$ on $\Lambda$, then $\partial_t u_1 = \partial_t u_2$, $\langle
\partial_\o u_1, \eta\rangle = \langle \partial_\o u_2, \eta\rangle$, $%
\langle \partial^2_{\omega\omega} u_1, (\eta, \eta)\rangle =
\langle\partial^2_{\omega\omega} u_2, (\eta, \eta)\rangle$ on $\Lambda$, for
all $\eta\in \Omega$.
\end{prop}

The proof is closely related to the functional It\^{o} formula below, so we
postpone it and combine with the proof of Theorem \ref{thm-FIto1}. We
believe it is possible to show that $\langle \partial _{\omega \omega
}^{2}u_{1},(\eta _{1},\eta _{2})\rangle =\langle \partial _{\omega \omega
}^{2}u_{2},(\eta _{1},\eta _{2})\rangle $ for all $\eta _{1},\eta _{2}\in
\Omega $ under possibly weaker regularity conditions on $u_1, u_2$.
 We do not pursue such generality in this paper.
We now define: 
\begin{defn}
\label{defn-C12} Let $C_{+}^{1,2}(\Lambda )$ denote the collection of
functions $u:\Lambda \rightarrow \mathbb{R}$ such that there exists $\tilde{u%
}\in C_{+}^{1,2}(\overline{\Lambda })$ satisfying $\tilde{u}=u$ on $\Lambda 
$. In this case we define the path derivatives: $\partial _{t}u :=\partial _{t}%
\tilde{u}$, $\partial _{\omega }u :=\partial _{\omega }\tilde{u}$, and $%
\partial _{\omega \omega }^{2}u :=\partial _{\omega \omega }^{2}\tilde{u}$ on $%
\Lambda $.
\end{defn}

By Proposition \ref{prop-path}, for any $\eta \in \Omega $, clearly $%
\partial _{t}u$, $\partial _{\omega }u$, and $\langle \partial _{\omega
\omega }^{2}u,(\eta ,\eta )\rangle $ are uniquely determined on $\Lambda $,
regardless of the choice of $\tilde{u}$ in Definition \ref{defn-C12}.

\subsection{Functional It\^{o} formula in the regular case\label{regular}}

As noted in the introduction, the use of fBm causes two difficulties: (i) it
is non-Markovian, because of the two-variable kernel $K(t,r)$; (ii) for $%
H<1/2$, the kernel is singular, i.e. $\lim_{r\rightarrow t}K(t,r)=\infty $.
The SDE \textrm{(\ref{X})} has the same issues, stemming from the same
properties of $b$ and $\sigma $ as functions of $\left( t,r\right) $. The
possible dependence of $b$ and $\sigma $ on the path $X$ (which we mostly
omit in the notation below) add to the path dependence. To understand the
problem better, in this subsection we focus on the lack of a Markov property
and additional path dependence, and postpone the
singularity issue to the next subsection; thus we assume $b$ and $\sigma $ have no singularity as $r$
tends to $t$. We remark that in this \textquotedblright
regular\textquotedblright\ case, $X$ is typically a semimartingale:
\begin{equation}
\label{Xsemimg}
dX_{t}=b(t;t)dt+\sigma (t; t)dW_{t}+\Big[\int_{0}^{t}\partial
_{t}b(t; r)dr+\int_{0}^{t}\partial _{t}\sigma (t; r)dW_{r}\Big]dt,
\end{equation}%
provided that $\partial _{t}b(t; r)$, $\partial _{t}\sigma (t; r)$ exist and
have good integrability properties near the time diagonal. We shall remark that,  even for  $H>{1\over 2}$, the fBM $B^H$ does not satisfy the above expression because the corresponding $\pa_t \si(t; r)$ is not square integrable. 

In this subsection we assume the following.
\begin{assum}
\label{assum-regular} $\partial _{t}b(t; r,\cdot ),\partial _{t}\sigma
(t; r,\cdot )$ exist for $t\in \lbrack r,T]$, and for $\varphi =b,\sigma
,\partial _{t}b,\partial _{t}\sigma $, 
\begin{equation}
|\varphi (t; r,\omega )|\leq C_{0}[1+\Vert \omega \Vert _{T}^{\kappa_{0}}]\quad 
\mbox{for some
constants}~C_{0}, \kappa_{0}>0. \label{growth}
\end{equation}
\end{assum}
\no Recall \reff{Th2} and \reff{Xsemimg}.  Under Assumptions \ref{assum-X} and \ref{assum-regular}, it is obvious that
\begin{equation}
\dbE\Big[\|X\otimes_t \Th^t\|_T^p\Big] \le C_p,\q\mbox{for all}~ p\ge 1, 0\le t\le T. \label{ETh}
\end{equation}

Our main result is the following  functional It\^{o} formula.

\begin{thm}
\label{thm-FIto1} Let Assumptions \ref{assum-X} and \ref{assum-regular} hold
and $u\in C_{+}^{1,2}(\Lambda )$. Then 
\begin{equation}
\left.\ba{c}
\dis du(t,X\otimes _{t}\Theta ^{t}) =\partial _{t}u(t,X\otimes _{t}\Theta
^{t})dt+{\frac{1}{2}}\langle \partial _{\omega \omega }^{2}u(t,X\otimes
_{t}\Theta ^{t}),(\sigma ^{t,X},\sigma ^{t,X})\rangle dt  \label{FIto} \\
\dis +\langle \partial _{\omega }u(t,X\otimes _{t}\Theta ^{t}),b^{t,X}\rangle
dt+\langle \partial _{\omega }u(t,X\otimes _{t}\Theta ^{t}),\sigma
^{t,X}\rangle dW_{t},\quad \mathbb{P}\mbox{-a.s.} 
\ea\right.
\end{equation}
where, for $\varphi =b,\sigma $,  $\varphi _{s}^{t,\omega }:=\varphi (s; t,\omega )$ emphasizes the
dependence on $s\in \lbrack t,T]$.
\end{thm}

The main idea of the proof follows that of Dupire \cite{Dupire} and Cont \& Fournie \cite{CF2}. However, here we have to deal with the 
  two time variables, and for that purpose we need a few technical lemmas. The first one  
is a direct consequence of the proof of Kolmogorov's continuity criterion, see e.g. Revuz \& Yor \cite[Chapter I, Theorem 2.1]{RY}.

\begin{lem}
\label{lem-kolmogorov} Let $\tilde X$ be a process on $[0, T]$ with $\tilde X_0=0$, and $\alpha, \beta>0$ be  constants. Assume, 
\begin{eqnarray*}
\mathbb{E}\Big[|\tilde X_t - \tilde X_{t^{\prime}}|^{2p}\Big] \le C_p
\beta^p |t-t^{\prime}|^{\alpha p},\quad \mbox{for all}~~ 0\le
t<t^{\prime}\le T,\quad p\ge 1.
\end{eqnarray*}
Then, for each $p\ge 1$,  there exists another constant $\tilde C_p>0$, which may depend on $T$, $\alpha$, $p$, 
the dimension $d$,  and the above $C_p$, but does not depend on $%
\beta$, such that 
\begin{eqnarray*}
\mathbb{E}\Big[ \|\tilde X\|_T^{2p}\Big] \le \tilde C_p \beta^p.
\end{eqnarray*}
\end{lem}

The following result will be crucial for the functional It\^{o} formula.

\begin{lem}
\label{lem-regXest} Let Assumptions \ref{assum-X} and \ref{assum-regular}
hold. Fix $n$ and set 
\begin{equation}
X_{t}^{n}:=\sum_{i=0}^{2^n-1}\Theta _{t}^{t_{i}}\mathbf{1}%
_{[t_{i},t_{i+1})}(t)+X_{T}\mathbf{1}_{\{T\}}(t),~ \mbox{where}~
h:=2^{-n}T,t_{i}:=ih,i=0,\cdots ,2^{n}.  \label{Xn}
\end{equation}%
Then, for any $p\geq 1$, 
\begin{equation}
\mathbb{E}\Big[\Vert X\otimes _{t}\Theta ^{t}-X\otimes _{t^{\prime }}\Theta
^{t^{\prime }}\Vert _{T}^{4p}\Big]\leq C_p|t^{\prime }-t|^{p},\quad \mathbb{E}%
[\Vert X-X^{n}\Vert _{T}^{8}]\leq C2^{-n}.  \label{regXest}
\end{equation}
\end{lem}

{\noindent \textbf{Proof\quad }} We start with the first inequality of \textrm{(\ref{regXest})}.  Assume $t<t^{\prime }$ and denote 
\begin{equation*}
\tilde{X}_{s}:=[X\otimes _{t^{\prime }}\Theta
^{t^{\prime }} - X\otimes _{t}\Theta ^{t}]_{s}=\int_{s\wedge t}^{s\wedge t^{\prime
}}b(s;r)dr+\int_{s\wedge t}^{s\wedge t^{\prime }}\sigma (s;r)dW_{r}.
\end{equation*}%
We claim that, for any $s<s^{\prime }$ and any $p\geq 1$, 
\begin{equation}
I_{p}:=\mathbb{E}\Big[|\tilde{X}_{s}-\tilde{X}_{s^{\prime }}|^{2p}\Big]\leq
C_{p}(t^{\prime }-t)^{\frac{p}{2}}(s^{\prime }-s)^{\frac{p}{2}}.
\label{DXest}
\end{equation}%
Then, the first inequality of \textrm{(\ref{regXest})} follows from Lemma %
\ref{lem-kolmogorov}. We prove \textrm{(\ref{DXest})} in three cases.

\textit{Case 1.} $s\le t$. Then $\tilde X_s=0$ and thus, by Assumptions \ref{assum-X} and \ref{assum-regular},
\begin{eqnarray*}
I_p &=& \mathbb{E}[|\tilde X_{s^{\prime}}|^{2p}] \le C_p \mathbb{E}\Big[ %
\Big|\int_{s^{\prime}\wedge t}^{s^{\prime}\wedge t^{\prime}} b(s';r)
dr\Big|^{2p} + \Big|\int_{s^{\prime}\wedge t}^{s^{\prime}\wedge t^{\prime}}
|\sigma(s';r)|^2 dr\Big|^{p}\Big] \\
&\le&C_p \mathbb{E}\Big[ \Big|\int_{s^{\prime}\wedge t}^{s^{\prime}\wedge
t^{\prime}} [1+\|X\|_T^{\k_0}] dr\Big|^{2p} + \Big|\int_{s^{\prime}\wedge
t}^{s^{\prime}\wedge t^{\prime}} [1+\|X\|_T^{2\k_0}] dr\Big|^{p}\Big] \\
&\le& C_p[s^{\prime}\wedge t^{\prime}-s^{\prime}\wedge t]^{p}\le
C_p(t^{\prime}-t)^{\frac{p}{2}} (s^{\prime}-s)^{\frac{p}{2}}.
\end{eqnarray*}

\textit{Case 2.} $t<s \le t^{\prime}$. Then 
\begin{eqnarray*}
I_p &=& \mathbb{E}\Big[\Big| \int_{t}^{s} \big[b(s;r) dr+ \sigma(s;r) dW_r %
\big] - \int_{t}^{s^{\prime}\wedge t^{\prime}}\big[ b(s';r)
dr+\sigma(s';r) dW_r\big]\Big|^{2p}\Big] \\
&\le& C_p\mathbb{E}\Big[\Big|\int_t^s [b(s;r)-b(s';r)] dr + \int_t^s
[\sigma(s;r)-\sigma(s';r)] dW_r\Big|^{2p} \\
&&\qquad + \Big|\int_s^{s^{\prime}\wedge t^{\prime}} b(s';r) dr+
\int_s^{s^{\prime}\wedge t^{\prime}} \sigma(s';r) dW_r\Big|^{2p}\Big]
\\
&\le& C_p\mathbb{E}\Big[\Big|\int_t^s [s^{\prime}-s][1+\|X\|_T^{\k_0}]dr\Big|%
^{2p} + \Big| \int_t^s [s^{\prime}-s]^2[1+\|X\|_T^{2\k_0}] dr\Big|^{p} \\
&&\qquad + \Big|\int_s^{s^{\prime}\wedge t^{\prime}} [1+\|X\|_T^{\k_0}] dr%
\Big|^{2p}+\Big| \int_s^{s^{\prime}\wedge t^{\prime}}[1+\|X\|_T^{2\k_0}] dr%
\Big|^{p}\Big] \\
&\le& C_p\Big[(s-t)^p(s^{\prime}-s)^{2p} + (s^{\prime}\wedge t^{\prime}-s)^p%
\Big]\le C_p(t^{\prime}-t)^{\frac{p}{2}} (s^{\prime}-s)^{\frac{p}{2}}.
\end{eqnarray*}

\textit{Case 3.} $s>t^{\prime }$. Then 
\begin{eqnarray*}
I_{p} &=&\mathbb{E}\Big[\Big|\int_{t}^{t^{\prime }}\big[b(s;r)dr+\sigma
(s;r)dW_{r}\big]-\int_{t}^{t^{\prime }}\big[b(s^{\prime }; r)dr+\sigma
(s^{\prime }; r)dW_{r}\big]\Big|^{2p}\Big] \\
&\leq &C_{p}\mathbb{E}\Big[\Big|\int_{t}^{t^{\prime }}[s^{\prime
}-s][1+\Vert X\Vert _{T}^{\k_{0}}]dr\Big|^{2p}+\Big|\int_{t}^{t^{\prime
}}[s^{\prime }-s]^{2}[1+\Vert X\Vert _{T}^{2\k_{0}}]dr\Big|^{p}\Big] \\
&\leq &C_{p}(t^{\prime }-t)^{p}(s^{\prime }-s)^{2p}\leq C_{p}(t^{\prime
}-t)^{\frac{p}{2}}(s^{\prime }-s)^{\frac{p}{2}}.
\end{eqnarray*}%
So in all the cases, we have proved \textrm{(\ref{DXest})}.

To see the second inequality in \textrm{(\ref{regXest})}, for each $i$, note that
\begin{equation*}
X_{t} - \Th^{t_i}_t =\int_{t_i}^{t}b(t;r)dr+\int_{t_i}^{t}\sigma (t;r)dW_{r},\quad t\ge t_i.
\end{equation*}
By Case 2 above we see that 
\begin{equation*}
\mathbb{E}\Big[\big|[X_{t}-\Theta _{t}^{t_{i}}]-[X_{t^{\prime }}-\Theta
_{t^{\prime }}^{t_{i}}]\big|^{2p}\Big]\leq C_{p}h^{\frac{p}{2}}|t-t^{\prime
}|^{\frac{p}{2}},\quad t_{i}\leq t<t^{\prime }\leq t_{i+1}.
\end{equation*}%
Then by Lemma \ref{lem-kolmogorov} we have 
$\mathbb{E}\Big[\sup_{t_{i}\leq t\leq t_{i+1}}|X_{t}-\Theta _{t}^{t_{i}}|^{4p}%
\Big]\leq C_{p}h^{p}.
$
Thus 
\begin{equation*}
\mathbb{E}[\Vert X-X^{n}\Vert _{T}^{8}]\leq \sum_{i=0}^{2^{n}-1}\mathbb{E}%
\Big[\sup_{t_{i}\leq t\leq t_{i+1}}|X_{t}-X_{t}^{n}|^{8}\Big]%
\le Ch^{2}2^{n}=C2^{-n},
\end{equation*} 
completing the proof.  \hfill \vrule width.25cm height.25cm depth0cm\smallskip

We need another technical lemma dealing with two variable functions/processes.

\begin{lem}
\label{lem-linearoperator}
Let Assumptions \ref{assum-X} and \ref{assum-regular} hold,  $u\in C^0(\ol\L)$, and $0\le t_1 < t_2\le T$. 

(i) If $\pa_\o u$ is continuous and has polynomial growth, then
\begin{equation}
\left.\begin{array}{c}
\dis \langle \pa_\o u(t_2, X\otimes_{t_1} \Th^{t_1}), \int_{t_1}^{t_2} b(\cd ~; r) dr\rangle  =\int_{t_1}^{t_2}  \langle \pa_\o u(t_2, X\otimes_{t_1} \Th^{t_1}), b^{r, X} \rangle dr;\\
\dis \langle \pa_\o u(t_2, X\otimes_{t_1} \Th^{t_1}), \int_{t_1}^{t_2} \si(\cd ~; r) dW_r\rangle  =\int_{t_1}^{t_2}  \langle \pa_\o u(t_2, X\otimes_{t_1} \Th^{t_1}), \si^{r, X} \rangle dW_r.
\end{array}
\right. \label{paolinear}
\end{equation}
Here, assuming $W$ is $k$-dimensional, $\langle \pa_\o u, \si\rangle dW_r := \sum_{i=1}^k \langle \pa_\o u, \si^i\rangle dW^i_r$, where $\si^i$ is the $i$-th column of $\si$.

(ii) If $\pa^2_{\o\o} u$ is continuous and has polynomial growth, then
\begin{eqnarray}
&& \Big\langle \pa^2_{\o\o} u(t_2, X\otimes_{t_1} \Th^{t_1}), \big(\int_{t_1}^{t_2} \si(\cd  ~; r) dW_r, \int_{t_1}^{t_2} \si(\cd ~; r) dW_r\big)\Big\rangle  \nonumber\\
 &=&  \sum_{i=1}^k \int_{t_1}^{t_2} \Big\langle \pa^2_{\o\o} u(t_2, X\otimes_{t_1} \Th^{t_1}), (\si^{i,t, X}, \si^{i, t,X}) \Big\rangle dt \\
&& + \int_{t_1}^{t_2}  \Big\langle \pa^2_{\o\o} u(t_2, X\otimes_{t_1} \Th^{t_1}), \big(\int_{t_1}^t \si(\cd ~; r) dW_r,  \si^{t,X}\big)\Big\rangle dW_t\nonumber \\
&& +  \int_{t_1}^{t_2}  \Big\langle \pa^2_{\o\o} u(t_2, X\otimes_{t_1} \Th^{t_1}), \big( \si^{t,X}, \int_{t_1}^t \si(\cd ~; r) dW_r\big)\Big\rangle dW_t.\nonumber
 \label{paobilinear}
\end{eqnarray}
\end{lem}
\proof \q For notational simplicity, we assume $d=k=1$, and  omit the variable $(t_2, X\otimes_{t_1} \Th^{t_1})$ inside $\pa_\o u$ and $\pa^2_{\o\o} u$.  

 (i) We prove the second equality in three steps. The first one follows similar  arguments.
 
 {\it Step 1.} Assume $\si(s; r) = \sum_{i=0}^{n-1} \si(s; r_i)\1_{[r_i, r_{i+1})}$ for some $t_1 = r_0< \cds<r_n=t_2$. Since $\pa_\o u$ is linear, the second equality of \reff{paolinear} is obvious.
 
 {\it Step 2.} Assume, for some constants $C$, $\k>0$ and for all $t_1 \le r<r'\le t_2$,
 \bea
 \label{sireg}
  |\si(s; r) - \si(s; r')| +  |\pa_s\si(s; r) - \pa_s\si(s; r')|\le C[1+\|\o\|_T^\k]|r-r'|.
  \eea
    Denote $\si_n(s;r) := \sum_{i=0}^{2^n-1} \si(s; r_i)\1_{[r_i, r_{i+1})}$, where $r_i := t_1 + i (t_2-t_1) 2^{-n}$, $i=0,\cds, 2^{-n}$.  Then $\sup_{t_2\le s\le T}[ |\si_n(s; r) - \si(s;r)| +  |\pa_s\si_n(s; r) - \pa_s\si(s;r)|\le C [1+\|\o\|_T^\k] 2^{-n}$ for all $r\in [t_1, t_2]$.   By \reff{paoubound}, this implies 
    $
    \lim_{n\to\infty}  \langle \pa_\o u, \si_n^{r, X} \rangle=  \langle \pa_\o u, \si^{r, X}\rangle,
    $
    which, together with the dominated convergence theorem, implies further that
    \bea
    \label{sinconv}
    \lim_{n\to\infty} \dbE\Big[\Big|\int_{t_1}^{t_2}  \langle \pa_\o u, \si_n^{r, X} \rangle dW_r - \int_{t_1}^{t_2}  \langle \pa_\o u, \si^{r, X} \rangle dW_r\Big|^2\Big]=0.
    \eea
Moreover,   denote $\tilde X_s := \int_{t_1}^{t_2} [\si_n(s; r) - \si(s;r)] dW_r$. Then,  for $t_2 \le s<s'\le T$ and $p\ge 1$,
 \beaa
&& \dbE\Big[\Big|\tilde X_s - \tilde X_{s'}|^{2p} \Big] = \dbE\Big[\Big| \int_{t_1}^{t_2} \int_s^{s'} \pa_s \si_n(l; r) - \pa_s \si(l; r) dl dW_r\Big|^{2p} \Big] \\
 &\le& C_p  \dbE\Big[\Big| \int_{t_1}^{t_2} \Big|\int_s^{s'} \pa_s \si_n(l; r) - \pa_s \si(l; r) dl\Big|^2 dr\Big|^{p} \Big] \le C_p 2^{-2pn} (s'-s)^{2p}.
 \eeaa 
 Applying Lemma \ref{lem-kolmogorov}  we get  $\dbE\Big[\sup_{t_2\le s\le T} |\tilde X_s|^2\Big] \le C 2^{-2n}$. Then $\lim_{n\to\infty}\sup_{t_2\le s\le T} |\tilde X_s|=0$, $\dbP$-a.s. and thus 
 \bea
 \label{sinconv2}
\lim_{n\to\infty} \langle \pa_\o u, \int_{t_1}^{t_2} \si_n(\cd ~; r) dW_r\rangle = \langle \pa_\o u, \int_{t_1}^{t_2} \si(\cd ~; r) dW_r\rangle,\q\dbP\mbox{-a.s.}
\eea
By Step 1, the second equality of \reff{paolinear} holds for $\si_n$. Then  by \reff{sinconv}   and \reff{sinconv2} we obtain the desired equality for $\si$.

 {\it Step 3.} Denote $\si_\e(s; r) := {1\over \e} \int_{(r-\e)^+}^r \si(s; l) dl$, and thus $\pa_s \si_\e(s; r) := {1\over \e} \int_{(r-\e)^+}^r \pa_s\si(s; l) dl$.  It is clear that $\lim_{\e\to 0}\dbE\big[\int_{t_1}^{t_2} [|\si_\e(s; r) - \si(s;r)|^p + | \pa_s\si_\e(s; r) - \pa_s\si(s;r)|^p] dr\big]=0$ for all $p\ge 1$.  Fix some $p$ large enough, then there exists $\e_n \downarrow 0$ such that 
 \beaa
\dbE\big[\int_{t_1}^{t_2} [|\si_{\e_n}(s; r) - \si(s;r)|^p + | \pa_s\si_{\e_n}(s; r) - \pa_s\si(s;r)|^p] dr\big] \le 2^{-n}.
\eeaa
Now following the arguments of Step 2 as well as that of Lemma \ref{lem-kolmogorov}, one can show that
\beaa
&   \lim_{n\to\infty} \dbE\big[\big|\int_{t_1}^{t_2}  \langle \pa_\o u, \si_{\e_n}^{r, X} \rangle dW_r - \int_{t_1}^{t_2}  \langle \pa_\o u, \si^{r, X} \rangle dW_r\big|^2\big]=0;&\\
&   \lim_{n\to\infty} \langle \pa_\o u, \int_{t_1}^{t_2} \si_{\e_n}(\cd ~; r) dW_r\rangle = \langle \pa_\o u, \int_{t_1}^{t_2} \si(\cd ~; r) dW_r\rangle,\q\mbox{a.s.}&
   \eeaa
Clearly $\si_\e$ satisfies the conditions in Step 2 and thus the second equality of \reff{paolinear} holds for each $\si_\e$. Then the above limits imply the desired equality for $\si$.

 (ii) Combining the arguments in Steps 2 and 3 in (i) above, it suffices to prove \reff{paobilinear} in the case $\si$ is piecewise constant in $r$: $\si(s; r) = \sum_{i=0}^{n-1} \si(s; r_i)\1_{[r_i, r_{i+1})}$ for some $t_1 = r_0< \cds<r_n=t_2$. In this case,   
$
 \int_{t_1}^{t_2} \si(s; r) dW_r =  \sum_{i=0}^{n-1} [W_{r_{i+1}}-W_{r_i}] \si(s; r_i).
$
Denote $W_{s, t} := W_t - W_s$ and $I_{ij} :=\langle \pa^2_{\o\o} u, (\si^{r_i, X}, \si^{r_j,X})\rangle$.  Since $\pa^2_{\o\o} u$ is bilinear, we see that
\beaa
&&\int_{t_1}^{t_2}  \Big\langle \pa^2_{\o\o} u, \big(\int_{t_1}^t \si(\cd ~; r) dW_r,  \si^{t,X}\big)\Big\rangle dW_t\\
&=& \sum_{i=0}^{n-1} \int_{r_i}^{r_{i+1}}  \Big\langle \pa^2_{\o\o} u, \big(\sum_{j=0}^{i-1} \si^{r_j, X} W_{r_j, r_{j+1}} + \si^{r_i, X} W_{r_i, t},  \si^{r_i,X}\big)\Big\rangle dW_t\\
&=& \sum_{i=0}^{n-1} \Big[\sum_{j=0}^{i-1} I_{ji} W_{r_i, r_{i+1}} W_{r_j, r_{j+1}} +  I_{ii} \int_{r_i}^{r_{i+1}} W_{r_i, t} dW_t\Big].
\eeaa
Then, by similar arguments for the last term of  \reff{paobilinear}, the right hand side of  \reff{paobilinear}  becomes
\beaa
&&\sum_{i=0}^{n-1} I_{ii} [r_{i+1}-r_i] +\sum_{i=0}^{n-1} \Big[\sum_{j=0}^{i-1} [I_{ji} + I_{ij}] W_{r_i, r_{i+1}} W_{r_j, r_{j+1}} + 2 I_{ii} \int_{r_i}^{r_{i+1}} W_{r_i, t} dW_t\Big]\\
&=&\sum_{i=0}^{n-1} I_{ii} W_{r_i, r_{i+1}}^2 + \sum_{0\le j< i\le n-1} W_{r_i, r_{i+1}}W_{r_j, r_{j+1}}[I_{ij} + I_{ji}]=\sum_{i=0}^{n-1}\sum_{j=0}^{n-1} W_{r_i, r_{i+1}}W_{r_j, r_{j+1}}I_{ij},
\eeaa
 which is equal to the left side of \reff{paobilinear}.
 \qed

\ms

We are now ready to prove our main result.

\no{\bf Proof of Theorem  \ref{thm-FIto1} and Proposition \ref{prop-path}. } 
As announced, we shall prove these two results together.
This does not create a circular argument. The first step is to prove Theorem  \ref{thm-FIto1}  on $\overline{%
\Lambda }$, which is not related to Proposition \ref{prop-path},
 then we prove Proposition \ref{prop-path}, and finally we
invoke Proposition \ref{prop-path} to draw the conclusion of Theorem  \ref{thm-FIto1}  as
it applies to $\Lambda $ instead of $\overline{\Lambda }$, where the
derivatives in (\ref{Fintro}) are uniquely determined because of Proposition %
\ref{prop-path}. For notational simplicity in this proof we assume all processes are scalar. The multidimensional case can be proved without any significant difficulty. Moreover, we emphasize again that we denote by $\kappa$ the generic polynomial growth order which may vary from line to line.

\textit{Step 1.} By abusing the notation slightly, in this step we assume $%
u\in C_{+}^{1,2}(\overline{\Lambda })$ and prove \textrm{(\ref{FIto})} for
such a function. This step does not refer to Proposition \ref{prop-path}
since it works in $C_{+}^{1,2}(\overline{\Lambda })$. Without loss of
generality, we shall only prove the result for $u(T,X)-u(0,0)$. Fix $n$ and
consider the setting in \textrm{(\ref{Xn})}. Then 
\begin{eqnarray}
&\displaystyle u(T,X)-u(0,0)=u(T,X)-u(T,X^{n})+%
\sum_{i=0}^{2^n-1}[I_{i}^{1}+I_{i}^{2}],\quad \mbox{where}&  \label{uT-0} \\
&I_{i}^{1}:=u(t_{i+1},X^{n}\otimes _{t_{i}}\Theta
^{t_{i}})-u(t_{i},X^{n}\otimes _{t_{i}}\Theta ^{t_{i}});&  \notag \\
&I_{i}^{2}:=u(t_{i+1},X^{n}\otimes _{t_{i+1}}\Theta
^{t_{i+1}})-u(t_{i+1},X^{n}\otimes _{t_{i}}\Theta ^{t_{i}}).&  \notag
\end{eqnarray}

First, by the second inequality in \textrm{(\ref{regXest})} we have 
\begin{equation*}
\mathbb{E}\Big[\sum_{n=1}^{\infty }\Vert X-X^{n}\Vert _{T}^{8}\Big]\leq
C\sum_{n=1}2^{-n}<\infty \quad \mbox{and thus}\quad \lim_{n\rightarrow
\infty }\Vert X-X^{n}\Vert _{T}=0,~\mathbb{P}\mbox{-a.s.}
\end{equation*}%
Since $u$ is continuous, we have 
\begin{equation}
\lim_{n\rightarrow \infty }[u(T,X)-u(T,X^{n})]=0,\quad \mathbb{P}\mbox{-a.s.}
\label{conv1}
\end{equation}

Next, by definition \textrm{(\ref{patu})} and the continuity of $\pa_t u$
we have $I_{i}^{1}=\int_{t_{i}}^{t_{i+1}}\partial _{t}u(t,X^{n}\otimes _{t_{i}}\Theta
^{t_{i}})dt.
$
By \textrm{(\ref{regXest})}, one can easily show that 
\begin{equation}
\lim_{n\rightarrow \infty } \sum_{i=0}^{2^n-1} \int_{t_i}^{t_{i+1}}\Vert X^{n}\otimes _{t_{i}}\Theta
^{t_{i}}-X\otimes _{t}\Theta ^{t}\Vert _{T} dt =0,\quad \mathbb{P}\mbox{-a.s.} \label{convint}
\end{equation}%
and thus, again by the continuity of $\pa_t u$ together with polynomial growth of $\pa_t u$ and \reff{ETh},
\begin{equation}
\lim_{n\rightarrow \infty
}\sum_{i=0}^{2^{n}-1}I_{i}^{1}=\int_{0}^{T}\partial _{t}u(t,X\otimes
_{t}\Theta ^{t})dt,\quad \mathbb{P}\mbox{-a.s.}  \label{I1conv}
\end{equation}

Moreover, note that 
\bea
\label{Xnti}
X^{n}\otimes _{t_{i}}\Theta ^{t_{i}}=X^{n}\otimes_{t_{i+1}}\Theta ^{t_{i}} =: X^{n, i}.
\eea
 Denote $\Delta \Theta ^{t_{i}}:=\Theta
^{t_{i+1}}-\Theta ^{t_{i}}$, then 
\bea
I_{i}^{2} &=& u(t_{i+1},X^{n}\otimes _{t_{i+1}}\Theta
^{t_{i+1}})-u(t_{i+1},X^{n,i}) \\
&=& \int_{0}^{1}\Big\la\partial _{\o }u\big(t_{i+1},X^{n}\otimes_{t_{i+1}}[\Theta ^{t_{i}} + \l \D \Th^{t_i}]\big),~\Delta
\Theta ^{t_{i}}\Big\ra d\lambda  = I_{i}^{2,1}+I_{i}^{2,2}+I_{i}^{2,3}, \nonumber  \label{I2}
\eea
where
\beaa
I_{i}^{2,1}& := &\Big\la\partial _{\o }u(t_{i+1},X^{n,i}),\Delta \Theta ^{t_{i}}\Big\ra  \\
I_{i}^{2,2}&:= & {\frac{1}{2}}\Big\la\partial _{\omega \omega }^{2}u(t_{i+1},X^{n,i}),~\big(\Delta \Theta
^{t_{i}},\Delta \Theta ^{t_{i}}\big)\Big\ra \\
I_{i}^{2,3} &:= &\Big\la\partial _{\omega \omega }^{2}u\big(t_{i+1},X^{n}\otimes _{t_{i+1}}[\Theta ^{t_{i}}+\lambda^* \Delta \Theta
^{t_{i}}]\big)-\partial _{\omega \omega }^{2}u(t_{i+1},X^{n,i}),~\big(\Delta \Theta ^{t_{i}},\Delta \Theta
^{t_{i}}\big)\Big\ra,
\eeaa
for some appropriate (random) $\lambda^*$ taking values on $[0,1]$. Note that 
\begin{equation}
\Delta \Theta
_{s}^{t_{i}}=\int_{t_{i}}^{t_{i+1}}b(s;r)dr+\int_{t_{i}}^{t_{i+1}}\sigma
(s;r)dW_{r},\quad s\geq t_{i+1}.\label{DTh}
\end{equation}%
Since $X^{n}\otimes
_{t_{i+1}}\Theta ^{t_{i}} = X^{n}\otimes _{t_{i}}\Theta ^{t_{i}}$ is $\mathcal{F}_{t_{i}}$-measurable, similar to Lemma \ref{lem-linearoperator} (i) we have 
\begin{equation*}
I_{i}^{2,1}=\int_{t_{i}}^{t_{i+1}}\big\la\partial _{\o }u\big(%
t_{i+1},X^{n,i}\big),~b^{r,X}\Big\ra %
dr+\int_{t_{i}}^{t_{i+1}}\Big\la\partial _{\o }u\big(t_{i+1},X^{n,i}\big),~\sigma ^{r,X}\big\ra dW_{r}.
\end{equation*}%
Recall \reff{Xnti} and  \textrm{(\ref{convint})}. By Assumptions \ref{assum-X} and \ref{assum-regular}, since $\pa_\o u$ is continuous, we have 
\begin{eqnarray*}
&&\lim_{n\rightarrow \infty }\sum_{i=0}^{2^{n}-1}%
\int_{t_{i}}^{t_{i+1}}\Big|\big\la\partial _{\o }u\big(t_{i+1},X^{n,i}\big),~b^{r,X}\big\ra-\big\la\partial _{\o %
}u(r,X\otimes _{r}\Theta ^{r}),~b^{r,X}\big\ra\Big|dr=0; \\
&&\lim_{n\rightarrow \infty }\sum_{i=0}^{2^{n}-1}%
\int_{t_{i}}^{t_{i+1}}\Big|\big\la\partial _{\o }u\big(t_{i+1},X^{n,i}\big),~\sigma ^{r,X}\big\ra-\big\la\partial _{\o %
}u(r,X\otimes _{r}\Theta ^{r}),~\sigma ^{r,X}\big\ra\Big|^{2}dr=0.
\end{eqnarray*}%
Therefore, with convergence in  $\dbL^2$,
\begin{equation}
\lim_{n\rightarrow \infty }\sum_{i=0}^{2^{n}-1}I_{i}^{2,1}=\int_{0}^{T}%
\big\la\partial _{\o }u(t,X\otimes _{t}\Theta ^{t}),~b^{t,X}\big\ra %
dt+\int_{0}^{T}\big\la\partial _{\o }u(t,X\otimes _{t}\Theta ^{t}),~\sigma
^{t,X}\big\ra dW_{t}.  \label{I21conv}
\end{equation}%

We now consider $I^{2,2}_i$.  In the spirit of  Lemma \ref{lem-linearoperator} (ii) we can prove
\bea
&\dis 2I_{i}^{2,2} = I^{2,2,1}_i + I^{2,2,2}_i + I^{2,2,3}_i,\q\mbox{where}&\label{I22}\\
&\dis I^{2,2,1}_i := \Big\la\partial _{\omega \omega }^{2}u( t_{i+1},X^{n,i}),~\big(\int_{t_{i}}^{t_{i+1}}b(\cd ~;r)dr, \int_{t_{i}}^{t_{i+1}}b(\cd ~;r)dr\big)\Big\ra&\nonumber\\
&\dis \qq\qq\q +\Big\la\partial _{\omega \omega }^{2}u( t_{i+1},X^{n,i}),~\big(\int_{t_{i}}^{t_{i+1}}b(\cd ~;r)dr, \int_{t_{i}}^{t_{i+1}}\si(\cd ~;r)dW_r\big)\Big\ra&\nonumber\\
&\dis\qq\qq\q+\Big\la\partial _{\omega \omega }^{2}u( t_{i+1},X^{n,i}),~\big(\int_{t_{i}}^{t_{i+1}}\si(\cd ~;r)dW_r, \int_{t_{i}}^{t_{i+1}}b(\cd ~;r)dr\big)\Big\ra&\nonumber\\
&\dis I^{2,2,2}_i:=\int_{t_i}^{t_{i+1}}\Big\la\partial _{\omega \omega }^{2}u\big(t_{i+1},X^{n,i}\big),~\big(\int_{t_{i}}^t \si(\cd ~;r)dW_r, \si^{t, X}\big)\Big\ra dW_t&\nonumber\\
&\dis\qq\qq+\int_{t_i}^{t_{i+1}}\Big\la\partial _{\omega \omega }^{2}u\big(t_{i+1},X^{n,i}\big),~\big(\si^{t, X}, \int_{t_{i}}^t \si(\cd ~;r)dW_r\big)\Big\ra dW_t&\nonumber\\
&\dis I^{2,2,3}_i:=\int_{t_i}^{t_{i+1}}\Big\la\partial _{\omega \omega }^{2}u\big(t_{i+1},X^{n,i}\big),~\big(\si^{t, X}, \si^{t,X}\big)\Big\ra dt.\qq\qq\q &\nonumber
\eea
One can similarly show that, with convergence in $\dbL^2$,
\begin{equation}
\lim_{n\to\infty} \sum_{i=0}^{2^n-1}  I^{2,2,3}_i = \int_0^T \Big\la\partial _{\omega \omega }^{2}u\big(t,X\otimes
_t\Theta ^t\big),~\big(\si^{t, X}, \si^{t,X}\big)\Big\ra dt. \label{I223}
\end{equation}
By the martingale property,
\begin{eqnarray*}
 \dbE\Big[\Big|\sum_{i=0}^{2^n-1}I^{2,2,2}_i\Big|^2\Big] 
&=& \sum_{i=0}^{2^n-1}\dbE\Big[ \int_{t_i}^{t_{i+1}}\Big|\Big\la\partial _{\omega \omega }^{2}u\big(t_{i+1},X^{n,i}\big),~\big(\int_{t_{i}}^t \si(\cd ~;r)dW_r, \si^{t, X}\big)\Big\ra\Big|^2 dt \\
&&+\int_{t_i}^{t_{i+1}}\Big|\Big\la\partial _{\omega \omega }^{2}u\big(t_{i+1},X^{n,i}\big),~\big(\si^{t, X}, \int_{t_{i}}^t \si(\cd ~;r)dW_r\big)\Big\ra\Big|^2 dt\Big] \\
&\le&C \sum_{i=0}^{2^n-1}\dbE\Big[[1+\|X^{n,i}\|_T^\kappa+\|X\|_T^{\kappa}] \int_{t_i}^{t_{i+1}} \sup_{t_{i+1}\le s\le T} \big|\int_{t_i}^t \si(s;r)dW_r\big|^2 dt\Big]\\
&\le& C\sum_{i=0}^{2^n-1}\Big(\dbE\Big[ 2^{-n} \int_{t_i}^{t_{i+1}} \sup_{t_{i+1}\le s\le T} \big|\int_{t_i}^t \si(s;r)dW_r\big|^4 dt\Big]\Big)^{1\over 2}.
\end{eqnarray*}
Note that, for $t_i \le t\le t_{i+1} \le s<s' \le T$ and $p\ge 1$,
\begin{eqnarray*}
&&  \dbE\Big[\big|\int_{t_i}^t \si(s;r)dW_r - \int_{t_i}^t \si(s';r)dW_r\big|^{2p} \Big]\le C_p\dbE\Big[\big[\int_{t_i}^t |\si(s;r)- \si(s';r)|^2dr\big]^{p}\Big] \\
&&\le C_p\dbE\Big[\big[\int_{t_i}^t [\int_s^{s'} |\pa_l\si(l,r)|dl]^2dr\big]^{p}\Big]\le C_p(t-t_i)^p (s'-s)^{2p}.
\end{eqnarray*}
Applying Lemma \ref{lem-kolmogorov}, we have
\begin{equation}
\dbE\Big[\sup_{t_{i+1}\le s\le T} \big|\int_{t_i}^t \si(s;r)dW_r\big|^4 dt\Big] \le C (t-t_i)^2 \le C2^{-2n}. \label{SFT}
\end{equation}
Then
\begin{equation}
 \dbE\Big[\Big|\sum_{i=0}^{2^n-1}I^{2,2,2}_i\Big|^2\Big] \le C\sum_{i=0}^{2^n-1}\Big( 2^{-n} \int_{t_i}^{t_{i+1}} 2^{-2n} dt\Big)^{1\over 2}= C\sum_{i=0}^{2^n-1} 2^{-2n} = C2^{-n} \to 0. \label{I222}
\end{equation}
 Moreover,  by \reff{SFT} again,
\begin{eqnarray*}
&& \dbE\Big[\Big|\sum_{i=0}^{2^n-1}I^{2,2,1}_i\Big|^2\Big] \le 2^n \sum_{i=0}^{2^n-1}\dbE[I^{2,2,1}_i|^2]\\
&\le& C 2^n \sum_{i=0}^{2^n-1}\dbE\Big[[1+\|X^{n,i}\|_T^\kappa+\|X\|_T^{\kappa}] \big[2^{-4n} + 2^{-2n}\sup_{t_{i+1}\le s\le T} |\int_{t_i}^{t_{i+1}} \si(s; r) dW_r|^2\big] \Big]\\
&\le& C2^{-2n} + C2^{-n} \sum_{i=0}^{2^n-1} \Big(\dbE\Big[\sup_{t_{i+1}\le s\le T} |\int_{t_i}^{t_{i+1}} \si(s; r) dW_r|^4 \Big]\Big)^{1\over 2}\\
&\le&  C2^{-2n} + C2^{-n} \sum_{i=0}^{2^n-1} 2^{-n} \le C2^{-n}.
\end{eqnarray*}
Plug this and \reff{I223}, \reff{I222} into \reff{I22}, we obtain 
\begin{equation}
\lim_{n\to\infty} \sum_{i=0}^{2^n-1}  I^{2,2}_i = {1\over 2}\int_0^T \Big\la\partial _{\omega \omega }^{2}u\big(t,X\otimes
_t\Theta ^t\big),~\big(\si^{t, X}, \si^{t,X}\big)\Big\ra dt. \label{I22conv}
\end{equation}

Finally, denote $\|\D \Th_i\|:= \sup_{t_{i+1}\leq s\leq T}|\Delta
\Theta _{s}^{t_{i}}|$.  For any $\varepsilon >0$, by  \textrm{(\ref{paooucont})} we have
\begin{equation*}
|I_{i}^{2,3}| \leq \rho(\|\D\Th_i\|) \|\D\Th_i\|^2  \le \rho(\e)  \|\D\Th_i\|^2 + C\e^{-1} \|\D \Th_i\|^3.
\end{equation*}
By \reff{DTh}, similar to \reff{SFT} we have
\begin{equation*}
\dbE[|I_{i}^{2,3}|]  \le \rho(\e)  \dbE[\|\D\Th_i\|^2] + C\e^{-1} \dbE[\|\D \Th_i\|^3]\le  C\rho(\e)  2^{-n} + C\e^{-1} 2^{-{3\over 2}n}.
\end{equation*}
Then
\begin{equation*}
\dbE\Big[\sum_{i=0}^{2^n-1}|I_{i}^{2,3}|\Big]  \le   C\rho(\e)   + C\e^{-1} 2^{-{n\over 2}}.
\end{equation*}
By first sending $n\to \infty$ and then $\e\to 0$, we obtain $\lim_{n\to\infty} \dbE\Big[\sum_{i=0}^{2^n-1}|I_{i}^{2,3}|\Big]  =0$.  
Plug this and  \reff{I21conv}, \reff{I22conv} into \textrm{(\ref{I2})}, we get 
\begin{eqnarray*}
\lim_{n\rightarrow \infty }\sum_{i=0}^{2^{n}-1}I_{i}^{2} &=&\int_{0}^{T}%
\big\la\partial _{\o }u(t,X\otimes _{t}\Theta ^{t}),~b^{t,X}\big\ra %
dt+\int_{0}^{T}\big\la\partial _{\o }u(t,X\otimes _{t}\Theta ^{t}),~\sigma
^{t,X}\big\ra dW_{t} \\
&&+{\frac{1}{2}}\int_{0}^{T}\big\la\partial _{\omega \omega
}^{2}u(t,X\otimes _{t}\Theta ^{t}),~(\sigma ^{t,X},\sigma ^{t,X})\big\ra %
dt,\quad \mathbb{P}\mbox{-a.s.}
\end{eqnarray*}%
This, together with \textrm{(\ref{conv1})} and \textrm{(\ref{I1conv})},
proves \textrm{(\ref{FIto})} for $u\in C_+^{1,2}(\overline{\Lambda })$.

\textit{Step 2.} We now prove Proposition \ref{prop-path}. Set $%
u:=u^{1}-u^{2}\in C_{+}^{1,2}(\overline{\Lambda })$ and thus $u=0$ on $%
\Lambda $. By definition in \textrm{(\ref{patu}),} it is clear that $%
\partial _{t}u=0$ on $\Lambda $. Fix $(t_{0},\omega ^{0})\in \Lambda $ and $%
\eta \in \Omega _{t_{0}}$. Define 
\beaa
\tilde{X}_{t}:=\omega _{t}^{0}\mathbf{1}_{[0,t_{0})}(t)+\big[\omega _{t_{0}}^{0}+\int_{t_{0}}^{t}\eta _{r}dW_{r}\big]\mathbf{1}_{[t_{0},T]}(t),
~ \tilde{\Theta}_{s}^{t}:=\omega _{s}^{0}+\int_{t_{0}}^{t}\eta_{r}dW_{r},~ t_{0}\leq t\leq s\leq T.
\eeaa
Note that $\tilde{X}\otimes _{t}\tilde{\Theta}^{t}$ is continuous, thus $u(t,%
\tilde{X}\otimes _{t}\tilde{\Theta}^{t})=0$, $t_{0}\leq t\leq T$. By Step 1,
which does not require Proposition \ref{prop-path}, the process $u(t,\tilde{X%
}\otimes _{t}\tilde{\Theta}^{t})$ satisfies \textrm{(\ref{FIto})} on $%
[t_{0},T]$ and thus 
\begin{equation*}
\langle \partial _{\o }u(t,\tilde{X}\otimes _{t}\tilde{\Theta}^{t}),\eta
\rangle =0,\quad \langle \partial _{\omega \omega }^{2}u(t,\tilde{X}\otimes
_{t}\tilde{\Theta}^{t}),(\eta ,\eta )\rangle =0,\quad t_{0}\leq t\leq T,%
\mathbb{P}\mbox{-a.s.}
\end{equation*}%
In particular, noting that $\tilde{X}\otimes _{t_{0}}\tilde{\Theta}%
^{t_{0}}=\omega ^{0}$, then for $t=t_{0}$ we have 
\begin{equation*}
\langle \partial _{\o }u(t_{0},\omega ^{0}),\eta \rangle =0,\quad \langle
\partial _{\omega \omega }^{2}u(t_{0},\omega ^{0}),(\eta ,\eta )\rangle
=0,\quad \mathbb{P}\mbox{-a.s.}
\end{equation*}%
Since $\eta $ is arbitrary,  we prove Proposition \ref{prop-path}.

\textit{Step 3}. Finally it is clear that \textrm{(\ref{FIto})} holds for $%
u\in C_{+}^{1,2}(\Lambda )$. In particular, by Proposition \ref{prop-path} (or say Step 2 above),
the path derivatives in the right hand side of \textrm{(\ref{FIto})} do not
depend on the choice of $\tilde{u}\in C_+^{1,2}(\overline{\Lambda })$. \hfill %
\vrule width.25cm height.25cm depth0cm\smallskip 

\begin{rem}
\label{rem-C12} \textrm{An alternate way to define path derivatives is
directly through (\ref{FIto}) by positing that this functional It\^{o}
formula holds, and the derivatives will be uniquely defined in appropriate
sense. This is the approach in \cite{ETZ1, ETZ2} for PPDEs in a
semmartingale framework. In this way we may avoid involving the {c\`{a}dl%
\`{a}g} space }$\overline{\Lambda }$\textrm{. \hfill \vrule width.25cm
height.25cm depth0cm\smallskip }
\end{rem}

\subsection{The singular case\label{singular}}

We now consider the case where $b(t; t)$ and $\sigma (t; t)$ may blow up. 
We shall assume the following growth condition which is
modeled after the behavior of the kernel for Gaussian processes with
self-similarity parameter $H\in (0,1/2)$, and is therefore satisfied by fBm
with $H$ in that range (and is more true for fBM with larger hurst parameter by setting $H={1\over 2}$ at below).

\begin{assum}
\label{assum-irregular}   For $\varphi =b,\sigma $,  $\partial _{t}\f(t; s,\cdot )$ exists for $t\in (s,T]$, and there exists $0<H<{\frac{1}{2}}$
such that, for any $0\leq s<t\leq T$, 
\begin{equation}
|\varphi (t; s,\omega )|\leq C_{0}[1+\Vert \omega \Vert _{T}^{\kappa_{0}}](t-s)^{H-%
{\frac{1}{2}}},\quad |\partial _{t}\varphi (t; s,\omega )|\leq C_{0}[1+\Vert
\omega \Vert _{T}^{\kappa_{0}}](t-s)^{H-{\frac{3}{2}}}.  \label{KH}
\end{equation}
\end{assum}

We remark that, in this case $b^{t,X},\sigma ^{t,X}$ are not in $\Omega _{t}$
and thus cannot serve as the test function in the right side of \textrm{(\ref%
{FIto})}. To overcome this difficulty, we assume the following conditions on 
$u$ which roughly mean that $u(t,\omega )$ does not depend on $\{\omega_s
\}_{t\leq s\leq t+\delta }$ for some small $\delta >0$, or depends very
weakly on the paths in the sense that $u$'s derivatives become increasingly
smaller as one approaches the diagonal.

\begin{defn}
\label{defn-C12a}
We say $u\in C_{+}^{1,2}(\Lambda )$ vanishes diagonally with rate $\a\in (0,1)$, denoted as $u\in C_{+,\a}^{1,2}(\Lambda )$ if there exist an extension of $u$ in $C_{+}^{1,2}(\overline{\Lambda })$, still denoted as $u$, a polynomial growth order $\kappa$, and a bounded modulus of continuity function $\rho$   satisfying: for any $0\leq t<T$, $0<\delta
\leq T-t$, and $\eta ,\eta _{1},\eta _{2}\in \Omega _{t}$ with the supports
of $\eta $, $\eta _{1}$, and $\eta _{2}$ contained in $[t,t+\delta ]$,

(i) for any $\omega \in \overline{\O }$ such that $\omega \mathbf{1}%
_{[t,T]}\in \Omega _{t}$, 
\begin{equation}
\left.\begin{array}{c}
\displaystyle\big|\langle \partial _{\o }u(t,\omega ),\eta \rangle \big|\leq
C[1+\|\o\|_T^\kappa] \Vert \eta \Vert _T\delta ^{\alpha },\\
\displaystyle \big|\langle \partial
_{\omega \omega }^{2}u(t,\omega ),(\eta _{1},\eta _{2})\rangle \big|\leq
C[1+\|\o\|_T^\kappa]\big\Vert |\eta _{1}||\eta _{2}|\big\Vert _T\delta ^{2\alpha }.
\end{array}\right.
\label{vanish}
\end{equation}

(ii) for any other $\omega ^{\prime }\in \overline{\O }$ such that $%
\omega ^{\prime }\mathbf{1}_{[t,T]}\in \Omega _{t}$
\bea
&\Big|\langle \partial _{\o }u(t,\omega )-\partial _{\o }u(t,\omega ^{\prime }),\eta \rangle \Big|
\leq [1+\|\o\|_T^\kappa + \|\o'\|_T^\kappa] \Vert \eta \Vert _T\rho(\Vert \omega -\omega
^{\prime }\Vert _{T})\delta
^{\alpha }, \label{vanish2} \\ 
&\Big|\langle \partial _{\omega \omega }^{2}u(t,\omega
)-\partial _{\omega \omega }^{2}u(t,\omega ^{\prime }),(\eta _{1},\eta
_{2})\rangle \Big|\leq   [1+\|\o\|_T^\kappa + \|\o'\|_T^\kappa] \big\Vert |\eta _{1}||\eta _{2}|\big\Vert _T \rho(\Vert \omega -\omega ^{\prime }\Vert
_{T})\delta
^{2\alpha }.\nonumber
\eea
\end{defn}

These conditions will allow us to truncate the coefficients $b,\sigma $ near
the diagonal, and control the error made by this truncation. For $\varphi
=b,\sigma $ and $\delta >0$, we introduce the truncated functions: 
\begin{equation*}
\varphi ^{\delta }(t; s,\omega ):=\varphi (t\vee (s+\delta ); s,\omega ).
\end{equation*}%
We also again use the notation $\varphi ^{\delta ,s,\omega }$ for the path $%
t \in [s, T]\mapsto \varphi ^{\delta }(t; s,\omega )$. Another consequence of using
these truncated coefficients is that the notion of time and path derivatives
must be understood as limits when the truncation parameter $\delta $ tends
to $0$. Specifically, we prove the following functional It\^{o} formula, where
in particular, the said limits exist.

\begin{thm}
\label{thm-FIto2} Let Assumptions \ref{assum-X} and \ref{assum-irregular} hold. Assume $u\in C^{1,2}_{+,\a}(\L)$ with $\b:=\alpha +H-{\frac{1}{2}}>0$. Then the functional
It\^{o} formula \textrm{(\ref{FIto})} still holds true, where 
\begin{equation}
\left. 
\begin{array}{c}
\displaystyle\langle \partial _{\o }u(t,\omega ),\varphi ^{t,\omega }\rangle
:=\lim_{\delta \downarrow 0}\langle \partial _{\o }u(t,\omega ),\varphi
^{\delta ,t,\omega }\rangle ,\quad \varphi =b,\sigma  \\ 
\displaystyle\langle \partial _{\omega \omega }^{2}u(t,\omega ),(\sigma
^{t,\omega},\sigma ^{t,\omega })\rangle :=\lim_{\delta \downarrow 0}\langle
\partial _{\omega \omega }^{2}u(t,\omega ),(\sigma ^{\delta ,t,\omega
},\sigma ^{\delta ,t,\omega })\rangle .%
\end{array}%
\right.   \label{pathu-singular}
\end{equation}
\end{thm}

{\noindent \textbf{Proof\quad}} We proceed in three steps.

\textit{Step 1.} We first show that the limits in \textrm{(\ref%
{pathu-singular})} exist. We shall only prove it for $\sigma $, and the
result for $b$ follow the same arguments. Denote 
\begin{equation}
\delta _{n}:={\frac{1}{2^{n}}},t_{n}:=t+\delta _{n},~\psi _{n}(s):={\frac{%
s-t _{n+1}}{\delta _{n}-\delta _{n+1}}}\mathbf{1}%
_{(t_{n+1},t_{n}]}(s)+{\frac{t _{n-1}-s}{\delta _{n-1}-\delta _{n}}}%
\mathbf{1}_{(t_{n},t_{n-1})}(s).  \label{dn}
\end{equation}%
Then $\psi _{n}$ is continuous, with support $(t_{n+1},t_{n-1})$, and $%
\psi _{n}+\psi _{n+1}=1$ on $[t_{n+1},t_{n}]$. Now for any $\delta
^{\prime }<\delta $, assume $\delta ^{\prime }\in \lbrack \d_{n+1}, \d_n)$
and $\delta \in \lbrack \d_{m+1},\d_m)$ for some $m\leq n$. Consider the
following decomposition into continuous functions for the constant $1$: 
\begin{equation}
\mathbf{1}_{[t,T]}=[1-\psi _{m}]\mathbf{1}_{[t_{m},T]}+\sum_{k=m}^{n}\psi
_{k}+[1-\psi _{n}]\mathbf{1}_{[t,t_{n}]}.  \label{fn}
\end{equation}%
Note that $\sigma _{s}^{\delta ,t,\omega }=\sigma _{s}^{\delta ^{\prime
},t,\omega }$ for $s\in \lbrack t_{m},T]$. Then, for $s\geq t$, 
\begin{equation}
\sigma _{s}^{\delta ,t,\omega }-\sigma _{s}^{\delta ^{\prime },t,\omega
}=\sum_{k=m}^{n}[\sigma _{s}^{\delta ,t,\omega }-\sigma _{s}^{\delta
^{\prime },t,\omega }]\psi _{k}(s)+[\sigma _{s}^{\delta ,t,\omega
}-\sigma _{s}^{\delta ^{\prime },t,\omega }][1-\psi _{n}]\mathbf{1}%
_{[t,t_{n}]}.  \label{sid}
\end{equation}%
Thus, by  the first inequalities of \textrm{(\ref{vanish})} and  \textrm{(\ref{KH})},
\begin{eqnarray}
&&\Big|\langle \partial _{\o }u(t,\omega ),\sigma ^{\delta ,t,\omega
}\rangle -\langle \partial _{\o }u(t,\omega ),\sigma ^{\delta ^{\prime
},t,\omega }\rangle \Big|  \notag \\
&\leq &\sum_{k=m}^{n}\Big|\langle \partial _{\o }u(t,\omega ),\psi
_{k}[\sigma ^{\delta ,t,\omega }-\sigma ^{\delta ^{\prime },t,\omega
}]\rangle \Big|+\Big|\langle \partial _{\o }u(t,\omega ),[1-\psi _{n}]%
\mathbf{1}_{[t,t_{n}]}[\sigma ^{\delta ,t,\omega }-\sigma ^{\delta ^{\prime
},t,\omega }]\rangle \Big|  \notag \\
&\leq &C[1+\Vert \omega \Vert _{T}^\kappa] \Big[\sum_{k=m}^{n}\sup_{t_{k+1}\leq s\leq t_{k-1}}[|\sigma _{s}^{\delta
,t,\omega }|+|\sigma _{s}^{\delta ^{\prime },t,\omega }|]\delta
_{k-1}^{\alpha }+\sup_{t\leq s\leq t_{n-1}}[|\sigma _{s}^{\delta ,t,\omega
}|+|\sigma _{s}^{\delta ^{\prime },t,\omega }|]\delta _{n-1}^{\alpha } \Big]
\notag \\
&\leq& C[1+\Vert \omega \Vert _{T}^\k]\Big[\sum_{k=m}^{n}\delta
_{k+1}^{H-{\frac{1}{2}}}\delta _{k-1}^{\alpha }+\delta _{n+1}^{H-{\frac{1}{2}%
}}\delta _{n-1}^{\alpha }\Big]\leq C[1+\Vert \omega \Vert
_{T}^\k]\sum_{k=m}^{\infty }2^{-\b k}  \notag \\
&\leq &C[1+\Vert \omega \Vert _{T}^\kappa]2^{-\b m}\leq C[1+\Vert \omega
\Vert _{T}^\k]\delta ^{\b}\rightarrow 0,\quad \mbox{as}~\delta
\rightarrow 0.  \label{pawconv}
\end{eqnarray}
Similarly, by the second inequality of \textrm{(\ref{vanish})} and the first inequality of  \textrm{(\ref{KH})}, we have
\begin{eqnarray}
&& \!\!  \Big|\langle \partial _{\omega \omega }^{2}u(t,\omega ),(\sigma ^{\delta
,t,\omega },\sigma ^{\delta ,t,\omega })\rangle -\langle \partial _{\omega
\omega }^{2}u(t,\omega ),(\sigma ^{\delta ^{\prime },t,\omega },\sigma
^{\delta ^{\prime },t,\omega })\rangle \Big|  \notag \\
&\leq & \!\! \Big|\langle \partial _{\omega \omega }^{2}u(t,\omega ),(\sigma
^{\delta ,t,\omega },\sigma ^{\delta ,t,\omega }-\sigma ^{\delta ^{\prime
},t,\omega })\rangle \Big|+\Big|\langle \partial _{\omega \omega
}^{2}u(t,\omega ),(\sigma ^{\delta ,t,\omega }-\sigma ^{\delta ^{\prime
},t,\omega },\sigma ^{\delta ^{\prime },t,\omega })\rangle \Big|  \notag \\
&\leq & \!\! C[1+\Vert \omega \Vert _{T}^\kappa] \Big[\sum_{k=m}^{n}\sup_{t_{k+1}\leq s\leq t_{k-1}}[|\sigma _{s}^{\delta
,t,\omega }|+|\sigma _{s}^{\delta ^{\prime },t,\omega }|]^{2}\delta
_{k-1}^{2\alpha }+\!\! \sup_{t\leq s\leq t_{n-1}}[|\sigma _{s}^{\delta ,t,\omega
}|+|\sigma _{s}^{\delta ^{\prime },t,\omega }|]^{2}\delta _{n-1}^{2\alpha } \Big]
\notag \\
&\leq & \!\!  C[1+\Vert \omega \Vert _{T}^\k]\Big[\sum_{k=m}^{n}\delta
_{k+1}^{2H-1}\delta _{k-1}^{2\alpha } + \delta _{n+1}^{2H-1}\delta
_{n-1}^{2\alpha }\Big]  \notag \\
&\leq & \!\!  C[1+\Vert \omega \Vert _{T}^\k]\sum_{k=m}^{\infty }2^{-2\b k}\leq
C[1+\Vert \omega \Vert _{T}^\k]\delta ^{2\b}\rightarrow 0,\quad %
\mbox{as}~\delta \rightarrow 0.  \label{pawconv2}
\end{eqnarray}%
This, together with \textrm{(\ref{pawconv})}, implies \textrm{(\ref%
{pathu-singular})}. Moreover, by sending $\delta ^{\prime }\rightarrow 0$,
we obtain the following estimates which are stronger than \textrm{(\ref%
{pathu-singular})}: 
\begin{equation}
\left. 
\begin{array}{c}
\Big|\langle \partial _{\o }u(t,\omega ),\varphi ^{\delta ,t,\omega }\rangle
-\langle\partial _{\o }u(t,\omega ),\varphi ^{t,\omega }\rangle \Big|\leq C[1+\Vert
\omega \Vert _{T}^\k]\delta ^{\b},\quad \varphi =b,\sigma , \\ 
\Big|\langle \partial _{\omega \omega }^{2}u(t,\omega ),(\sigma ^{\delta
,t,\omega },\sigma ^{\delta ,t,\omega })\rangle -\langle\partial _{\omega \omega
}^{2}u(t,\omega ),(\sigma^{t,\omega },\sigma ^{t,\omega })\rangle \Big|%
\leq C[1+\Vert \omega \Vert _{T}^\k]\delta ^{2\b}.%
\end{array}%
\right.   \label{pauKconv}
\end{equation}

\textit{Step 2.} Denote%
\bea
\label{XTHd}
\left. \begin{array}{c}
\displaystyle X_{t}^{\delta }:=x+\int_{0}^{t}b^{\delta }(t; r,X)dr+\int_{0}^{t}\sigma ^{\delta }(t; r,X)dW_{r}, \\ 
\displaystyle\Theta _{s}^{\delta ,t}:=x+\int_{0}^{t}b^{\delta
}(s;r,X)dr+\int_{0}^{t}\sigma ^{\delta }(s;r,X)dW_{r}.
\end{array} \right. 
\eea
We emphasize that in the above definitions, $b^{\delta }$, $\sigma ^{\delta }
$ depend on $X$, not $X^{\delta }$, since the truncation occurs in the first
two parameters of $b,\sigma $ only. In particular, $X^{\delta }$ is explicit
given by $X$, and does not solve an SDE. For notational simplicity at below we
shall still omit $X$ in the coefficients $b,\sigma $. In this step we prove 
\begin{equation}
\mathbb{E}[\Vert X^{\delta }\otimes _{t}\Theta ^{\delta ,t}-X\otimes
_{t}\Theta ^{t}\Vert _{T}^{2p}]\leq C_{p}\delta ^{pH},\quad \mbox{for any
$0\le t\le T$, $p\ge 1$, and $\d>0$}.  \label{Xdconv}
\end{equation}

We first estimate the difference of $\Th$. By stochastic Fubini theorem we have
\begin{eqnarray*}
\sup_{t\leq s\leq T}|\Theta _{s}^{\delta ,t}-\Theta _{s}^{t}| 
\!\!\!\!&=&\!\!\!\! \sup_{t\leq s\leq t+\delta }\Big|\int_{s-\delta }^{t}\Big[[b(r+\delta; r)-b(s;r)]dr+[\sigma (r+\delta; r)-\sigma (s;r)]dW_{r}\Big]\Big| \\
\!\!\!\!&=&\!\!\!\! \sup_{t\leq s\leq t+\delta }\Big|\int_{s-\delta }^{t}\int_{s}^{r+\delta }%
\Big[\partial _{t}b(\lambda; r)d\lambda dr+\partial _{t}\sigma (\lambda; r)d\lambda dW_{r}\Big]\Big| \\
\!\!\!\!&=&\!\!\!\! \sup_{t\leq s\leq t+\delta }\Big|\int_{s}^{t+\delta }\int_{s-\delta
}^{\lambda -\delta }\big[\partial _{t}b(\lambda; r)dr+\partial _{t}\sigma
(\lambda; r)dW_{r}\big]d\lambda \Big| \\
\!\!\!\!&\le&\!\!\!\! \int_{t}^{t+\delta }\sup_{t-\delta \leq l\leq \lambda -\delta }\Big|%
\int_{l}^{\lambda -\delta }\big[\partial _{t}b(\lambda; r)dr+\partial
_{t}\sigma (\lambda; r)dW_{r}\big]\Big|d\lambda 
\end{eqnarray*}%
Then, for any $p\geq 1$,  by Burkholder-Davis-Gundy inequality and the second inequality of  \textrm{(\ref{KH})} we obtain
\beaa
&&\mathbb{E}\Big[\sup_{t\leq s\leq T}|\Theta _{s}^{\delta ,t}-\Theta
_{s}^{t}|^{2p}\Big]  \notag \\
&\leq &C_{p}\delta ^{2p-1}\int_{t}^{t+\delta }\mathbb{E}\Big[\sup_{t-\delta
\leq l\leq \lambda -\delta }\big|\int_{l}^{\lambda -\delta }\big[\partial
_{t}b(\l;r)dr+\partial _{t}\sigma (\l;r)dW_{r}\big]\big|^{2p}%
\Big]d\lambda   \notag \\
&\leq &C_{p}\delta ^{2p-1}\int_{t}^{t+\delta }\mathbb{E}\Big[\big(\int_{t-\delta }^{\lambda -\delta }[|\partial _{t}b (\lambda; r)|^{2}+|\partial _{t}\sigma (\lambda; r)|^{2}]dr\big)^{p}\Big]d\lambda   \notag \\
&\leq &C_{p}\delta ^{2p-1}\int_{t}^{t+\delta }\big(
\int_{t-\delta }^{\lambda -\delta }(\lambda -r)^{2H-3}dr\big)^{p}d\lambda   \notag \\
&=&C_{p}\delta ^{2p-1}\int_{t}^{t+\delta }\big(\delta
^{2H-2}-(\delta +\lambda -t)^{2H-2}\big)^{p}d\lambda.
\eeaa
By changing variable, this implies
\bea
\mathbb{E}\Big[\sup_{t\leq s\leq T}|\Theta _{s}^{\delta ,t}-\Theta
_{s}^{t}|^{2p}\Big]  &\leq &C_{p}\delta ^{2pH}\int_{0}^{1}\big(1-(1+\lambda )^{2H-2}\big)^{p}d\lambda =C_{p}\delta ^{2pH}.  \label{Thdconv}
\eea

We next estimate the difference of $X$.  For any $t<t^{\prime }$, note that 
\begin{eqnarray}
&\dis \Big|(X_{t}^{\delta }-X_{t})-(X_{t^{\prime }}^{\delta }-X_{t^{\prime
}})\Big|\leq I_{1}+I_{2},\quad \mbox{where} & \label{I=I12}\\
&\dis I_{1}:= \Big|\int_{t}^{t^{\prime }}[b^{\delta }(t';r)-b(t^{\prime
}; r)]dr+\int_{t}^{t^{\prime }}[\sigma ^{\delta }(t';r)-\sigma
(t';r)]dW_{r}\Big|,&\nonumber \\
&\dis I_{2}:= \Big|\int_{0}^{t}[b^{\delta }(t;r)-b(t;r)-b^{\delta }(t^{\prime
}; r)+b(t';r)]dr\qq\qq\qq&\nonumber \\
&\dis \qq\qq\qq +\int_{0}^{t}[\sigma ^{\delta }(t;r)-\sigma (t;r)-\sigma ^{\delta
}(t';r)+\sigma (t';r)]dW_{r}\Big|.&\nonumber
\end{eqnarray}%
Denote $\delta ^{\prime }:=t^{\prime }-t$. For any $p\geq 1$, by the Burkholder-Davis-Gundy inequality and \reff{KH},
\begin{eqnarray}
\mathbb{E}[I_{1}^{2p}] 
&\leq & C_{p}\mathbb{E}\Big[\Big(\int_{t}^{t^{\prime }} [|b^{\delta }(t';s)-b(t';s)|^2+|\sigma ^{\delta }(t';s)-\sigma (t^{\prime}; s)|^{2}]ds\Big)^{p}\Big]\nonumber \\
&=&C_{p}\mathbb{E}\Big[\Big(\int_{t\vee
(t^{\prime }-\delta )}^{t^{\prime }}[|b(s+\delta; s)-b(t';s)|^2+|\sigma (s+\delta; s)-\sigma (t^{\prime
}; s)|^{2}]ds\Big)^{p}\Big]\nonumber\\
&\le&C_{p}\mathbb{E}\Big[\Big(\int_{t\vee (t^{\prime }-\delta )}^{t^{\prime
}} \big[\int_{t'}^{s+\d} [|\pa_t b(r; s)| + |\pa_t \si(r; s)|]dr\big]^2 ds \Big)^{p}\Big]\nonumber \\
&\leq &C_{p}\Big(\int_{t\vee (t^{\prime }-\delta )}^{t^{\prime
}}|\int_{t^{\prime }}^{s+\delta }(r-s)^{H-{\frac{3}{2}}}dr|^{2}ds\Big)^{p} =C_{p}\Big(\int_{0}^{\delta \wedge \delta ^{\prime }}[r^{H-{\frac{1}{2}}%
}-\delta ^{H-{\frac{1}{2}}}]^{2}dr\Big)^{p} \nonumber\\
&\leq &C_{p}\Big(\int_{0}^{\delta \wedge \delta ^{\prime }}r^{2H-1}dr\Big)%
^{p}\leq C_{p}(\delta \wedge \delta ^{\prime })^{2pH}. \label{I1est}
\end{eqnarray}
To estimate $I_2$, when $\delta ^{\prime }\geq \delta $,  by \reff{Thdconv} we have 
\bea
\dbE[|I_{2}|^{2p}] &=&  \dbE\Big[\Big|\int_{0}^{t}[b^{\delta }(t;s)-b(t;s)]ds+\int_{0}^{t}[\sigma
^{\delta }(t;s)-\sigma (t;s)]dW_{s}\Big|^{2p}\Big]\nonumber\\
&=& \dbE\Big[|\Theta _{t}^{\delta ,t}-\Theta
_{t}|^{2p}\Big]  \leq C_{p}\delta ^{2pH}. \label{I2est1}
\eea
When $\delta ^{\prime }<\delta $, one can check straightforwardly that 
\begin{equation*}
I_{2}=\Big|\int_{t-\delta }^{t}[b(t^{\prime }\wedge (s+\delta
); s)-b(t;s)]ds+\int_{t-\delta }^{t}[\sigma (t^{\prime }\wedge (s+\delta
); s)-\sigma (t;s)]dW_{s}\Big|.
\end{equation*}%
Then, again by the Burkholder-Davis-Gundy inequality and \reff{KH}, 
\bea
\mathbb{E}[I_{2}^{2p}] &\leq & 
C_{p}\Big(\int_{t-\delta }^{t}\big[\int_{t}^{t^{\prime }\wedge
(s+\delta )}(r-s)^{H-{\frac{3}{2}}}dr\big]^{2}ds\Big)^{p}  \nonumber\\
&\leq &C_{p}\Big(\int_{0}^\d{{\big[}}r^{H-{\frac{1}{2}}}-[(r+\delta ^{\prime
})\wedge \delta ]^{H-{\frac{1}{2}}}\big]^{2}dr\Big)^{p} \nonumber\\
&=&C_{p}\Big(\int_{0}^{\delta -\delta ^{\prime }}\big[r^{H-{\frac{1}{2}}%
}-(r+\delta ^{\prime })^{H-{\frac{1}{2}}}\big]^{2}dr\Big)^{p}+C_{p}\Big(%
\int_{\delta -\delta ^{\prime }}^\d{{\big[}}r^{H-{\frac{1}{2}}}-\delta ^{H-{%
\frac{1}{2}}}\big]^{2}dr\Big)^{p} \nonumber\\
&\leq &C_{p}(\delta ^{\prime })^{2pH}\Big(\int_{0}^{\infty }\big[r^{H-{\frac{%
1}{2}}}-(r+1)^{H-{\frac{1}{2}}}\big]^{2}dr\Big)^{p}+C_{p}\Big(\int_{\delta
-\delta ^{\prime }}^\d {r}{}^{2H-1}dr\Big)^{p}.\nonumber \\
&\leq &C_{p}(\delta ^{\prime })^{2pH}+C_{p}\big[\delta ^{2H}-(\delta -\delta^{\prime })^{2H}\big]^{p} ~\le~ C_{p}(\delta ^{\prime })^{2pH},  \label{I2est2}
\eea
where the last inequality thanks to the assumption that $H<{1\over 2}$.
 Plug \reff{I1est}, \reff{I2est1}, and \reff{I2est2} into \reff{I=I12},  we get
\begin{equation*}
\mathbb{E}\Big[\big|(X_{t}^{\delta }-X_{t})-(X_{t^{\prime }}^{\delta
}-X_{t^{\prime }})\big|^{2p}\Big]\leq C_{p}(\delta \wedge
\delta ^{\prime })^{2pH}=C_{p}\delta ^{pH}(\delta ^{\prime
})^{pH}=C_{p}\delta ^{pH}|t^{\prime }-t|^{pH}.
\end{equation*}%
Then by Lemma \ref{lem-kolmogorov} we see that 
$\mathbb{E}\big[\Vert X^{\delta }-X\Vert _{T}^{2p}\big]\leq C_{p}\delta ^{pH}$.
This, together with \textrm{(\ref{Thdconv})}, proves \textrm{(\ref{Xdconv})}.

\textit{Step 3.} We now prove \textrm{(\ref{FIto})}. We first note that $%
u(t,X^{\delta }\otimes _{t}\Theta ^{\delta ,t})$ falls short of satisfying
the conditions in Theorem \ref{thm-FIto1}. In fact, $X^{\delta }$ is not a
solution to the SDE \textrm{(\ref{X})} with coefficients $(b^{\delta },\sigma
^{\delta })$, and $b^{\d}$, $\sigma ^{\delta }$ are not
differentiable at $t=s+\delta $. However, by checking the arguments of the
proof of Theorem \ref{thm-FIto1} one can see that these differences do not cause any trouble and thus the conclusion still holds true. Therefore, $u(t,X^{\delta }\otimes _{t}\Theta ^{\delta
,t})$ satisfies \textrm{(\ref{FIto})}: 
\begin{eqnarray}
&&du(t,X^{\delta }\otimes _{t}\Theta ^{\delta ,t})=\partial
_{t}u(t,X^{\delta }\otimes _{t}\Theta ^{\delta ,t})dt+{\frac{1}{2}}%
\langle \partial _{\omega \omega }^{2}u(t,X^{\delta }\otimes _{t}\Theta
^{\delta ,t}),(\sigma ^{\delta ,t,X},\sigma ^{\delta ,t,X})\rangle dt  \notag
\\
&&\qquad +\langle \partial _{\o }u(t,X^{\delta }\otimes _{t}\Theta
^{\delta ,t}),b^{\delta ,t,X}\rangle dt+\langle \partial _{\o %
}u(t,X^{\delta }\otimes _{t}\Theta ^{\delta ,t}),\sigma ^{\delta
,t,X}\rangle dW_{t},\quad \mathbb{P}\mbox{-a.s.}  \label{FItod}
\end{eqnarray}%
Then, by \textrm{(\ref{Xdconv})} and by the continuity of $u$ and $\partial
_{t}u$, we have 
\begin{equation}
\lim_{\delta \rightarrow 0}\mathbb{E}\Big[|u(t,X^{\delta }\otimes
_{t}\Theta ^{\delta ,t})-u(t,X\otimes _{t}\Theta ^{t})|+|\partial
_{t}u(t,X^{\delta }\otimes _{t}\Theta ^{\delta ,t})-\partial
_{t}u(t,X\otimes _{t}\Theta ^{t})|\Big]=0.  \label{udconv}
\end{equation}%
Moreover, recall the notations in \textrm{(\ref{dn})} and assume $\delta_{n+1}\leq
\delta <\delta_{n}$. Set $m=1$ in \textrm{(\ref{fn})}: 
\begin{equation*}
\mathbf{1}_{[t_{1},T]}=[1-\psi _{1}]\mathbf{1}_{[t_{1},T]}+\sum_{k=1}^{n}\psi
_{k}+[1-\psi_{n}]\mathbf{1}_{[t,t_{n}]}.
\end{equation*}%
Then, denoting $\tilde{X}:=X\otimes _{t}\Theta ^{t},\tilde{X}%
^\d:=X^{\delta }\otimes _{t}\Theta ^{\delta ,t}$, by \textrm{(\ref%
{pauKconv})} and \reff{vanish2} we have, for $\varphi =b,\sigma $, 
\begin{eqnarray*}
&&\Big|\langle \partial _{\o }u(t,\tilde{X}^{\delta }),\varphi ^{\delta
,t,X}\rangle -\langle \partial _{\o }u(t,\tilde{X}),\varphi ^{t,X}\rangle %
\Big| \\
&\leq &\Big|\langle \partial _{\o }u(t,\tilde{X}^{\delta })-\partial _{\o %
}u(t,\tilde{X}),\varphi ^{\delta ,t,X}\rangle \Big|+\Big|\langle \partial _{%
\o }u(t,\tilde{X}),\varphi ^{\delta ,t,X}-\varphi ^{t,X}\rangle \Big| \\
&\leq &\Big|\langle \partial _{\o }u(t,\tilde{X}^{\delta })-\partial _{\o %
}u(t,\tilde{X}),\big[[1-\psi _{1}]\mathbf{1}_{[t_{1},T]}+\sum_{k=1}^{n}%
\psi _{k}+[1-\psi_{n}]\mathbf{1}_{[t,t_{n}]}\big]\varphi ^{\delta
,t,X}\rangle \Big| \\
&&+C[1+\Vert \tilde X\Vert _{T}^\k]\delta ^{\b} \\
&\leq &C[1+\|\tilde X\|_T^\k+ \Vert\tilde X^\d\Vert _{T}^\k]\Big[\rho(\Vert \tilde{X}^{\delta }-\tilde{X}\Vert _{T})\big[1+\sum_{k=1}^{n}\delta _{k}^\b+\delta_{n}^\b\big]+\delta ^{\b} \Big]\\
&\leq &C[1+\|\tilde X\|_T^\k+ \Vert\tilde X^\d\Vert _{T}^\k] \Big[\rho(\Vert \tilde{X}^{\delta }-\tilde{X}\Vert _{T}) + \delta ^{\b}\Big].
\end{eqnarray*}%
Similarly,  by \textrm{(\ref{pauKconv})} and \reff{vanish2} we have
\begin{eqnarray*}
&&\Big|\langle \partial _{\omega \omega }^{2}u(t,\tilde{X}^{\delta
}),(\sigma ^{\delta ,t,X},\sigma ^{\delta ,t,X})\rangle -\langle \partial
_{\omega \omega }^{2}u(t,\tilde{X}),(\sigma ^{t,X},\sigma ^{t,X})\rangle %
\Big| \\
&\leq &\Big|\langle \partial _{\omega \omega }^{2}u(t,\tilde{X}^{\delta
})-\partial _{\omega \omega }^{2}u(t,\tilde{X}),(\sigma ^{\delta
,t,X},\sigma ^{\delta ,t,X})\rangle \Big| \\
&&+\Big|\langle \partial _{\omega \omega }^{2}u(t,\tilde{X}),(\sigma
^{\delta ,t,X},\sigma ^{\delta ,t,X})\rangle -\langle \partial _{\omega
\omega }^{2}u(t,\tilde{X}),(\sigma ^{t,X},\sigma ^{t,X})\rangle \Big| \\
&\leq &\Big|\langle \partial^2 _{\o\o }u(t,\tilde{X}^{\delta })-\partial^2 _{\o\o %
}u(t,\tilde{X}),\Big(\sigma ^{\delta ,t,X},\big[[1-\psi _{1}]\mathbf{1}%
_{[t_{1},T]}+\sum_{k=1}^{n}\psi _{k}+[1-\psi_{n}]\mathbf{1}%
_{[t,t_{n}]}\big]\sigma ^{\delta ,t,X}\Big)\rangle \Big| \\
&&+C[1+\Vert X\Vert _{T}^\k]\delta ^{2\b} \\
&\leq &C[1+\|\tilde X\|_T^\k+ \Vert\tilde X^\d\Vert _{T}^\k]\Big[\rho(\Vert \tilde{X}^{\d}-\tilde{X}\Vert
_{T})\big[1+\sum_{k=1}^{n}\delta _{k}^{2\b}+\delta
_{n}^{2\b)}\big]+\delta ^{2\b}\Big] \\
&\leq &C[1+\|\tilde X\|_T^\k+ \Vert\tilde X^\d\Vert _{T}^\k]\Big[\rho(\Vert \tilde{X}^{{\d}}-\tilde{X}%
\Vert _{T})+\delta ^{2\b}\Big].
\end{eqnarray*}%
Put together, we derive from \textrm{(\ref{Xdconv})} that 
\begin{eqnarray*}
&\dis \lim_{\d\to 0}\mathbb{E}\Big[\Big|\langle \partial _{\o }u(t,\tilde{X}^{\delta
}),\varphi ^{\delta ,t,X}\rangle -\langle \partial _{\o }u(t,\tilde{X}%
),\varphi ^{t,X}\rangle \Big|^{2}\Big] =0,& \\ 
&\dis \lim_{\d\to 0}\mathbb{E}\Big[\Big|\langle \partial _{\omega \omega }^{2}u(t,\tilde{X}%
^{\delta }),(\sigma ^{\delta ,t,X},\sigma ^{\delta ,t,X})\rangle -\langle
\partial _{\omega \omega }^{2}u(t,\tilde{X}),(\sigma ^{t,X},\sigma
^{t,X})\rangle \Big|\Big]=0.& 
\end{eqnarray*}%
Plug this and \textrm{(\ref{udconv})} into \textrm{(\ref{FItod})}, we obtain 
\textrm{(\ref{FIto})} in this singular case. \hfill \vrule width.25cm
height.25cm depth0cm\smallskip 

\section{Path dependent PDEs and Feynman-Kac formulae}
\label{sect-PPDE} \setcounter{equation}{0} 

\subsection{The linear case}
\label{sect-linear}
We first apply Theorem \ref{thm-FIto2}  to  the
path-dependent example of Section \ref{ss-path}. This works because the
kernel of fBm satisfies all the hypotheses on $b$ and $\sigma $ in Theorem %
\ref{thm-FIto2}, including the ones where $u$ is weakly dependent on paths
near its time diagonal. In fact, we find that we may choose $\alpha =1/2$ in
Definition \ref{defn-C12a}.

\begin{thm}
\label{thm-rep} Consider the setting in Subsection \ref{ss-path}, and
denote 
\begin{equation}
u(t,\omega ):=\mathbb{E}\Big[g\big(\omega _{T}+\int_{t}^{T}K(T,r)dW_{r}\big)%
+\int_{t}^{T}f\big(s,\omega _{s}+\int_{t}^{s}K(s,r)dW_{r}\big)ds\Big].
\label{linearu}
\end{equation}%
Assume $f$ is continuous in $t$; and for $\f = g, f(t,\cd)$, $\f\in C^{2
}(\mathbb{R})$ such that all the derivatives have polynomial growth with $|\f''(t,x) - \f''(t,\tilde x)|\le C[1+|x|^\kappa + |\tilde x|^\kappa] \rho(|x-\tilde x|)$.  Then

(i) $u$ evaluated at $B^{H}\otimes _{t}\Theta ^{t}$ coincides with the
conditional expectation: 
\begin{equation}
Y_{t}:=\mathbb{E}\Big[g(B_{T}^{H})+\int_{t}^{T}f(s,B_{s}^{H})ds\Big|\mathcal{%
F}_{t}\Big]=u(t,B^{H}\otimes _{t}\Theta ^{t}).  \label{linearY}
\end{equation}

(ii) $u\in C_{+}^{1,2}(\overline{\Lambda })$ with path derivatives: 
\begin{eqnarray}
&&\langle \partial _{\o }u(t,\omega ),\eta \rangle =\mathbb{E}\Big[g^{\prime
}\big(\omega _{T}+\int_{t}^{T}K(T,r)dW_{r}\big)\eta
_{T}+\int_{t}^{T}f^{\prime }\big(s,\omega _{s}+\int_{t}^{s}K(s,r)dW_{r}\big)%
\eta _{s}ds\Big],  \notag \\
&&\langle \partial _{\omega \omega }^{2}u(t,\omega ),(\eta _{1},\eta
_{2})\rangle =\mathbb{E}\Big[g^{\prime \prime }\big(\omega
_{T}+\int_{t}^{T}K(T,r)dW_{r}\big)\eta _{1}(T)\eta _{2}(T)  \label{pau-rep}
\\
&&\qquad \qquad \qquad \qquad \qquad +\int_{t}^{T}f^{\prime \prime }\big(%
s,\omega _{s}+\int_{t}^{s}K(s,r)dW_{r}\big)\eta _{1}(s)\eta _{2}(s)ds\Big] 
\notag
\end{eqnarray}

(iii) $u$  vanishes diagonally with rate $\alpha={\frac{1}{2}}$, in the sense of Definition \ref{defn-C12a}.
Consequently, the functional It\^{o} formula \textrm{(\ref{FIto})} holds
true for all $H\in (0,1)$.

(iv) $u$ is a classical solution to the following linear PPDE which,
together with \textrm{(\ref{pau-rep})}, provides a representation for $%
\partial _{t}u$: 
\begin{equation}
\partial _{t}u(t,\omega )+{\frac{1}{2}}\langle \partial _{\omega \omega
}^{2}u(t,\omega ),(K^{t},K^{t})\rangle +f(t,\omega _{t})=0,~(t,\omega )\in
\Lambda ;\quad u(T,\omega )=g(\omega _{T}).  \label{linearPPDE}
\end{equation}
\end{thm}

{\noindent \textbf{Proof\quad }} We shall only prove the irregular case $H<{%
\frac{1}{2}}$. The regular case $H\geq {\frac{1}{2}}$ follows by similar but
easier arguments.

First, \textrm{(\ref{linearY})} follows directly from the arguments in
Subsection \ref{ss-state}.

Next, applying \textrm{(\ref{pathu})} and \textrm{(\ref{path2u})} on \textrm{%
(\ref{linearu})}, one may easily verify \textrm{(\ref{pau-rep})}, and that $\partial_{\omega} u$, $\partial^2_{\omega\omega} u$ are
 continuous, have polynomial growth, and satisfy \textrm{(\ref{paooucont})}, \textrm{(\ref{vanish}%
)}, and \textrm{(\ref{vanish2})} with $\alpha={\frac{1}{2}}$.

Moreover, denote $\overline{\sigma }^{2}(s,t):=\int_{t}^{s}K^{2}(s,r)dr<%
\infty $, $0\leq t<s\leq T$. Then $\int_{t}^{s}K(s,r)dW_{r}$ has 
distribution Normal$(0, \overline{\sigma }^{2}(s,t))$. Thus,
denoting by $\mathcal{N}$ a standard normal distribution, 
\begin{equation*}
u(t,\omega )=\mathbb{E}\Big[g\big(\omega _{T}+\overline{\sigma }(T,t)%
\mathcal{N})\big)+\int_{t}^{T}f\big(s,\omega _{s}+\overline{\sigma }(s,t)%
\mathcal{N}\big)ds\Big].
\end{equation*}%
Note that $\partial _{t}\overline{\sigma }(s,t)=-{\frac{K^{2}(s,t)}{2%
\overline{\sigma }(s,t)}}$. Then one can easily see that 
\begin{eqnarray}
\partial _{t}u(t,\omega ) &=&-{\frac{1}{2}}\mathbb{E}\Big[\big[g^{\prime }%
\big(\omega _{T}+\overline{\sigma }(T,t)\mathcal{N}\big)-g^{\prime }(\omega
_{T})\big]{\frac{K^{2}(T,t)}{\overline{\sigma }(T,t)}}\mathcal{N}  \notag \\
&&\qquad +\int_{t}^{T}\big[f^{\prime }\big(s,\omega _{s}+\overline{\sigma }%
(s,t)\mathcal{N}\big)-f^{\prime }(s,\omega _{s})\big]{\frac{K^{2}(s,t)}{%
\overline{\sigma }(s,t)}}\mathcal{N}ds\Big]  \notag \\
&=&-{\frac{1}{2}}\int_{0}^{1}\mathbb{E}\Big[g^{\prime \prime }\big(\omega
_{T}+\lambda \overline{\sigma }(T,t)\mathcal{N}\big)K^{2}(T,t)\mathcal{N}^{2}
\label{patu-rep} \\
&&\qquad +\int_{t}^{T}f^{\prime \prime }\big(s,\omega _{s}+\lambda \overline{%
\sigma }(s,t)\mathcal{N}\big)K^{2}(s,t)\mathcal{N}^{2}ds\Big]d\lambda . 
\notag
\end{eqnarray}%
where, to justify the integrability in the right side above, we note that 
\begin{equation*}
|\partial _{t}u(t,\omega )|\leq C\Big[K^{2}(T,t)+\int_{t}^{T}K^{2}(s,t)ds%
\Big][1+\|\o\|_T^\kappa] \leq C[1+\|\o\|_T^\kappa].
\end{equation*}%
This means that $\partial _{t}u$ exists and has polynomial growth. By \textrm{(\ref%
{patu-rep})} one can also see that $\partial _{t}u$ is continuous. Then $%
u\in C_+^{1,2}(\overline{\Lambda })$.

Finally, note that 
\begin{equation*}
u(t,B^{H}\otimes _{t}\Theta ^{t})+\int_{0}^{t}f(s,B_{s}^{H})ds=\mathbb{E}%
\Big[g(B_{T}^{H})+\int_{0}^{T}f(s,B_{s}^{H})ds\Big|\mathcal{F}_{t}\Big]
\end{equation*}%
is a martingale. Applying the functional It\^{o} formula \textrm{(\ref{FIto})%
} on $u(t,B^{H}\otimes _{t}\Theta ^{t})$, we see that \textrm{(\ref%
{linearPPDE})} holds on $B^{H}\otimes _{t}\Theta ^{t}$, $\mathbb{P}$-a.s. In
particular, \textrm{(\ref{linearPPDE})} holds at $(0,0)$. Given $(t,\omega
)\in \overline{\Lambda }$, apply the same arguments on the system starting
with $(t,\omega )$, in the spirit of the proof of Proposition  \ref{prop-path} (Step 2 of that joint proof) ,
we can see that \textrm{(\ref{linearPPDE})} holds at $(t,\omega )$ as well.
\hfill \vrule width.25cm height.25cm depth0cm\smallskip 

\subsection{The semilinear case}

In this subsection we consider the
following BSDE: 
\begin{equation}
Y_{t}=g(X_{\cd})+\int_{t}^{T}f(s,X_{\cd},Y_{s},Z_{s})ds-%
\int_{t}^{T}Z_{s}dW_{s};\quad 0\leq t\leq T,  \label{BSDE}
\end{equation}%
where $f=f(s,X_{s\wedge \cdot },y,z)$ is adapted. We emphasize again that,
as we saw in the previous section, even if the coefficients $b,\sigma ,f,g$
are state dependent (namely depending only on $X_{s}$ at time $s$), the BSDE is not
Markovian as soon as $X$ is not a Markov process. When $X$ is a strong
solution to SDE \textrm{(\ref{X})}, namely $X$ is $\mathbb{F}^{W}$%
-progressively measurable, then it follows from the seminal work Pardoux \&
Peng \cite{PP} that the above BSDE is well-posed, provided $f$ is uniformly
Lipschitz-continuous in $(y,z)$. When $X$ is a weak solution and no strong
solution exists, then typically one needs to introduce an orthogonal
martingale term in the BSDE \textrm{(\ref{BSDE})}. We avoid this situation
in the sequel, though the uniqueness of a strong solution to (\ref{BSDE}) is
not a requirement.

This BSDE is closely related to the following semilinear PPDE, where the
notation is that of Section \ref{Derivatives}: 
\begin{equation}
\left. 
\begin{array}{lll}
\partial _{t}u(t,\omega )+{\frac{1}{2}}\langle \partial _{\omega \omega
}^{2}u(t,\omega ),(\sigma ^{t,\omega },\sigma ^{t,\omega })\rangle +\langle
\partial _{\o }u(t,\omega ),b^{t,\omega }\rangle  &  &  \\ 
\qquad \qquad +f\big(t,\omega ,u(t,\omega ),\langle \partial _{\o %
}u(t,\omega ),\sigma ^{t,\omega }\rangle \big)=0,\quad (t,\omega )\in
\Lambda , &  &  \\ 
u(T,\omega )=g(\omega ). &  & 
\end{array}%
\right.   \label{semilinearPPDE}
\end{equation}%
We have the following Feynman-Kac formula.

\begin{thm}
\label{thm-BSDE} Let Assumption \ref{assum-X} hold and assume the semilinear
PPDE (\ref{semilinearPPDE}) has a classical solution $u\in
C_{+}^{1,2}(\Lambda )$. Assume further that either Assumption \ref%
{assum-regular} holds or Assumption \ref{assum-irregular} 
holds and $u\in C_{+,\a}^{1,2}(\Lambda )$ for some $\a > {1\over 2} -H$. Then the following provides a strong solution to the BSDE \textrm{(\ref%
{BSDE})}: 
\begin{equation}
Y_{t}:=u(t,X\otimes _{t}\Theta ^{t}),\quad Z_{t}:=\langle \partial _{\o %
}u(t,X\otimes _{t}\Theta ^{t}),\sigma ^{t,X}\rangle . \label{FK}
\end{equation}
\end{thm}

{\noindent \textbf{Proof\quad }} By our assumptions $u(t,X\otimes _{t}\Theta
^{t})$ satisfies the functional It\^{o} formula \textrm{(\ref{FIto})}. Then it
is straightforward to verify that the process $(Y,Z)$ defined by \textrm{(%
\ref{FK})} satisfies \textrm{(\ref{BSDE})}. \hfill \vrule width.25cm
height.25cm depth0cm\smallskip 

We note that, when the PPDE has a classical solution, \textrm{(\ref{FK})}
provides a solution to the BSDE \textrm{(\ref{BSDE})} even if $X$ is a weak
solution to \textrm{(\ref{X})}. However, in this case typically $(Y,Z)$ are
also not $\mathbb{F}^{W}$-progressively measurable.

We now proceed in the opposite direction, namely to provide a representation
for the solution of PPDE \textrm{(\ref{semilinearPPDE})} through the BSDE 
\textrm{(\ref{BSDE})}. For each $(t,\omega )\in \overline{\Lambda }$, define 
\begin{eqnarray}
\displaystyle &u(t,\omega ):=Y_{t}^{t,\omega },\quad \mbox{where} & \label{semilinearu}\\ 
&\left. 
\begin{array}{c}
\displaystyle X_{s}^{t,\omega }=\omega _{s}+\int_{t}^{s}b(s;r,\omega \otimes
_{t}X^{t,\omega })dr+\int_{t}^{s}\sigma (s;r,\omega \otimes _{t}X^{t,\omega
})dW_{r} \\ 
\displaystyle Y_{s}^{t,\omega }=g(\omega \otimes _{t}X^{t,\omega
})+\int_{s}^{T}f(r,\omega \otimes _{t}X^{t,\omega },Y_{r}^{t,\omega
},Z_{r}^{t,\omega })dr-\int_{s}^{T}Z_{r}^{t,\omega }dW_{r},
\end{array}%
\right.  t\leq s\leq T.& \nonumber
\end{eqnarray}%
Here we are assuming that the FBSDE in (\ref{semilinearu}) has a unique strong solution for all $(t,\omega )\in \overline{\Lambda }
$. With that assumption, for fixed $(t,\omega )\in \overline{\Lambda }$, the
pair of processes $\left( Y^{t,\omega },Z^{t,\omega }\right) $ is given
unambiguously by the FBSDE, and $u$ in (\ref{semilinearu}) is well defined
on $\overline{\Lambda }$. We avoid further technical discussion, stating the
representation result with comments only.%

\begin{rem}
\label{rem-nonlinear} \textrm{(i) Provided appropriate conditions on the
coefficients of (\ref{X}) and (\ref{BSDE}), one can show that the $u$
defined by (\ref{semilinearu}) is indeed smooth and 
is the classical solution to PPDE (\ref{semilinearPPDE}). See Peng
{\rm \&} Wang \cite{PW} for a related result in the Brownian motion framework.}

\textrm{(ii) Our BSDE (\ref{BSDE}) is time consistent: } the coefficient $f$%
\textrm{\  depends only on one time variable and $g$ is independent of  time. 
We refer to Yong 
\cite{Yong} and the references therein for  Volterra
type BSDEs. \hfill \vrule width.25cm height.25cm depth0cm\smallskip }
\end{rem}

\subsection{A strategy for control problems}

The framework of the previous subsection applies directly, as a slight
extension, if there is control involved. Formally one can easily write down
the path-dependent Hamilton-Jacobi-Bellman (HJB) equation. More precisely,
let $\mathcal{A}$ be an appropriate set of admissible controls taking values
in certain set $A$;
$X$ solve a controlled Volterra SDE;  and  $Y$ solve a controlled BSDE. We define the value function $u$ for the control problem as follows.
\begin{eqnarray}
&\displaystyle  u(t,\omega ):=\sup_{a\in \mathcal{A}}Y_{t}^{t,\omega ,a},\quad %
\mbox{where} & \label{nonlinearu}\\ 
&\left. 
\begin{array}{c}
\displaystyle X_{s}^{t,\omega ,a}=\omega _{s}+\int_{t}^{s}b(s;r,\omega
\otimes _{t}X^{t,\omega ,a},a_{s})dr+\int_{t}^{s}\sigma (s;r,\omega \otimes
_{t}X^{t,\omega ,a},a_{s})dW_{r} \\ 
\displaystyle Y_{s}^{t,\omega ,a}=g(\omega \otimes _{t}X^{t,\omega
,a})+\int_{s}^{T}f(r,\omega \otimes _{t}X^{t,\omega ,a},Y_{r}^{t,\omega
,a},Z_{r}^{t,\omega ,a},a_{s})dr-\int_{s}^{T}Z_{r}^{t,\omega ,a,}dW_{r}.%
\end{array}%
\right. &\nonumber
\end{eqnarray}%
Then formally $u$ should satisfy the following path dependent HJB equation:
\begin{equation}
\left. 
\begin{array}{lll}
\partial _{t}u(t,\omega )+\sup_{a\in A}\Big[{\frac{1}{2}}\langle \partial
_{\omega \omega }^{2}u(t,\omega ),(\sigma ^{t,\omega ,a},\sigma ^{t,\omega
,a})\rangle +\langle \partial _{\o }u(t,\omega ),b^{t,\omega ,a}\rangle  & 
&  \\ 
\qquad \qquad +f\big(t,\omega ,u(t,\omega ),\langle \partial _{\o %
}u(t,\omega ),\sigma ^{t,\omega ,a}\rangle ,a_{t}\big)\Big]=0,\quad
(t,\omega )\in \Lambda , &  &   
\end{array}%
\right.   \label{HJB}
\end{equation}%
with terminal condition $u(T,\omega )=g(\omega )$. Here,  for $\varphi =b,\sigma $, $\varphi _{s}^{t,\omega ,a} := \varphi
(s; t,\omega ,a)$. See Fouque \& Hu \cite{FH} for an application in this
direction.

When the path dependent HJB equation \textrm{(\ref{HJB})} has a classical
solution $u\in C_{+}^{1,2}(\Lambda )$ or when the value function $u$ defined by \textrm{(%
\ref{nonlinearu})} is indeed in $C_{+}^{1,2}(\Lambda )$ (or $u\in C_{+,\a}^{1,2}(\Lambda )$ for some appropriate $\a$ in the singular case) , it is
not difficult to prove that they are equal, as in the standard verification
theorem. However, in general it is difficult to expect a classical solution
for such a control problem, because of the  path dependence. We shall study
viscosity solutions, in the spirit of Ekren, Touzi \& Zhang \cite{ETZ1, ETZ2}%
, for these fully nonlinear PPDEs in our future research.

\section{An application to finance}
\label{sect-Appli}
\setcounter{equation}{0}

In the reference El Euch \& Rosenbaum \cite{ER}, the authors work with the so-called rough Heston
model, whose wellposedness was established in their previous publication \cite{ER0}.
In \cite{ER}, they show that options
on equities given by this model can be hedged if one assumes that the
volatility is observed. In fact, for an option on a given equity, they argue that, since  the
forward variance can be replicated in the market using liquid instruments, then all that is required for hedging purposes, is observation of that forward  variance and the equity's spot price. We will describe their model more precisely, how it fits in ours, and what more pricing and hedging questions can be reached in ours.

Recall Section \ref{Fintro} and in particular \reff{BS}. Consider the following rough Heston model: 
\bea
\left.\ba{lll} 
\dis S_t   =   S_0 + \int_0^t S_r \sqrt{V_{r}} ~\big[\sqrt{1-\rho^2} dW^1_r + \rho dW^2_r\big],   \\
\dis V_{t}  = V_{0}+\frac{1}{\Gamma \left( H+{1\over 2}\right) }\int_{0}^{t}\left( t-r\right) ^{H-{1\over 2}} \big[\lambda [ \theta -V_r] dr + \nu \sqrt{V_r} dW^2_{r}\big].
\ea\right.  \label{Heston}
\eea
Here $W = (W^1, W^2)^\top$ is a two dimensional  Brownian motion, $\rho \in [0, 1]$ is a correlation parameter, 
$\theta $ is a mean-reversion level, 
$\lambda $ is a mean-reversion rate, $\nu $ is a noise intensity, and the roughness parameter $H$ is typically in $(0, {1\over 2})$. We
leave aside the question of whether any of these parameters can be estimated
or calibrated from the data. We note instead that the term $\left(
t-r\right) ^{H-{1\over 2}}$ is similar to the kernel $K$ of fBm (it is in fact
identical to the kernel of the so-called Riemann-Liouville fBm).   By \cite{ER0}, the SDE \reff{Heston} has a unique weak solution  $X := (S, V)^\top$. 
The paper \cite{ER} asks the question of how to compute the conditional expectation of a
non-path-dependent contingent claim at any time prior to maturity:
\begin{equation}
Y_{t} :=C_{t}:=\mathbb{E}\left[ g\left( S_{T}\right) |\mathcal{F}_{t}\right]
\label{Ct}
\end{equation}
for some deterministic contract function $g$.  They express $C_t$ as a function of $S_t$ as well as the so called forward variance:
\bea
\label{FV}
\hat \Th^t_s := \dbE[ V_s |\cF_t],\q 0\le t \le s \le T.
\eea
Note that both the  above forward variance and the forward volatility $\mathbb{E}[\sqrt{V_{s}}|\mathcal{F}_{t}]$ introduced in
Section \ref{Fintro} are  financial indexes available in the market. 
A PPDE is derived for $C_t$ in this special case. 
Moreover, the forward variance can be approximated by using liquid variance swaps or vanilla options, and in this sense one may view the forward variance as a set of additional tradable assets. The main contribution of \cite{ER} is to provide a perfect hedge for the derivative $g(S_T)$ by using the stock  $S$ and the forward variance $\hat \Th$. The hedging portfolio relies on the Frechet derivative of $C_t$ and certain characteristic functions, which requires the special structure of \reff{Heston} and  that $C_T = g(S_T)$ is state dependent.

We now explain how our framework covers the above example and beyond. First note that, for $X = (S, V)^\top$,  \reff{Heston} is a Volterra SDE \reff{X} with 
\bea
\label{Hestoncoeff}
b(t; r, x_1, x_2)  = \left[ \ba{c}  \!\!\! 0 \\  \!\!\! {\l(t-r)^{H-{1\over 2}} [\th-x_2] \over \G(H+{1\over 2})}\ea  \!\!\! \right],  ~ \si(t;r, x_1, x_2) = \left[ \ba{lll}  \!\!\! \sqrt{1-\rho^2} x_1 \sqrt{x_2} &  \q\!\!\! \rho x_1 \sqrt{x_2}\\\qq  \!\!\! 0 & \!\!\! {\nu (t-r)^{H-{1\over 2}} \sqrt{x_2} \over \G(H+{1\over 2})}\ea  \!\!\! \right].
\eea
One may easily check that \reff{Heston} satisfies all the properties in Assumptions \ref{assum-X} and  \ref{assum-irregular}, needed
in Section \ref{singular} for $H\in (0,1/2)$, see Remark \ref{rem-Hestonintegrability} below. Note that the dynamics of $S$ is standard, without involving a two-time-variable kernel. While we may apply the results in previous sections directly on the two dimensional SDE \reff{Heston}, for simplicity we restrict the path dependence only to the dynamics of $V$. Therefore, recall \reff{Th}, we denote
\bea
\label{Heston-Th}
\Th^t_s :=  V_{0}+\frac{1}{\Gamma \left( H+{1\over 2}\right) }\int_{0}^{t}\left( s-r\right) ^{H-{1\over 2}} \Big[\lambda [ \theta -V_r] dr + \nu \sqrt{V_r} dW^2_{r}\Big],  \q t<s.
\eea
By the special structure of the rough Heston model, we can actually see that 
\bea
\label{Heston-Cu}
C_t = u\big(t, S_t,   \Th^t_{[t, T]}\big).
\eea
In particular, the dependence of $C_t$ on $S$ is only via $S_t$ and its dependence on $V$ does not involve  $V_{[0, t)}$. Denote $u$ as $u(t,x,\o)$  and we shall assume $g$ is smooth which would imply the smoothness of $u$. Now following the arguments in Section \ref{sect-linear}, in particular noting that $C$ is a martingale,  we see that $u$ satisfies the following PPDE:  
\bea
\label{Heston-PPDE}
\left.\ba{lll}
\dis \pa_t u +  {\l [\th- \o_t] \over \G(H+{1\over 2})} \la \pa_\o u, a^t\ra  + {x^2 \o_t \over 2} \pa^2_{xx} u+ {\rho \nu x \o_t \over \G(H+{1\over 2})} \la \pa_{\o} (\pa_x u), a^t\ra \\
\dis\qq     + {\nu^2  \o_t \over  2\G(H+{1\over 2})} \la \pa^2_{\o\o} u, (a^t, a^t)\ra = 0,\qq\mbox{where}\q a^t_s := (s-t)^{H-{1\over 2}}.
\ea\right.
\eea
Moreover, by Theorem \ref{thm-FIto2}, we have (recalling $V_t = \Th^t_t$)
\bea
\label{Heston-Ito}
d C_t = \pa_x u \big(t, S_t,    \Th^t_{[t, T]}\big) dS_t +  {\nu \sqrt{V_t} \over \G(H+{1\over 2})} \Big\la \pa_\o u\big(t, S_t,   \Th^t_{[t, T]}\big) , a^t\Big\ra dW^2_t.
\eea
The first term in the right side above obviously provides the $\D$-hedging in terms of the stock $S$. Note further that $t \mapsto \Th^t_T$ is a semi-martingale, and we have
\beaa
{\nu \sqrt{V_t} \over \G(H+{1\over 2})} dW^2_t = (T-t)^{{1\over 2}-H} d \Th^t_T ~ - ~  {\l [\th- V_t]\over \G(H+{1\over 2})} dt.
\eeaa
Then
\bea
\label{Heston-Ito2}
\left.\ba{c}
\dis d C_t = \pa_x u \big(t, S_t,    \Th^t_{[t, T]}\big) dS_t +  (T-t)^{{1\over 2}-H}  \Big\la \pa_\o u\big(t, S_t,    \Th^t_{[t, T]}\big) , a^t\Big\ra d\Th^t_T\\
\dis -  {\l [\th- V_t]\over \G(H+{1\over 2})} \Big\la \pa_\o u\big(t, S_t,    \Th^t_{[t, T]}\big) , a^t\Big\ra  dt.
\ea\right.
\eea
That is, provided that we could  replicate $\Th^t_T$ using market instruments, which we will discuss in details below,  then we may (perfectly) hedge $g(S_T)$ as claimed in \cite{ER}.

We mention that our $\Th^t$ in \reff{Heston-Th} is different from the forward variance $\hat \Th^t$ in \reff{FV}. However, it can easily be replicated  by using $\hat \Th^t$, which can further be replicated (approximately) by variance swaps. Indeed, by \reff{Heston-Th} and  taking conditional expectation on the dynamics of $V$ in \reff{Heston}, we see that
\bea
\label{Heston-equivalent}
\hat \Th^t_s = \Th^t_s +  \frac{1}{\Gamma \left( H+{1\over 2}\right) } \int_t^s\left( s-r\right) ^{H-{1\over 2}} \lambda [ \theta - \hat \Th^t_r] dr,\q t\le s\le T.
\eea
 For any fixed $t$,  clearly $\Th^t_s$ is uniquely determined by $\{\hat\Th^t_r\}_{t\le r\le s}$:
\bea
\label{Heston-equivalent1}
\Th^t_s  = \hat \Th^t_s -  \frac{1}{\Gamma \left( H+{1\over 2}\right) } \int_t^s\left( s-r\right) ^{H-{1\over 2}} \lambda [ \theta - \hat \Th^t_r] dr.
\eea
In particular, this implies that, provided we observe the forward variance $\hat \Th^t_s$,  the process $\Th^t_s$ is also observable at $t$. Moreover, as a function of $t$, 
\beaa
d \Th^t_T =  d\hat \Th^t_T +\frac{1}{\Gamma \left( H+{1\over 2}\right) } \left( T-t\right) ^{H-{1\over 2}} \lambda [ \theta - \hat \Th^t_t] dt.
\eeaa
Plug this into \reff{Heston-Ito2} and note that $\hat \Th^t_t=V_t$, we obtain
\bea
\label{Heston-Ito3}
d C_t = \pa_x u \big(t, S_t,   \Th^t_{[t, T]}\big) dS_t +  (T-t)^{{1\over 2}-H}  \big\la \pa_\o u\big(t, S_t,   \Th^t_{[t, T]}\big), a^t\big\ra d\hat \Th^t_T.
\eea
That is, $C_T$ can be replicated by using $S_t$ and $\hat \Th^t_T$, with the corresponding hedging portfolios $\pa_x u$ and $(T-t)^{{1\over 2}-H} \big \la \pa_\o u , a^t\big\ra$,  respectively.

\begin{rem}
\label{rem-hedging}
We notice that, to hedge $C_T = g(S_T)$ in the rough Heston model, it is sufficient to use $S$ and $\hat \Th^\cd_T$. However, if we want to hedge $C_T = g(S_T) + \int_0^T f(t, S_t)dt$ (or even more general path dependent contingent claims, which is not covered by \cite{ER}), then we shall need $S$ and $\{\hat \Th^\cd_s\}_{0\le s\le T}$. Indeed, in this case we will have $C_t := \dbE[g(S_T) + \int_t^T f(s, S_s) ds |\cF_t] =  u(t, S_t, \Th^t_{[t,T]})$,  and, provided $u$ is smooth,
\beaa
dC_t &=& \pa_x u \big(t, S_t,   \Th^t_{[t, T]}\big) dS_t +  (T-t)^{{1\over 2}-H}  \big\la \pa_\o u\big(t, S_t,   \Th^t_{[t, T]}\big), a^t\big\ra d\hat \Th^t_T\\
&&+ \int_t^T (s-t)^{{1\over 2}-H}  \big\la \pa_\o u\big(t, S_t,   \Th^t_{[t, T]}\big), a^t\big\ra d\hat \Th^t_s ds - f(t, S_t) dt.
\eeaa
In other words, the portfolio of $\hat \Th^\cd_s$ at time $t$ is $(s-t)^{{1\over 2}-H}  \big\la \pa_\o u\big(t, S_t,   \Th^t_{[t, T]}\big), a^t\big\ra ds$.
\qed
\end{rem}
 We would like to comment further on how to replicate $\Th^t$ by using $\hat \Th^t$ in more general cases.  Mathematically, as we saw in previous sections, $\Th^t$ is intrinsically more appropriate for this framework. In fact, recall \reff{u12} and the discussion afterwards. In the general model \reff{X}, if we use  $\hat \Th^t_s:= \dbE[V_s|\cF_t]$ as our "state variable", it is not clear if one would be able to derive a sensible PPDE. However, it is clear that $\hat \Th = \Th$ when $b=0$.  For the rough Heston model \reff{Heston}, thanks to the fact that its drift $b$ is linear in $V$, $\Th^t$ and $\hat \Th^t$ are still equivalent in the following sense.  Given $\Th^t_{[t, T]}$, \reff{Heston-equivalent} is a linear convolution ODE which, by Laplace transformation,  has a unique solution $\hat \Th^t$:  denoting $\a:= H+{1\over 2}$,
 \bea
 \label{Heston-equivalent2}
  \hat \Th^t_s = \Th^t_s + {\l (s-t)^\a\over \G(\a+1)}  + \int_t^s \big[\sum_{n=1}^\infty  {(-\l)^n\over \G(n\a)} (s-r)^{n\a-{1\over 2}}\big] \big[\Th^t_r + {\l (r-t)^{\a}\over \G(\a+1)}\big] dr.
 \eea
Together with \reff{Heston-equivalent1}, we have a one-to-one mapping between the paths $\Th^t_{[t, T]}$ and $\hat\Th^t_{[t, T]}$.    In this sense, it is conceivable to write $C_t$ as a function of $\hat \Th^t$ in the rough Heston model, and we believe this is the underlying reason that a PPDE could be derived in \cite{ER}. The same arguments would work for the affine Volterra process in  Abi Jaber,   Larsson, \& Pulido \cite{ALP}, where $V$ satisfies the following convolution type of Volterral SDE: 
\bea
\label{affine}
V_t  = V_0+ \int_{0}^{t} K(t-r) \Big[ [b_0 + b_1 V_r] dr + \sqrt{a_0 + a_1 V_r} dW^2_{r}\Big].
\eea

However, we emphasize  that one cannot extend \reff{Heston-equivalent} when  the volatility $V$  satisfies the following general model with nonlinear $b$ :
\bea
\label{V2}
V_t = V_0 + \int_0^t b(t; r, V_r) dr + \int_0^t \si(t; r, V_r) dW^2_r.
\eea
In this case, as before we denote
\bea
\label{VTh}
\Th^t_s := V_0 + \int_0^t b(s; r, V_r) dr + \int_0^t \si(s; r, V_r) dW^2_r.
\eea
Then we have
\bea
\label{VTh2}
\Th^t_s = \dbE [ V_s | \cF_t] - \int_t^s \dbE [ b(s; r, V_r) |\cF_t] dr.
\eea
As we mentioned above, the forward variance   $\dbE [ V_s | \cF_t]$  can be replicated by using variance swaps. For nonlinear $b$, under technical conditions, by Carr \& Madan \cite{CM}  one may replicate $\dbE [ b(s; r, V_r) |\cF_t]$ and hence $\Th^t_s$ provided one can replicate the variance options $\dbE[ (V_s - K)^+ | \cF_t]$, or the volatility options $\dbE[ (\sqrt{V_s} - K)^+ | \cF_t]$,  for all $K$. We note again that a wide range of volatility options are available in the financial market,  at least for the S\&P 500: see the VIX options in Gatheral \cite{Gatheral}.  

To conclude this article, we point out that our framework can cover much more general models than the rough Heston model \reff{Heston}. As already mentioned, we allow for nonlinear $b$ (and $\si$) in \reff{V2} and we can still  derive the PPDE and provide a perfect hedge for $g(S_T)$, as long as the PPDE has a classical solution and $\Th$ can be replicated as we discussed above. In addition, our framework also covers the fractional Stein-Stein model, where $\sqrt{V}$ is Gaussian and is the fractional Ornstein-Uhlenbeck process; see Comte \& Renault  \cite{CR}, Chronopoulou \& Viens \cite{CV}, and Gulisashvili, Viens, \& Zhang \cite{GVZ}, and references therein. 
Besides the generality of the underlying model, we also allow for more general derivatives. On the one hand, the derivatives can be path dependent in our framework; for instance, we discussed the special case $C_T = g(S_T) + \int_0^T f(t, S_t)dt$ in Section \ref{ss-path} and Remark \ref{rem-hedging}. On the other hand we can allow for nonlinear pricing (e.g. when the borrowing and lending interest rates are different) as a solution of the BSDE \reff{BSDE}. We leave these details to the interested readers and further investigations.

\begin{rem}
\label{rem-Bergomi} In this remark we provide more details concerning the rough Bergomi model considered in Bayer,  Friz, $\&$ Gatheral \cite{BFG}. Here we shall only formally discuss the hedging issues, and  leave some technical issues in Remark \ref{rem-Bergomiintegrability} below. Let $S$ be as in \reff{Heston}, but the variance $V$ is replaced with 
\bea
\label{Bergomi}
V_t := V_0 \exp(M_t - {1\over 2} \l^2 t^{2H}),\q M_{t}  =\l \sqrt{2H}\int_{0}^{t}\left( t-r\right) ^{H-{1\over 2}}  dW^2_{r}.
\eea
We note that the variance $V$ is not in the form \reff{V2}, so the situation here is slightly different from above. However, clearly the dynamics of $X=(S, M)$ is in the form of Volterra SDE \reff{X} with two dimensional $W$ and $b=0$,
\beaa
 \si(t;r, x_1, x_2) = \left[ \ba{lll}  \sqrt{1-\rho^2} x_1 \sqrt{V_0 \exp(x_2- {1\over 2} \l^2 t^{2H})} &  \q \rho  x_1 \sqrt{V_0 \exp(x_2- {1\over 2} \l^2 t^{2H})} \\\qq\qq\qq  0 & \qq\q \l \sqrt{2H} (t-r)^{H-{1\over 2}}  \ea   \right].
\eeaa 
As in \reff{Heston},  the dynamics of $S$ is linear and thus has explicit representation: 
\bea
\label{Bergomi2}
 S_t =   S_0 \exp\big(\int_0^t \sqrt{V_s} d \tilde W_s - {1\over 2} \int_0^t V_s ds\big),\q\mbox{where}~  \tilde W_t := \sqrt{1-\rho^2} W^1_t + \rho dW^2_t.
 \eea

In this case $\Th^t_s = \dbE[M_s |\cF_t] = \l \sqrt{2H}\int_{0}^{t}\left( s-r\right) ^{H-{1\over 2}}  dW^2_{r}$ (again there is no need to introduce another component corresponding to $S$), which in particular is a martingale in this case.  The option price $C_t$ in \reff{Ct} still takes the form \reff{Heston-Cu}, while \reff{Heston-Ito} becomes
\bea
\label{Bergomi-Ito}
d C_t = \pa_x u \big(t, S_t,   \Th^t_{[t, T]}\big) dS_t +  \l\sqrt{2H} \big\la \pa_\o u\big(t, S_t,   \Th^t_{[t, T]}\big) , a^t\big\ra dW^2_t,
\eea
for the same $a^t$ as in \reff{Heston-PPDE}. Define the forward variance $\hat \Th^t_s$ as in \reff{FV}. By using the orthogonal decomposition of $M$ as in Section \ref{sect-fBM}, by \reff{Bergomi} we can easily have
\beaa
\hat\Th^t_T = V_0 \exp\Big(\Th^t_T +{1\over 2} \l^2 [(T-t)^{2H} - T^{2H}]\Big).
\eeaa
By straightforward computation, we obtain
\beaa
dW^2_t = {(T-t)^{{1\over 2}-H}\over \l\sqrt{2H}} d\Th^t_T = {(T-t)^{{1\over 2}-H}\over \l\sqrt{2H} \hat \Th^t_T} d\hat\Th^t_T.
\eeaa
Plug this into  \reff{Bergomi-Ito}, we have
\bea
\label{Bergomi-Ito2}
d C_t = \pa_x u \big(t, S_t,   \Th^t_{[t, T]}\big) dS_t +  {(T-t)^{{1\over 2}-H}\over  \hat \Th^t_T}  \big\la \pa_\o u\big(t, S_t,   \Th^t_{[t, T]}\big) ,a^t\big\ra d\hat\Th^t_T.
\eea
That is, we can replicate $C_T$  by using $S_t$ and $\hat \Th^t_T$, with the corresponding hedging portfolios $\pa_x u$ and ${(T-t)^{{1\over 2}-H}\over  \hat \Th^t_T}  \big\la \pa_\o u , ~a^t\big\ra$, respectively.
\end{rem}

\section{Appendix}
\label{sect-Appendix}
\setcounter{equation}{0}
In this section we provide some sufficient conditions for Assumption \ref{assum-X} (ii).  We first remark that, by examining our proofs carefully,  it is sufficient to assume that $X$ has the $p^*$-th moment for some finite  $p^*$ large enough, however, in that case we need to put corresponding constraints on the polynomial growth order $\k$ in  Definitions \ref{defn-C0}, \ref{defn-C12bar}, and \ref{defn-C12a}, as well as the $\k_0$ in Assumptions \ref{assum-regular} and \ref{assum-irregular}.  We also remark that, for many financial models like those we saw in the previous section, the dynamics of $S$ is typically a semimartinagle and is  linear in $S$, and thus we have a representation like \reff{Bergomi2}. Then we need much lower integrability for the $S$ part, as in the standard literature.  

The following result extends  Abi Jaber, Larsson, \& Pulido \cite[Lemma 3.1]{ALP}.
\begin{thm}
\label{thm-X}
Let $(X, W)$ be a weak solution to Volterra SDE \reff{X}.  Assume, for $\f = b, \si$, $\partial _{t}\f(t; s,\cdot )$ exists for $t\in (s,T]$, and there exists $0<H<{\frac{1}{2}}$
such that, for any $0\leq s<t\leq T$, 
\begin{equation}
|\varphi (t; s,\omega )|\leq C_{0}[1+\Vert \omega \Vert _{T}](t-s)^{H-%
{\frac{1}{2}}},\quad |\partial _{t}\varphi (t; s,\omega )|\leq C_{0}[1+\Vert
\omega \Vert _{T}](t-s)^{H-{\frac{3}{2}}}.  \label{KH2}
\end{equation}
Then Assumption \ref{assum-X} (ii) holds true.
\end{thm}

{\noindent \textbf{Proof\quad }} Fix $p\ge 2$. We first show that  it suffices to prove a priori estimates by using the standard truncation arguments. Indeed, for any $n$, denote $\t_n := \inf\{t: |X_T|\ge n\}\wedge T$, and
\beaa
X^n_t := X_{\t_n\wedge t},\q b^n(t;s,\o) := b(t; s, \o) \1_{\{\t(\o)\ge s\}},\q  \si^n(t;s,\o) := \si(t; s, \o) \1_{\{\t(\o)\ge s\}}.
\eeaa
Then $(X^n, W)$ satisfies \reff{X} with coefficients $(b^n, \si^n)$ and $(b^n, \si^n)$ satisfies \reff{KH2} with the same constants $H$ and $C_0$.  Note that $X^n$ is bounded, and we shall prove that there exists a constant $C_p$, independent of $n$, such that $\dbE[\|X^n\|_T^p] \le C_p[1+|x|^p]$.  Then by sending $n\to \infty$, we prove the theorem. 

We now assume $X$ is bounded and prove in two steps that 
\bea
\label{Xest}
\dbE[\|X\|^p] \le C_p[1+|x|^p].
\eea

{\it Step 1.} Assume $T\le \d_0$, for some small $\d_0>0$ which will be specified later. Then
\beaa
&&\dbE[|X_t-x|^p] \le C_p\dbE\Big[ \big|\int_0^t b(t;s,X_\cd) ds \big|^p +  \big|\int_0^t \si(t;s,X_\cd) dW_s \big|^p\Big]\\
&&\le C_p\dbE\Big[ \big[\int_0^t |b(t;s,X_\cd)| ds \big]^p +  \big[\int_0^t |\si(t;s,X_\cd)|^2 ds \big]^{p\over 2}\Big]\\
&&\le C_p\dbE\Big[\big[\int_0^t (t-s)^{2H-1} \|X\|_T^2 ds \big]^{p\over 2}\Big]\le C_p\dbE\Big[  \big[t^{2H} \|X\|_T^2  \big]^{p\over 2}\Big] = C_p t^{pH} \dbE[\|X\|_T^p].
\eeaa
Next, for $0\le t_1 < t_2 \le T$, denoting $\d:= t_2-t_1$. When $t_1 \le \d$, we have
\beaa
\dbE[|X_{t_2} - X_{t_1}|^p] &\le& C_p\dbE\Big[|X_{t_2} - x|^p + |X_{t_1} - x|^p\Big] \\
&\le& C_p \dbE[\|X\|_T^p] [t_2^{pH} + t_1^{pH}] \le C_p\dbE[\|X\|_T^p] \d^{pH}.
\eeaa
When $t_1 > \d$, we have
\beaa
&&X_{t_2} - X_{t_1}  = I_1 + I_2,\q \mbox{where}\\
&&I_1 :=  \int_0^{t_1-\d} \Big[  \int_{t_1}^{t_2}\pa_t b(t; s, X_\cd) dt  ds +\int_{t_1}^{t_2}\pa_t \si(t; s, X_\cd) dt dW_s\Big],\\
&&I_2:= \int_{t_1-\d}^{t_2} \Big[ b(t_2; s, X_\cd)  ds + \si(t_2; s, X_\cd)  dW_s\Big] + \int_{t_1-\d}^{t_1} \Big[ b(t_2; s, X_\cd)  ds + \si(t_2; s, X_\cd)  dW_s\Big].
\eeaa
Note that
\beaa
\dbE[|I_1|^p] &\le& C_p\dbE\Big[ \big(\int_0^{t_1-\d} \int_{t_1}^{t_2} |\pa_t b(t; s, X_\cd)| dtds\big)^p +\big(\int_0^{t_1-\d} \big( \int_{t_1}^{t_2} |\pa_t \si(t; s, X_\cd)| dt\big)^2ds\big)^{p\over 2}\Big]\\
&&\le C_p\dbE\Big[\big(\int_0^{t_1-\d} \big(\d (t_1-s)^{H-{3\over 2}} [1+\|X\|_T] \big)^2ds\big)^{p\over 2} \Big]\le C_p\dbE\big[ 1+ \|X\|_T^p\big] \d^{pH};\\
\dbE[|I_2|^p] &\le&  C_p \sum_{i=1}^2 \dbE\Big[ \big(\int_{t_1-\d}^{t_i} |b(t_i; s, X_\cd)|  ds\big)^p  +\big(\int_{t_1-\d}^{t_i} |\si(t_i; s, X_\cd)|^2  ds\big)^{p\over 2}\Big]\\
&\le& C_p\sum_{i=1}^2 \dbE\Big[\big(\int_{t_1-\d}^{t_i}\big((t_i-s)^{H-{1\over 2}} [1+\|X\|_T] \big)^2  ds\big)^{p\over 2}\Big]\\
&=& C_p\sum_{i=1}^2\dbE\big[ 1+ \|X\|_T^p\big] [t_i-(t_1-\d)]^{pH} \le C_p\dbE\big[ 1+ \|X\|_T^p\big] \d^{pH}.
\eeaa
Then
\beaa
\dbE\Big[|X_{t_2} - X_{t_1}|^p\Big] \le C_p\dbE\big[ 1+ \|X\|_T^p\big] \d^{pH}.
\eeaa
This implies that
\beaa
\dbE\Big[|X_{t_2} - X_{t_1}|^p\Big] \le  C_p\dbE\big[ 1+ \|X\|_T^p\big]  \d_0^{pH\over 2}  \d^{pH\over 2}.
\eeaa
By Lemma \ref{lem-kolmogorov}, we obtain
\beaa
\dbE[\|X\|_T^p] \le C_p\dbE\Big[|x|^p + \sup_{0\le t\le T} |X_t-x|^p\Big] \le C_p\dbE\Big[|x|^p + \big[ 1+ \|X\|_T^p\big]  \d_0^{pH\over 2}\Big].
\eeaa
By choosing $\d_0$ such that $C_p  \d_0^{pH\over 2} = {1\over 2}$, we obtain \reff{Xest}. We emphasize that $\d_0$ depends on $p, H$, but not on $x$.

{\it Step 2.} Now consider arbitrary $T$. Let $\d_0$ be as in Step 1, and denote $0=t_0 < t_1<\cds< t_m = T$ be such that ${\d_0\over 2} < t_{i+1} - t_i \le \d_0$.  Note that, for $t_i \le t\le t_{i+1}$,
\beaa
X_t = \Th^{t_i}_t + \int_{t_i}^t b(t; s, X_\cd) ds + \int_{t_i}^t \si(t; s, X_\cd) dW_s.
\eeaa
Following the same arguments as in Step 1, we obtain
\beaa
\dbE\Big[\sup_{t_i\le t\le t_{i+1}} |X_t - \Th^{t_i}_t|^p\Big] \le C_p \dbE\Big[  1+ \sup_{0\le t\le t_i} |X_t|^p\Big]. 
\eeaa
Notice further that, for $t_i \le t < t+\d \le t_{i+1}$, 
\beaa
\Th^{t_i}_{t+\d} - \Th^{t_i}_t = \int_0^{t_i} \int_t^{t+\d} \pa_t b(r; s, X_\cd) dr ds +  \int_0^{t_i} \int_t^{t+\d} \pa_t \si(r; s, X_\cd) dr dW_s.
\eeaa
Then
\beaa
&&\dbE\Big[|\Th^{t_i}_{t+\d} - \Th^{t_i}_t|^p\Big]\\
&& \le C_p\dbE\Big[   \big(\int_0^{t_i} \int_t^{t+\d} |\pa_t b(r; s, X_\cd)| dr ds\big)^p +  \big(\int_0^{t_i} \big(\int_t^{t+\d} |\pa_t \si(r; s, X_\cd)| dr\big)^2 ds\big)^{p\over 2}\Big]\\
&&\le C_p\dbE\big[1+\sup_{0\le t\le t_i} |X_t|^p\big] \Big( \int_0^{t_i}  \big(\int_t^{t+\d}  (r-s)^{H-{3\over 2}} dr\big)^2 ds\Big)^{p\over 2}
\eeaa
Note that
\beaa
&&\int_0^{t_i}  \big(\int_t^{t+\d}  (r-s)^{H-{3\over 2}} dr\big)^2 ds \le \int_0^{t_i}  \big(\int_{t_i}^{t_i+\d}  (r-s)^{H-{3\over 2}} dr\big)^2 ds\\
&& = \Big[\int_0^{t_i-\d^{2\over 3}} + \int_{t_i-\d^{2\over 3}}^{t_i}\Big]  \big(\int_{t_i}^{t_i+\d}  (r-s)^{H-{3\over 2}} dr\big)^2 ds\\
&&\le \int_0^{t_i-\d^{2\over 3}}  [ \d (\d^{2\over 3})^{H-{3\over 2}}]^2 ds + C \int_{t_i-\d^{2\over 3}}^{t_i}  (t_i-s)^{2H-1} ds \le  C\d^{4H\over 3}.
\eeaa
Then
\beaa
\dbE\Big[|\Th^{t_i}_{t+\d} - \Th^{t_i}_t|^p\Big] \le  C_p\dbE\big[1+\sup_{0\le t\le t_i} |X_t|^p\big] \d^{2pH\over 3}.
\eeaa
By Lemma \ref{lem-kolmogorov}, we have
\beaa
\dbE\Big[\sup_{t_i\le t\le t_{i+1}} | \Th^{t_i}_t|^p\Big] \le C_p \dbE\Big[  1+ \sup_{0\le t\le t_i} |X_t|^p\Big].
\eeaa
Then
\beaa
\dbE\Big[\sup_{t_i\le t\le t_{i+1}} |X_t|^p\Big] \le C_p \dbE\Big[  1+ \sup_{0\le t\le t_i} |X_t|^p\Big]. 
\eeaa
Now by induction one obtains \reff{Xest} immediately.
\qed

\begin{rem}
\label{rem-Hestonintegrability}
Consider the rough Heston model \reff{Heston}. By Theorem \ref{thm-X} it is clear that $\dbE\big[\sup_{0\le t\le T} V_t^p\big]<\infty$ for all $p\ge 1$. However, we note that   the coefficient of $S$ does not grow linearly, and thus we cannot apply Theorem \ref{thm-X} on $S$ directly. We shall instead utilize the representation formula \reff{Bergomi2}.  Note that $V\ge 0$, then
\beaa
V_t  \le C + c\int_{0}^{t}\left( t-r\right) ^{H-{1\over 2}}  \sqrt{V_r} dW^2_{r},
\eeaa
for some generic constants $C, c$. Thus
\beaa
0\le \int_0^t V_s ds  \le C + c\int_{0}^{t} \int_0^s \left( s-r\right) ^{H-{1\over 2}}  \sqrt{V_r} dW^2_{r}ds  = C + c\int_{0}^{t} \left( t-r\right) ^{H+{1\over 2}}  \sqrt{V_r} dW^2_{r}.
\eeaa
Therefore, for any $n\ge 1$,
\beaa
\dbE\Big[\big( \int_0^t V_s ds\big)^n\Big] \le C^n  + C^n \dbE\Big[ \big(\int_{0}^{t} \left( t-r\right) ^{2H+1}  V_r d{r}\big)^{n\over 2}\Big] \le C^n + C^n\Big( \dbE\Big[ \big(\int_{0}^{t} V_r d{r}\big)^n\Big]\Big)^{1\over 2}.
\eeaa
This implies that 
\beaa
\dbE\Big[\big( \int_0^t V_s ds\big)^n\Big] \le C^n,\q\mbox{and hence}\q \dbE\Big[\exp\big(p \int_0^t V_s ds\big)\Big] \le C_p<\infty.
\eeaa
Now by \reff{Bergomi2} we obtain immediately that $\dbE[S_t^p] \le C_p<\infty$. Finally, noticing again that $S$ is a standard diffusion, applying Burkholder-Davis-Gundy inequality on the first equation in \reff{Heston} we see that $\dbE[\sup_{0\le t\le T}S_t^p] \le C_p<\infty$.
\qed
\end{rem}

\begin{rem}
\label{rem-integrability}
For a rough volatility model, we may also denote the state process as $X = (\hat S, V)$ where $\hat S := \ln S$. Then, in the case of \reff{Heston},  we have
\beaa
\hat S_t   =   \hat S_0 + \int_0^t  \sqrt{V_{r}} ~\Big[\sqrt{1-\rho^2} dW^1_r + \rho dW^2_r\Big] -  {1\over 2} \int_0^t V_r dr.
\eeaa
This clearly satisfies \reff{KH2} and thus $(\hat S, V)$ satisfies Assumption \ref{assum-X} (ii). However, in this case we shall write $C_t = \hat u(t, \hat S_t, \Th^t_{[t, T]})$, where $\hat u(t,x, \th) := u(t, e^x, \th)$. If it turns out that $\hat u$ has the desired polynomial growth in $x$ (which in particular requires $\hat g(x) := g(e^x)$ has polynomial growth in $x$), then we may derive  the required results by using $(\hat S, V)$.  However, when $g$ has linear growth, $\hat g$ would grow exponentially and then the integrability in Assumption \ref{assum-X} (ii) will not be enough.
\qed
\end{rem}

\begin{rem}
\label{rem-Bergomiintegrability}
In this remark we discuss the integrability for the  rough Bergomi model in Remark \ref{rem-Bergomi}. Let $p^*$ denote the  largest moment  for $S$, as introduced by Lee \cite{Lee}:
\bea
\label{p*}
p^* := \sup\{p: \dbE[S_T^p] <\infty\}.
\eea

(i) When $\rho =0$, we have $p^*=1$. For simplicity, let's assume $\l = V_0=S_0=1$. Then
\beaa
M_{t}  =\sqrt{2H}\int_{0}^{t}\left( t-r\right) ^{H-{1\over 2}}  dW^2_{r}, \q V_t :=  e^{M_t - {1\over 2} t^{2H}},\q S_t = e^{\int_0^t \sqrt{V_s} dW^1_s - {1\over 2} \int_0^t V_s ds}.
\eeaa
In particular, $V$ and $W^1$ are independent. Clearly, 
\beaa
\dbE[S_T] = \dbE\big[ \dbE[S_T |\cF^V_T]\big] = \dbE[1] = 1.
\eeaa
However, for any $p>1$ and for some generic constant $c>0$,  denote $t_0 := {T\over 2}$, we have
\beaa
\dbE[S_T^p] &=& \dbE\Big[ \exp\big(p\int_0^T \sqrt{V_s} dW^1_s - {p\over 2} \int_0^T V_s ds\big)\Big]\\
&=&  \dbE\Big[ \dbE\big[ \exp\big(p\int_0^T \sqrt{V_s} dW^1_s - {p\over 2} \int_0^T V_s ds\big) \big|\cF^V_T\big]\Big]\\
&=&  \dbE\Big[ \exp\big( {p(p-1)\over 2}\int_0^T V_s ds \big) \Big] \ge  \dbE\Big[ \exp\big( c\int_{t_0}^T V_s ds \big) \Big]\\
&=&\dbE\Big[ \dbE\big[ \exp\big( c\int_{t_0}^T V_s ds \big) \big|\cF^{W^2}_{t_0}\big]\Big]\ge \dbE\Big[ \exp\big( c\int_{t_0}^T \dbE[  V_s |\cF^{W^2}_{t_0}] ds \big) \Big], 
\eeaa 
where the last inequality is due to Jensen's inequality. Note that
\beaa
 \dbE[  V_s |\cF^{W^2}_{t_0}]  \ge c  \dbE[  e^{M_s} |\cF^{W^2}_{t_0}]  \ge c \exp\big(\dbE[M_s |\cF^{W^2}_{t_0}] \big) = c \exp\big( c \int_0^{t_0} (s-r)^{H-{1\over 2}} dW^2_r\big).
 \eeaa
 Then, by  Jensen's inequality again,
\beaa
\int_{t_0}^T \dbE[  V_s |\cF^{W^2}_{t_0}] ds  &=& {c\over T-t_0} \int_{t_0}^T  \exp\big( c \int_0^{t_0} (s-r)^{H-{1\over 2}} dW^2_r\big) ds\\
 &\ge& c  \exp\Big( {c\over T-t_0} \int_{t_0}^T \big(\int_0^{t_0} (s-r)^{H-{1\over 2}} dW^2_r\big) ds\Big) \\
 &=&  c\exp\Big(c \int_0^{t_0} [(T-r)^{H+{1\over 2}} - (t_0-r)^{H+{1\over 2}}] dW^2_r\Big).
 \eeaa
Note that  $\int_0^{t_0} [(T-r)^{H+{1\over 2}} - (t_0-r)^{H+{1\over 2}}] dW^2_r \sim N(0, c_0)$ for some $c_0>0$. Then
\beaa
\dbE[S_T^p] \ge  c\dbE\big[ \exp( c e^{c N(0,c_0)})\big] = \infty.
\eeaa

(ii) The situation could be worse if $\rho\neq 0$. Assume for simplicity that $H = {1\over 2}$ and  $\l = V_0=S_0=1$.  Then
\beaa
V_t = \exp(W^2_t - {t\over 2}),\q S_t = \exp(\int_0^t \sqrt{V_s} [\sqrt{1-\rho^2} dW^1_s + \rho dW^2_s]- {1\over 2} \int_0^t V_s ds).
\eeaa
Thus
\beaa
\dbE[S_T] &=& \dbE\Big[\dbE[S_T |\cF^{W^2}_T]\Big] =\dbE\Big[ \exp\big(\int_0^T \sqrt{V_s}  \rho dW^2_s- {1\over 2} \int_0^T V_s ds\big) \exp\big( {1-\rho^2\over 2}\int_0^T V_s ds\big)\Big]\\
&=&\dbE\Big[ \exp\big( \rho \int_0^t \sqrt{V_s} dW^2_s- {\rho^2\over 2} \int_0^t V_s ds\big) \Big] \\
&=&\dbE\Big[ \exp\big( \rho \int_0^t e^{{1\over 2}W^2_s -{s\over 4}}dW^2_s- {\rho^2\over 2} \int_0^t e^{W^2_s-{s\over 2}} ds\big) \Big].
\eeaa
This is in the framework of Girsanov theorem, but the drift $\rho e^{{1\over 2}W^2_s -{s\over 4}}$ has exponential growth. 
While we don't have a rigorous proof, we suspect that the above integral is strictly less than $1$, and then $S$ would be a strict local martingale.
\qed
\end{rem}


\begin{thebibliography}{99}


\bibitem{ALP}
Abi Jaber, Eduardo;  Larsson, Martin; Pulido, Sergio \textit{Affine Volterra processes}, preprint, arXiv:1708.08796.

\bibitem{BC}
Baudoin, Fabrice; Coutin,  Laure  {\it Operators associated with a stochastic differential
equation driven by fractional Brownian motions}, {\sl Stochastic Processes and their Applications}, 117 (2007), 550-574.


\bibitem{BFG}
Bayer, Christian; Friz, Peter; Gatheral, Jim \textit{Pricing under rough volatility}. {\sl Quantitative Finance}, 16 (2016), no. 6, 887-904.

\bibitem{BLP}
Bennedsen, Mikkel;  Lunde, Asger;   Pakkanen, Mikko {\it Decoupling the short- and long-term behavior of stochastic volatility}, preprint, 	arXiv:1610.00332.



\bibitem{BM1} 
 Berger, Marc A.; Mizel, Victor J.  \textit{Volterra equation with It\^{o}
integrals, I}, \textsl{J Intergal Equations,}  2 (1980), no. 3 187-245.

\bibitem{BM2} 
 Berger, Marc A.; Mizel, Victor J.  \textit{Volterra equation with It\^{o}
integrals, II}, \textsl{J Intergal Equations,}  2 (1980), no. 4, 319-337.

\bibitem{CM}
Carr, Peter; Madan, Dilip {\it Towards a Theory of Volatility Trading},  {\sl Volatility, Risk Publications}, R. Jarrow, ed., 417-427. Reprinted in {\sl Option Pricing, Interest Rates, and Risk Management}, Musiella, Jouini, Cvitanic, ed., Cambridge University Press, 2001, 458-476.

\bibitem{CV}  Chronopoulou, Alexandra; Viens, Frederi G.  \textit{Stochastic volatility models with
long-memory in discrete and continuous time}. \textsl{Quantitative Finance}, 
12 (2012), no. 4, 635-649.

\bibitem{CR} Comte, Fabienne; Renault, Eric  \textit{Long Memory in
Continuous-time Stochastic Volatility Models.} \textsl{Mathematical Finance}, 
8 (1998), no. 4, 291-323.

 \bibitem{Cont}
 Cont, Rama {\it Functional It\^{o} calculus and functional Kolmogorov equations.} {\sl Stochastic integration by parts and functional It\^{o} calculus,} 115-207, Adv. Courses Math. CRM Barcelona, Birkh\"{a}user/Springer, [Cham], 2016.
 
\bibitem{CF1} Cont, Rama; Fournie, David-Antoine \textit{Change of variable
formulas for non-anticipative functionals on path space.} \textsl{J. Funct.
Anal.} 259 (2010), no. 4, 1043-1072.

\bibitem{CF2} Cont, Rama; Fournie, David-Antoine \textit{Functional It\^o
calculus and stochastic integral representation of martingales}, \textsl{%
Annals of Probability}, 41 (2013), no. 1, 109-133.

\bibitem{CR2}
Cosso Andrea; Russo, Francesco {\it Strong-viscosity solutions: semilinear parabolic PDEs and path-dependent PDEs.}  {\sl Osaka Journal of Mathematics}, to appear, arXiv:1505.02927.


\bibitem{CT}
Cuchiero, Christa;   Teichmann, Josef {\it Generalized Feller processes and Markovian lifts of stochastic Volterra processes: the affine case}, preprint, arXiv:1804.10450.


\bibitem{DU}
Decreusefond, Laurent; Ustunel, Ali Suleyman {\it Stochastic Analysis of the Fractional Brownian Motion},  \textsl{Potential Analysis}, 10 (1999), no. 2, 177-214.


\bibitem{Dupire} Dupire, Bruno \textit{Functional It\^{o} calculus.}
preprint, papers.ssrn.com, (2009).

\bibitem{EKTZ} 
 Ekren, Ibrahim; Keller, Christian; Touzi, Nizar; Zhang, Jianfeng {\it On viscosity solutions of path dependent PDEs.} {\sl Ann. Probab.} 42 (2014), no. 1, 204-236.
 
\bibitem{ETZ1} Ekren, Ibrahim; Touzi, Nizar; Zhang, Jianfeng \textit{%
Viscosity solutions of fully nonlinear parabolic path dependent PDEs: part I.%
} \textsl{Ann. Probab.} 44 (2016), no. 2, 1212-1253.

\bibitem{ETZ2} Ekren, Ibrahim; Touzi, Nizar; Zhang, Jianfeng \textit{%
Viscosity solutions of fully nonlinear parabolic path dependent PDEs: part
II.} \textsl{Ann. Probab.} 44 (2016), no. 4, 2507-2553.

\bibitem{ER0} El Euch, Omar; Rosenbaum, Mathieu \textit{The characteristic function of rough Heston models}, {\sl Mathematical Finance}, to appear, arXiv:1609.02108. 


\bibitem{ER} El Euch, Omar; Rosenbaum, Mathieu \textit{Perfect hedging in
rough Heston models}, preprint, arXiv:1703.05049.

\bibitem{FH} Fouque, Jean-Pierre; Hu, Ruimeng \textit{Optimal Portfolio under
Fractional Stochastic Environment}, preprint, arXiv:1703.06969.

\bibitem{Gatheral}
Gatheral, Jim {\it Consistent Modeling of SPX andVIX Options}, {\sl Fifth World Congress of the Bachelier Finance Society}, (2008). 


\bibitem{GJR} Gatheral, Jim;  Jaisson,  Thibault;  Rosenbaum, Mathieu \textit{ Volatility is rough}, {\sl Quantitative Finance},  18 (2018), 933-949.


\bibitem{GK}
Gatheral, Jim;  Keller-Ressel, Martin {\it Affine forward variance models}, 	preprint, arXiv:1801.06416.


\bibitem{GVZ} Gulisashvili,  Archil; Viens, Frederi;   Zhang, Xin \textit{ Small-time
asymptotics for Gaussian self-similar stochastic volatility models}.
{\sl Applied Mathematics and Optimization}, (2018), https://doi.org/10.1007/s00245-018-9497-6.


\bibitem{Keller}
Keller, Christian {\it Viscosity solutions of path-dependent integro-differential equations}, {\sl Stochastic Processes and their Applications},  126 (2016), 2665-2718. 


\bibitem{Lee}
Lee, Roger {\it The moment formula for implied volatility at extreme strikes}, {\sl Mathematical Finance}, 14 (2004), No. 3, 469-480.


\bibitem{Lukoyanov}
Lukoyanov, N. Yu. {\it On viscosity solution of functional Hamilton-Jacobi type equations for hereditary systems}, {\sl Proceedings of the Steklov Institute of Mathematics}, Suppl. 2 (2007),  190-200.

\bibitem{MoVi}  Mocioalca, Oana; Viens, Frederi  \textit{ Skorohod integration and stochastic
calculus beyond the fractional Brownian scale}. \emph{Journal of Functional
Analysis},  222 (2005), no. 2, 385-434.

\bibitem{Nbook} Nualart, David  \emph{Malliavin calculus and related topics},
2nd edition. Springer V., 2006.

\bibitem{PP} Pardoux, Etienne.; Peng, Shige {\it Adapted solution of a backward stochastic differential equation.} {\sl Syst.
Control Lett.} 14 (1990), no. 1, 55-61.

\bibitem{Peng1}
Peng, Shige {\it Probabilistic interpretation for systems of quasilinear parabolic partial differential
equations.} {\sl Stoch. Stoch. Rep.} 37 (1991), 61-74.

\bibitem{Peng}
 Peng, Shige {\it Backward stochastic differential equation, nonlinear expectation and their applications.} {\sl Proceedings of the International Congress of Mathematicians.} Volume I, 393-432, Hindustan Book Agency, New Delhi, 2010.
 
\bibitem{PS}
 Peng, Shige; Song, Yongsheng {\it G-expectation weighted Sobolev spaces, backward SDE and path dependent PDE.} {\sl  J. Math. Soc. Japan} 67 (2015), no. 4, 1725-1757.
 
\bibitem{PW} Peng, ShiGe; Wang, FaLei \textit{BSDE, path-dependent PDE and
nonlinear Feynman-Kac formula.} \textsl{Sci. China Math.} 59 (2016), no. 1,
19-36.

 \bibitem{RTZ}
Ren, Zhenjie; Touzi, Nizar; Zhang, Jianfeng {\it  Comparison of Viscosity Solutions of Fully Nonlinear Degenerate Parabolic Path-dependent PDEs}, {\sl SIAM Journal on Mathematical Analysis},  49 (2017), 4093-4116.

\bibitem{RY}
Revuz, Daniel; Yor, Marc {\sl Continuous Martingales and Brownian Motion}, Springer-Verlag, 1991.

\bibitem{Yong} Yong, Jiongmin \textit{Backward stochastic Volterra integral
equations -- a brief survey}. \textsl{Appl. Math. J. Chinese Univ. Ser. B}
28 (2013), no. 4, 383-394.

\bibitem{Zhang} Zhang, Jianfeng {\sl Backward Stochastic Differential Equations --- from linear to fully nonlinear theory}, Springer, 2017.
\end{thebibliography}
\end{document}